\documentclass[pdflatex,sn-mathphys]{sn-jnl}

\jyear{2022}%

\usepackage{bm}
\usepackage{appendix}
\usepackage{nicefrac}
\usepackage{cleveref}
\usepackage{multirow}
\usepackage{dsfont}
\usepackage{float}
\usepackage{bbm}
\usepackage{amsmath}
\usepackage{booktabs}
\usepackage{enumerate}
\usepackage{algorithm}
\usepackage{accents}
\usepackage{hyperref}
\usepackage{physics}

\usepackage{diagbox}
\theoremstyle{thmstyleone}%
\newtheorem{theorem}{Theorem}
\theoremstyle{thmstyletwo}%
\theoremstyle{thmstylethree}%
\newtheorem{definition}{Definition}%

\newtheorem{cor}{{\bf Corollary}}
\newtheorem{lemma}{{\bf Lemma}}

\newcommand{\ignore}[1]{}
\newcommand{\htheta}{\theta^*}
\newcommand{\hpsi}{\psi^*}

\newcommand{\eop}{{\hfill $\blacksquare$}}

\newcommand{\vcmnt}[1]{{\color{blue}#1}}
\newcommand{\indc}[1]{\mathbbm{1}_{\{#1\}}}

\newcommand{\og}[1]{g^*_{(#1)}}
\newcommand{\dirac}[1]{{\bm \delta}_{\{#1\}}}
\newcommand{\tp}{\tilde{p}}

\newcommand{\bin}{\texttt{Bin}}
\renewcommand{\cp}{{\cal{P}}}
\newcommand{\cs}{x} 

\newcommand{\fs}{\mbox{\fontsize{6.7}{8}\selectfont{$\mathbb S$}}}
\newcommand{\fv}{\mbox{\fontsize{6.7}{8}\selectfont{$\mathbb V$}}}

\newcommand{\tg}{\Tilde{g}}
\newcommand{\perf}[1]{\mbox{\bf Perf} (#1)}
\newcommand{\pij}{{\bm \pi}^{-j}}

\newcommand{\pit}{{\bm \pi_{t}^{-j}}}
\newcommand{\bpi}{{\bm \pi}}

\newcommand{\up}{\Upsilon}
\newcommand{\Ui}{U}
\newcommand{\Np}{\mathcal{N}_P}

 \newcommand{\Gmcs}{C_{\sigma}}
\newcommand{\Gmcm}{c^*}

\newcommand{\zt} 
{\mbox{\fontsize{6.5}{8}\selectfont $\mathbb{Z}$}}

\newcommand{\szt} 
{\mbox{\fontsize{4.8}{6.3}\selectfont $\mathbb{Z}$}}

\newcommand{\uphi}{\Phi}

\newcommand{\scalx} 
{\mbox{\fontsize{6.5}{8}\selectfont $\mathbb{X}$}}

\newcommand{\sscalx} 
{\mbox{\fontsize{4.8}{6.3}\selectfont $\mathbb{X}$}}

\newcommand{\sfs} 
{\mbox{\fontsize{4.8}{6.3}\selectfont $\mathbb{S}$}}

\newcommand{\sfv} 
{\mbox{\fontsize{4.8}{6.3}\selectfont $\mathbb{V}$}}

\newcommand{\teps}{\Tilde{\epsilon}}

\newcommand{\ttheta}{\Tilde{\theta}}
\newcommand{\tpsi}{\Tilde{\psi}}
\newcommand{\nnu}{\Bar{\nu}}
\newcommand{\nnus}{\Bar{\nu}^*}

\newcommand{\baz}{\overline{\zt}}

\newcommand{\sbaz}{\overline{\szt}}
\newcommand{\Zt}{\mathbb{Z}}
\newcommand{\ps}[1]{\psi^*_{#1}}
\newcommand{\tht}[1]{\theta^*_{#1}}
\newcommand{\bnu}{\bm \nu}

\newcommand{\mvs}{E[\bar{\xi}]}

\newcommand{\ce}{{\cal{E}}}
\newcommand{\bg}{{\bf g}}
\newcommand{\sd}{s_e}
\newcommand{\CI}{C_i}

\newcommand{\tgamma}{\Tilde{\gamma}}

\newcommand{\nv}{\mbox{\tiny $nv$}}
\newcommand{\nvdf}{{\theta}_{\nv}^*}

\newcommand{\TR}[2]{#2}  

\newcommand{\ttau}{{\tilde{\tau}}}

\begin{document}

\title[]{Stochastic vaccination game among influencers, leader and public}   


\author*{\fnm{Vartika} \sur{Singh}}\email{vsvartika@iitb.ac.in}
\author{\fnm{Veeraruna} \sur{Kavitha}} 
\email{vkavitha@iitb.ac.in}

\affil{\orgdiv{IEOR}, \orgname{IIT Bombay}, \orgaddress{
\city{Mumbai}, 
\country{India}}}

\abstract{

\footnote{The work of the first author is partially supported by the Prime minister research fellowship (PMRF), India.}Celebrities can   significantly influence the public towards any desired outcome. In a bid to tackle an infectious disease, a leader (government) exploits such influence towards motivating a fraction of public to get vaccinated, sufficient enough to ensure eradication. The leader  also aims to minimize the vaccinated fraction of public (that ensures eradication) and use minimal incentives to  motivate the influencers; it also controls vaccine-supply-rates. 
Towards this, we consider a three-layered Stackelberg game, with the leader at the top. A set of influencers at the middle layer are involved in a stochastic vaccination game driven by incentives. The public at the bottom layer is involved in an evolutionary game with respect to vaccine responses.

We prove the disease can always be eradicated once the  public is sufficiently sensitive towards the vaccination choices of the influencers -- with a minimal fraction of public  vaccinated.  This minimal fraction depends only on the disease characteristics and not on other aspects. Interestingly, there are many configurations to achieve eradication, each configuration is specified by a dynamic vaccine-supply-rate and a number -- this number represents the count of the influencers that needs to be vaccinated to achieve the desired influence.
Incentive schemes are optimal when  this
number equals all or just one; the former curbs free-riding among influencers while the latter minimizes the dependency on influencers.

}


\keywords{}



\maketitle

\section{Introduction}
An infectious disease can be eradicated once herd immunity is achieved, i.e., when a large proportion of population is immune (\cite{herd}, \cite{intervention}). Immunity can either be achieved through vaccination, or in  case of some diseases, the infection itself results in immunity. However, to minimize the mortality and morbidities,  the herd immunity is primarily achieved through vaccination of large population.  
This is possible only when public is eager towards the vaccine.

The public vaccine response is an embodiment of individual vaccination decisions. In most cases, the public is hesitant towards vaccines due to perceived side effects and lack of information (e.g., \cite{hesitancy}).  Any individual compares the risk and severity of infection with the side-effects of vaccines to make a decision.  When the risk of infection is high, then the choice of vaccination is a simple one. But once a considerable fraction of population gets vaccinated, the risk of infection becomes small and fear of infection reduces even if the disease is severe. Hence, even slight side-effects of vaccines  outweigh the possible cost of infection, and lead to a non-vaccination decision. 

These types of dynamics have been studied widely using  game-theoretic framework (\cite{bauch2004}- \cite{socialplanner}). Authors in \cite{bauch2004} show that it is impossible to completely eradicate a disease with voluntary vaccination without incentives. 
Authors in \cite{vaccination} use  evolutionary game theoretic framework, and prove that the complete eradication is not an evolutionary stable limiting state. 
It is shown in \cite{ebola} that the disease (Ebola) can be eradicated  through voluntary vaccination when the relative  cost of vaccination  is infinitesimally small as compared to the cost of infection. 
Our aim is to study if some middle men, like influencers or celebrities, can be utilised to achieve complete eradication irrespective of disease, vaccine and public characteristics; authors in \cite{influencers} observed a positive influence in public when  influencers promoted  flu vaccination on social media.


When more influencers vaccinate, majority of the public may also get inclined towards vaccination, and may actually feel uncomfortable upon missing the vaccine.  The public may also perceive a smaller risk of vaccine-side-effects.  A leader can leverage upon this influence to ensure eradication. However, the influencers also consider various factors like severity of infection, side-effects etc., before vaccinating; further they may have access to more information. Thus, the leader should first  motivate the influencers  towards vaccination, possibly using incentives,  which can subsequently result in public getting motivated.

When the vaccines are first introduced, the influencers are provided incentives to get vaccinated  for a short duration. The vaccine supplies are continued on regular basis over a substantial period of time. At any vaccine availability epoch, any individual (among public) bases it's decision on  the fraction already vaccinated, disease characteristics and the influencers status (the public dynamics are modelled by suitably extending the framework of \cite{vaccination}).
In some cases,  some of the undesirable public vaccine responses can become evolutionary stable, under which the disease may not get eradicated. 
The leader aims to design a game among influencers that ensures that any such undesirable public vaccine response is not evolutionary stable. 

In all, we consider a three-layered Stackelberg game. The leader on the top announces incentive scheme for the influencers and controls dynamic vaccine supply rates. This induces a stochastic game among influencers at the middle layer which further depends upon the anticipated evolutionary stable public vaccine response and the resulting eradication status. Under the corresponding evolutionary stable vaccine response, the dynamics  of the public at the bottom layer  settle to a certain  infected and vaccinated fraction
 of  public, which also defines the utilities of various game components.  
 
The disease can be eradicated with any desired level of certainty, when public is sufficiently sensitive towards influencers.   We also prove that the vaccinated fraction required for eradication is lower-bounded by a constant, that depends only upon the disease parameters; a higher constant when disease is more infectious. The aim of the leader is to ensure complete eradication with minimal fraction of vaccinated public while optimizing the incentive scheme.  Thus we consider two notions of optimality: a) incentive optimality -- when cost of incentives required to ensure eradication is minimum, and, b) vaccine optimality -- vaccinated fraction of public is minimum at eradication.  


We show there always exists a vaccine-optimal policy. In other words, there exists an incentive scheme and dynamic vaccine supply rate such that the disease gets eradicated with  required level of certainty and with minimal fraction of public vaccinated.
Interestingly, if the leader wants to optimize the incentives, while ensuring complete eradication, it either has to: i) vaccinate all the influencers or ii) 
vaccinate just one of them (recall \cite{vaccination} shows that eradication is not possible without influencers).
However, the existence of such incentive optimal policy depends upon disease and public characteristics.

When an incentive optimal policy exists, the leader can either choose the same  which may lead to a higher vaccinated fraction of public, or  choose 
a vaccine-optimal policy (the cost of incentives will be sub-optimal). We also observe that both objectives can be optimized simultaneously when sensitivity of public towards influencers  is either high or  low. It is more expensive  to ensure eradication when  the sensitivity is moderate.

The organisation of the paper is as follows. Section \ref{sec_prob} contains the problem statement. Sections \ref{sec_pop_dyn}, \ref{sec_inf_game} and \ref{sec_lead_game} discuss public dynamics and evolutionary game, stochastic game among influencers, and leader's optimization respectively. We also provide some numerical results in sections \ref{sec_lead_game} and \ref{sec_num_res}. The notations are summarised in table \ref{table_notations}.

\section{Problem statement}\label{sec_prob}
We consider a leader who aims to encourage the public towards vaccines to eventually eradicate the disease, and a number of influencers who  can potentially impact the public response. \textit{We assume that  the perceived cost of vaccine-side-effects reduces and the perceived cost of insecurity upon missing the vaccine increases (for any individual in the public) as more influencers vaccinate.} 
If the number of vaccinated influencers is large enough, the public may vaccinate eagerly, eradicating the disease. 
Thus the leader may  consider incentivizing the influencers with an aim to encourage them towards vaccination. The leader may also control  vaccine supply/availability (VA) rates  dynamically, based on the responses of the public and the influencers.

We model this situation using  a Stackelberg game with  three layers of agents: the leader on the top, the influencers on the middle, and the public on the bottom layer.  The aim of the leader is to ensure that the disease gets eradicated with probability greater than, say $1-\delta$. Towards this, the leader announces an incentive scheme and provides $T$ days to the influencers to make a vaccination decision from the day of the introduction of vaccines  (the  incentives may include monetary benefits, goodwill, etc.).

The incentive scheme leads to a \textit{$T$-stage stochastic game} among the influencers. On any day, any susceptible influencer can either vaccinate with some probability or procrastinate the vaccination decision to the next day. The influencers make the vaccination decision after considering various factors: the incentives, the vaccine-side-effects, and the number of other vaccinated influencers; they also estimate the \textit{eventual risk of infection based upon anticipated public dynamics} -- they are aware that the public  dynamics depend upon the eventual number of vaccinated influencers, i.e., at the end of $T$ days. 
At day one, the influencers have only partial information about the side-effects of the vaccines; every day, they receive  new random data about the side-effects from the experiences of people getting vaccinated  worldwide; this improves their estimate of cost of side-effects (common public may not have access to precise details of this information). 

The public dynamics are modelled by  evolutionary framework as in \cite{vaccination}. 

Since the outcome of the stochastic game also depends upon the incentives offered by the leader, and on the corresponding response of the public, i.e., the outcome of the public evolutionary game,  and since the leader aims to minimize the cost of incentives  while keeping the probability of eradication above $1-\delta$ (for a desired $\delta<1$), we have a three-layered \textit{Stackelberg game}.  The leader also aims to minimize the fraction of public vaccinated.

To summarise, the public dynamics on the bottom layer are modeled by an evolutionary  framework,  depending upon the actions of the agents on top two layers; influencers' behavior is captured by a stochastic game on the middle layer, which depends upon the equilibrium achieved by the bottom layer, and the incentives from the leader; and leader on the top has an (multi-objective) optimization problem based on the nested equilibria achieved by the actions of the lower layer agents.

We begin with the evolutionary framework of the  population game which  depends upon $\zt$, the realisation of the number of influencers ($\Zt$) vaccinated  by the end of $T$ days, and the dynamic VA rate $\bnu$ controlled by the leader.





\section{Population Game}\label{sec_pop_dyn}
We study the public dynamics for given $(\zt,\bnu)$ by suitably extending the evolutionary framework of \cite{vaccination} -- here we  briefly introduce the framework, while more details are in \cite{vaccination}. The extension is to include the impact of dynamic vaccine supplies and the strategies of the leader and the influencers.

Let $N(t)$ represent the total population at time $t$, and let $S(t)$, $V(t)$ and $I(t)$ respectively represent the number of susceptible, vaccinated and infected individuals. Let $\phi(t)$, $ \theta(t)$ and $\psi(t)$ represent the corresponding proportions,
\begin{equation}
    \phi(t)= \frac{S(t)}{N(t)}, \ \ \psi(t)= \frac{V(t)}{N(t)}, \ \mbox{ and } \ \theta(t)= \frac{I(t)}{N(t)}. 
\end{equation}

The public evolves over time because of a new infection, recovery, birth, vaccination or death. These events occur at random times, and we model the same using exponential distribution with respective rates; basically the process evolves like a continuous time jump process (\cite{vaccination}). The infection rate is given  by $\lambda/N(t)$ (as is usually normalized in epidemic models, e.g., \cite{armbruster2017elementary}), recovery rate by $r$, birth rate by $b$ and natural death-rate by $d$. The public gets vaccinated at vaccination availability (VA) epochs, the inter VA epochs are again exponentially distributed; authors in  \cite{vaccination} consider constant VA rate while we have dynamic VA rate depending upon $\psi(t)$ as described below. 

\noindent{\bf Dynamic VA rate:} The leader may prefer to generate controlled supply of  vaccines depending upon the response of influencers and the public. If a high number of influencers get vaccinated, then the public may get confident and  show urgency towards the vaccines; the leader might then use a higher VA rate. Similarly, if a bigger proportion of public vaccinates, the remaining population  gets encouraged, and the VA rate can be increased further. To capture these, we consider the dynamic VA rate to be $\nu_b+\nu_e\psi$. Here, $\nu_b>0$ is the basic VA rate, and $\nu_e \ge 0$ is the additional VA rate, both of which can further depend upon $\zt$. 
\textit{The VA rate strategy of the leader is given by $\bnu = (\nu_b,\nu_e)$.}

\noindent
{\bf Public vaccine-response:}
At any VA epoch $t$,   any susceptible individual accepts to vaccinate with a probability depending upon $( \theta(t),\psi(t))$.
In \cite{vaccination}, many dynamic vaccine-responses are discussed, for simplicity we consider \textit{follow-the-crowd} behaviour\footnote{ \label{footnote_vaccine_response}It is observed in \cite{vaccination} that at any ESS stable against static mutation (see definition \ref{def_ess}), any individual (among public) either  vaccinates with probability 0 or 1 at any VA epoch after reaching equilibrium; the actual response of the population only influences the transient behaviour; in appendix \ref{appen_ode}, we will show that the extended framework of this paper also has similar characteristics. Thus it is sufficient to discuss one type of response.}-- as the vaccinated proportion among the population increases, the remaining susceptible individuals are encouraged to  vaccinate with a higher probability -- we say the population is using  vaccine-response $\varpi$ if any susceptible individual at any VA epoch vaccinates with probability $\varpi(\psi)=\min\{1,\beta\psi\}$; here $\beta\ge 0$ is the sensitivity parameter.

For a given vaccine-response $\varpi$, VA rate strategy $\bnu$ and the number of vaccinated influencers $\zt$, the population represented by proportions $(\theta(t),\psi(t), \eta(t))$, with $\eta(t) := \nicefrac{N(t)}{t}$, evolves according to a continuous time jump process, whose trajectory can be approximated by the solution of the ODE \eqref{eqn_ODE} given below. Such ODEs are well known in the literature, and a justification  of using the same is provided in \cite{vaccination} along with others. We provide the same via theorem \ref{thm_ode_conv} of appendix  \ref{appen_ode} after including the impact of $(\zt,\bnu)$, and the approximating ODE is given by ($\phi = 1-\theta-\psi$):
\begin{eqnarray}\label{eqn_ODE}
\dot{\theta}&=& \frac{\theta \lambda}{ \eta\varrho}\left(1-\theta -\psi - \frac{1}{\rho}\right),  \  \ \varrho =  b  + d  +  \lambda \theta \phi+ (\nu_b  + \nu_e \psi)\phi + r \theta,\nonumber\\
\dot{\psi}&=&  \frac{1}{\eta\varrho }\Big((1-\theta -\psi)\varpi(\psi)(\nu_b+\nu_e\psi)-b\psi \Big), \mbox{ and, }  \\
   \dot{\eta} &=& \frac{b-d }{\varrho} - \eta, \mbox{ with, }\nonumber  
    \ \rho:=\frac{\lambda}{r+b}. \nonumber
\end{eqnarray}
Under any  vaccine-response  $\varpi$, the population reaches an equilibrium and the corresponding  limit proportions of the population  are derived using the attractors\footnote{To completely justify using the attractors, one needs to prove that the system satisfies the assumption {\bf A}.3 of appendix \ref{appen_ode}, so that part (ii) of theorem \ref{thm_ode_conv} is applicable.
However our focus in this paper is on the three layer game.} (locally asymptotically stable equilibria in Lyapunov sense) of the above ODE. 
 

\textit{The population is strategic and thus it is appropriate to consider only the  equilibria that arise from vaccine-responses stable against  mutations}. We begin with the definition of ESSS, the vaccine-response that is Evolutionary Stable Strategy against Static mutations (introduced as ESS-AS in \cite{vaccination}); we reproduce the definition in the following.

Define  $q^* := \varpi(\psi^*)$  
as the probability with which the individuals get vaccinated after the system reaches equilibrium $(\theta^*,\psi^*,\eta^*)$ under vaccine-response $\varpi$.
Let $\varpi_\epsilon(q) := \epsilon q +(1-\epsilon)\varpi $ represent the     combined vaccine-response  of public when $\epsilon$-fraction of them deviate (at equilibrium) from $\varpi$ to   some static policy $q$ (i.e., vaccinate with probability $q$ irrespective of other things). 
A vaccine-response is said to be stable (i.e., ESSS), when the system after reaching equilibrium can not be `invaded' by mutants using any static policy $q$ (\cite{webb}):
\begin{definition}{\bf[ESSS]}\label{def_ess}
Let
$u(q,\varpi )$ represent the anticipated utility of an agent vaccinating with probability $q$ at equilibrium corresponding to public vaccine-response $\varpi$.
A vaccine-response $\varpi$ is said to be ESSS,  
i) if $\{q^* \} = {\cal B} (\varpi)$, where the static-best response set against $\varpi$ is defined  below,
\begin{equation}\label{eqn_BR_ess}
    {\cal B} (\varpi ) := \arg \min_{q \in [0,1]}  u(q,   \varpi);
\end{equation}  and 
ii) there exists an $\bar{\epsilon}>0$ such that $\{q^* \} = {\cal B} (\varpi_\epsilon(q)  )  $, for any $\epsilon \le {\bar \epsilon}(q)$ and any $q$. 
\end{definition}By  (i) an individual person (or mutant) using any static policy $q \ne q^* = \varpi(\psi^*)$ performs strictly inferior  in terms of the anticipated utility (after the system reaches equilibrium); while by (ii) mutants perform inferior even when noticeable $\epsilon$ fraction deviates to a $q \ne q^*$.

Next, we describe the anticipated utility $ u(q, \varpi) =  u(q, \varpi;\zt, \bnu) $ that drives the evolutionary game, and which depends upon the realisation $\zt$ of the number of vaccinated influencers, and $\bnu$, the VA rate strategy of the leader. 

\noindent{\bf Anticipated utility:} If any individual in population attempts to get  itself vaccinated with probability $q$, then the anticipated utility of that particular individual equals $q\times$(expected cost of vaccine)$+(1-q)\times$(expected cost of (infection + insecurity without vaccination)). 

 The expected cost of vaccine   has two components (as in \cite{vaccination}), a fixed cost component $c_{v_1}$ and a perceived cost of vaccine-side-effects. The perceived cost reduces with $\zt$, more specifically with the fraction of vaccinated influencers $\nicefrac{\zt}{M}$, where  $M$ is the total number of influencers; it also reduces   with   $\psi$, the fraction of population already vaccinated. In all, the expected cost of vaccine  is given by  $c_{v_1}+\left(1-\nicefrac{\zt}{M}\right)\min\left\{\bar{c}_{v_2}, \nicefrac{c_{v_2}}{{\psi}}\right\}$, where    $\bar{c}_{v_2}$ represents an upper bound on perceived cost of vaccine-side-effects.

 The  cost of infection  is due to the  inconvenience that an infected individual undergoes. The expected value of it is given by  $i_p({\theta},\psi)c_{i} $ where  $c_{i}$ is the cost due to infection and $i_p({\theta},\psi)$ is the probability of getting infected before next VA epoch  --  because of exponentially distributed  events\footnote{The next VA epoch is after exponentially distributed time with parameter $(\nu_b + \nu_e \psi)$,  an individual can get infected if it comes in contact with any one of the $I = \theta N$ infected individuals within that time, and any contact is after exponentially distributed time with parameter  $\lambda/N$.},  $i_p(\theta,\psi)= \frac{\lambda \theta}{\lambda \theta + \nu_b + \nu_e \psi}$.

 An individual may feel insecure without vaccination, and this perceived loss can increase with the number of vaccinated influencers $\zt$; let    $c_f(\zt)$  (an increasing function) represent the perceived loss of not vaccinating when $\zt$ influencers are vaccinated.  Further a higher value of $c_f(\zt)$  for the same  $\zt$ implies the public is more sensitive towards influencers.  
Thus in all, we model  the anticipated utility by:

\vspace{-3mm}
{\small
\begin{eqnarray}\label{eqn_utility}
  u(q, \varpi) &:=& q\left(c_{v_1}+\left(1-\frac{\zt}{M}\right)\min\left\{\bar{c}_{v_2}, \frac{c_{v_2}}{{\psi}}\right\}\right)+(1-q)[i_p({\theta},\psi)c_{i}+c_f(\zt)],\nonumber\\
  &=& q h(\varpi)  + i_p({\theta},\psi)c_{i}+c_f(\zt), \mbox{ where }\\
  h(\varpi
  )&=& h({\theta}, {\psi}) := c_{v_1} +\left(1-\frac{\zt}{M}\right)\min\left\{\bar{c}_{v_2}, \frac{c_{v_2}}{{\psi}} \right\} - i_p({\theta},\psi)c_{i}-c_f(\zt).\nonumber
\end{eqnarray}}

\textit{The resultant of the bottom level evolutionary game is the limit proportion  of  infected and vaccinated population at equilibrium,  under  an ESSS vaccine-response (briefly referred to as ES-limit); such a vaccine-response must be stable against static mutations as in definition \ref{def_ess}, when the anticipated utilities are given by \eqref{eqn_utility}.} Our immediate quest is to identify the components of the ES-limit $(\theta^*, \psi^*, q^*)$  using the attractors of the corresponding ODE \eqref{eqn_ODE} for any given $(\zt,\bnu)$, where  $q^*= \varpi(\psi^*)$ is now the vaccine-response at an  ES-limit.

There are many attractors of the ODE \eqref{eqn_ODE} depending upon the vaccine-response $\varpi$. From the structure of the best response in \eqref{eqn_BR_ess} and the  utility function \eqref{eqn_utility} specific to stability against static mutation, it is clear that   $q^*$ component of the ES-limit can take value either 0 or 1 (see \cite[Lemma 1, Theorem 5]{vaccination} for similar details). 
 To begin with,  we summarise all such  candidate\footnote{We actually restrict to those candidates such that if $\varpi(\psi^*)=1$, then it remains 1 in some neighbourhood, i.e., $\beta\psi^* \ne 1$. Similarly, we consider  $\beta \ne \nicefrac{b\rho}{\nu_b}$ and $\nu_e \ne b\rho -\nicefrac{\nu_b}{\nvdf}$ (see table \ref{table_ess_candidate}). The focus  of the paper is on the complicated three-layered game,  and we avoid these corner cases to keep discussions simple -- it is not easy to predict the ODE behaviour for them. \label{foot_remains1}} ODE-attractors (i.e.,  with $q^* \in \{0,1\}$)  
in table \ref{table_ess_candidate};
these are the equilibrium points of ODE \eqref{eqn_ODE}, i.e., the zeros of the RHS of the ODE 
after replacing $\varpi(\psi) $ with 0 or 1,   which are also locally asymptotically stable; 
the proof of local asymptotic stability parallels that in \cite{vaccination} and is provided in  appendix \ref{appen_ode}.
\vspace{-4mm}
\begin{lemma} \label{lem_self_erad}
Consider $\beta<\frac{b}{\nu_b}$, then $(0,0, \eta^*)$ is a locally asymptotically stable attractor of the ODE \eqref{eqn_ODE} with $\eta^*=\frac{b-d}{b+d+\nu_b}$ if and only if $\rho\le 1$. \eop
\end{lemma}%
As in \cite{vaccination}, one can call $(0,0,\eta^*)$ as \textit{self-eradicating} attractor. By the above lemma, such an attractor exists only when $\rho\le 1$. Further, by  lemma \ref{lem_ess} given later,  there exists a public vaccine response such that self-eradication is evolutionary stable. \textit{Thus leader's intervention is not required for $\rho \le 1$. We hence consider $\rho>1$} and first identify the candidate ES-limits (proof in appendix~\ref{appen_ode}):
\begin{lemma}\label{lem_attractors}
Consider $\rho>1$ and any equilibrium point  $(\theta^*,\psi^*, \eta^*)$ of ODE \eqref{eqn_ODE} with $(\theta^*,\psi^*)$ components as listed in table \ref{table_ess_candidate} and $\eta^*=\frac{b-d}{\varrho^*}$ (where $\varrho^*$ equals $\varrho$ of \eqref{eqn_ODE} calculated at $(\theta^*,\psi^*)$), which further satisfies the conditions 
 of footnote \ref{foot_remains1}. Then $(\theta^*,\psi^*, \eta^*)$  is an attractor that satisfies  
 $q^*\in\{0,1\}$ if
and only if the conditions in the first column of  table \ref{table_ess_candidate} are satisfied. 
\eop
\end{lemma}
%
%
%
\renewcommand{\arraystretch}{2}
{\begin{table}[h]
\vspace{-5mm}
\centering
\setlength\tabcolsep{3pt}
\begin{minipage}{0.4\textwidth}

{\small
\begin{tabular}{|c|l|l|l|}
\hline 
Condition & $q^*$ & $(\theta^*,\psi^*)$\\
\hline
$\beta< \frac{b\rho}{\nu_b}$ & $0$ & $(\nvdf,0)$\\ \hline 
\renewcommand{\arraystretch}{1.5}
  \begin{tabular}{l}
  $\beta \ps{e} >1$, \\
    $^{\dagger}\nu_e > b \rho -\frac{\nu_b}{\nvdf}$
  \end{tabular}  & $ 1$&$\left(0,  \ps{e}\right)$  \\ \hline
  \renewcommand{\arraystretch}{1.5}
\begin{tabular}{l}
  $\beta \ps{o} >1$, \\
       $  0\le \nu_e < b \rho -\frac{\nu_b }{\nvdf}$ 
  \end{tabular}           &$1 $& $(\tht{o},\ps{o})$  \\ \hline
\end{tabular}}
\end{minipage}
\hspace{7mm}
\begin{minipage}{0.5\textwidth}

    \vspace{-7mm}    
    {\small \begin{eqnarray} \label{eqn_attractors}
    \nvdf &:=& 1-\frac{1}{\rho}, \nonumber\\
        \ps{e} &:=&  \hspace{-2mm}
        \left \{ \begin{array}{ll}
             \frac{-(b +\nu_b-\nu_e)+\sqrt{(b +\nu_b-\nu_e)^2 + 4 \nu_e \nu_b}}{2 \nu_e }  & \\
          & \hspace{-15mm} \mbox{ if } \nu_e > 0,  \\
          \frac{\nu_b}{b+\nu_b}   & \hspace{-15mm}
          \mbox{ if } \nu_e = 0 ,
        \end{array} 
        \right . \nonumber  \\
        \ps{o}& := &\frac{\nu_b}{b\rho-\nu_e}, 
    \ \  \tht{o} :=\nvdf -\ps{o} .
    \end{eqnarray}}
\end{minipage}
\vspace{2mm}
\caption{Candidate ES-limits (ones with $q^*\in\{0,1\}$):  $^{\dagger}$When $\nu_b > b \rho \nvdf$, then $\nu_e \ge 0$.   \label{table_ess_candidate} 
}
\vspace{-5mm}
\end{table}
}

Thus given the disease characteristic 
$(\lambda, r, b)$ and VA rate strategy $\bnu$, there are many attractors,  we next identify the ones that are evolutionary stable.

\renewcommand{\arraystretch}{1.2}
\begin{table}[h]
    \centering
    \begin{tabular}{|l|l|l|}
    \hline
       Costs of  &   $\CI$, $\{C_{\sd,t}\}$, $C_v$ cost of infection, time-dependent vaccine side-effects\\  
     influencers  & and vaccination. $c$ is a realization of $C_{\sd,t}$ 
     \\
       \hline
        Costs of  public &   $c_{v_1}, c_{v_2}, {\bar c}_{v_2}$ cost components related  to vaccination, $c_i$ cost of infection, \\
        &  $c_f(\cdot)$ cost of insecurity upon missing out the vaccine \\
       \hline
       $\lambda$, $r$, $b$ & infection rate, recovery rate,  birth rate,  $\rho :=\nicefrac{\lambda}{r+b}$\\
       \hline
      $ \Zt,\zt$  &   number of vaccinated influencers after $T$, $\zt$ is a realisation of $\Zt$ \\
      $\baz$ & number of vaccinated influencers that ensure complete eradication \\
\hline
       $\bnu =(\nu_b, \nu_e)$ & basic rate of vaccination, additional vaccination rate \\ 
       $\bg = (g_0, \dots, g_{\sbaz})$ & incentives offered by leader \\
  \hline
      $\nvdf$ &  maximum infected fraction (at non-vaccinating ESSS) \\ 

      $\ce$, $\ps{e}$  &  complete eradication, vaccinated fraction at eradicating ESSS\\ 
      
      $\tht{o}$, $\ps{o}$ & infected and vaccinated fraction at co-occuring ESSS\\
       \hline
    \end{tabular}
        \caption{\centering Notations  \label{table_notations}}
        \vspace{-6mm}
\end{table}

\subsection{Evolutionary stable limits }
For any attractor of table \ref{table_ess_candidate} to be an ES-limit, the corresponding   vaccine-response $\varpi$ should be an ESSS as in definition \ref{def_ess}. We begin with few definitions related to  \eqref{eqn_utility} (recall $i_p(\theta^*,\psi^*)= \frac{\lambda \theta^*}{\lambda \theta^* + \nu_b + \nu_e \psi^*}$): 

\vspace{-3mm}
{\small
\begin{eqnarray} \label{eqn_h_func}
h_i(\zt;\bnu)&:=& h\left(\nvdf,0\right)= c_{v_1} + \left(1-\frac{\zt}{M}\right)\bar{c}_{v_2} - \frac{\lambda \nvdf  c_{i}}{\lambda \nvdf+\nu_b} - c_f(\zt),\nonumber\\
h_v(\zt;\bnu)&:=& h\left(0,\ps{e}\right)=  c_{v_1} + \left(1-\frac{\zt}{M}\right)\min\left\{\bar{c}_{v_2}, \frac{c_{v_2}}{{\ps{e}}} \right\} - c_f(\zt)\mbox{ and, }\\
h_v^o(\zt;\bnu)&:=& h(\tht{o},\ps{o}) =c_{v_1} + \left(1-\frac{\zt}{M}\right)\min\left\{\bar{c}_{v_2}, \frac{c_{v_2}}{{\ps{o}}} \right\}- \frac{\lambda \tht{o} c_{i}}{\lambda \tht{o}+\nu_b+\nu_e\ps{o}}  -c_f(\zt).\nonumber
\end{eqnarray}}Note that, $\ps{e}$ and $\ps{o}$ depend upon $\bnu$, however, we omit mentioning  the explicit dependency when the context is clear. We  immediately have the following result regarding the ESSS (proof is in appendix \ref{appen_ode}).

\begin{lemma} 
\label{lem_ess}
i) {\bf Self-eradicating ESSS:} If $\rho\le 1$  then there exists an ESSS with $q^*= 0$ and ES-limit $(0, 0)$; now consider $\rho>1$:\\
ii) {\bf Non-vaccinating ESSS:} If $h_i(\zt;\bnu)>0$, then there exists an ESSS with $q^*= 0$ and ES-limit $(\nvdf, 0)$;\\
iii) {\bf Eradicating ESSS:}  If  $h_v(\zt;\bnu)<0$, $ \nu_e > b \rho -\frac{\nu_b}{\nvdf
}$ and $\nu_e \ge 0$, then there exists an ESSS  with $q^*=1$ and ES-limit $(0,\ps{e})$, thus  the disease gets eradicated;\\
iv)  {\bf Co-occuring ESSS:} If $h^o_v(\zt;\bnu)<0$ and $0\le \nu_e < b \rho -\frac{\nu_b}{\nvdf
}$, then there exists an ESSS with $(\tht{o},\ps{o})$ as ES-limit and $q^*=1$; and \\
v) if none of the above are satisfied, there is no ESSS. \eop
\end{lemma}
The above lemma summarises various possible ES-limits for any given strategy and response of the top level agents,  to be specific given $(\zt, \bnu)$. 
At the non-vaccinating ESSS, the disease is not  eradicated and the public is not willing to vaccinate at all ($q^*=0$) at the equilibrium. The  eradicating ESSS is the desired one, where the disease is eradicated and the fraction of vaccinated public reaches $\ps{e}$ of table \ref{table_ess_candidate}. The co-occurring ESSS $(\tht{o},\ps{o})$ is an intermediate limit, where the disease is eradicated partially. There is  a possibility of multiple ES-limits for some $(\zt,\bnu)$, but the eradicating and the  co-occurring ESSS can not coexist. It is also possible that there is no ESSS for some $(\zt,\bnu)$.

Observe eradicating ESSS is not possible with any $\zt$, if  $c_{v_1} - c_f(M) \ge 0$; recall $c_f(\zt)$  increases with $\zt$. In other words, \textit{if the public is not sufficiently sensitive towards the influencers, it is not possible to eradicate the disease.} Thus we assume the following: 

  \begin{enumerate}[{\bf A.}1]
    \item The maximum cost of insecurity upon missing out the vaccines is more than the fixed cost of vaccines, $c_{v_1} - c_f(M)<0$. 
\end{enumerate}


\noindent
{\bf Complete Eradication:}
The public anticipatory utility functions
$h_i,h_v,h_o$ governing the evolutionary stability of various attractors of table \ref{table_ess_candidate} depend upon
 disease characteristics $(\lambda,r,b)$, number of vaccinated influencers $\zt$, and VA strategy $\bnu$. Towards ensuring eradication, the leader's aim is two-fold: i) to ensure none of the unwanted attractors  are evolutionary stable (attractors with non-zero infected fraction of population, $(\nvdf,0)$ and $(\tht{o},\ps{o})$), and ii) at least one of the vaccine-responses leading to eradicating attractor $(0,\ps{e})$  is evolutionary stable (attractor with zero infected fraction of population). Let  $\ce$ represent such an event, briefly referred to as  \textit{complete eradication}.
 
The leader attempts to achieve complete eradication by controlling $(\zt,\bnu)$, by appropriately designing influencer's game discussed in section \ref{sec_inf_game}, VA strategy $\bnu$, and using the knowledge of $(\lambda,r,b)$.
Thus, it considers a design such that non-vaccinating and co-occurring attractors are never stable,  by having $h_i \le 0$  and $\nu_e > b \rho - \frac{\nu_b}{ \nvdf}$ respectively.
%
Further, it has to ensure that the attractor with $q^*=1$ is eradicating ESSS by having $h_v<0$.
Thus conditioned on  the outcome of the influencers game, $\zt$ (number of vaccinated influencers), 
the leader considers the following as the conditional probability of complete eradication,
\begin{equation} \label{eqn_condi_prob_erad}
P (\ce \mid \Zt =
\zt)   =  
\indc{h_i(\szt;\bnu)\le 0} \indc{h_v(\szt;\bnu)<0}  \indc{ \nu_e > b \rho -\frac{\nu_b }{\nvdf}, \ \nu_e \ge 0}. 
 \end{equation}%
By considering such a conditional probability into the optimal design, the leader ensures the emergence of an evolutionary stable vaccine response, under which at equilibrium the fraction of infected population is zero; and none of the responses whose attractors have non-zero infected fraction are stable.

Recall, the aim of the leader is to ensure complete  eradication with probability above $(1-\delta)$,  to be more explicit $P(\ce(\Zt;\bnu)) > 1-\delta$. Thus any $\bnu$ negating, 
\begin{equation}\label{eqn_eradication_condition}
   \nu_e >   b \rho -\frac{\nu_b }{\nvdf}, \mbox{ and } \nu_e\ge 0,
\end{equation}
   becomes infeasible -- as then $P(\ce(\Zt;\bnu)) = 0$ irrespective of $\Zt$ (see \eqref{eqn_condi_prob_erad}); 
 we say $\bnu$ is \textbf{admissible} if \eqref{eqn_eradication_condition} is satisfied. 
For any admissible $\bnu$, we now identify a threshold $\baz(\bnu)$ on number of vaccinated influencers such that, complete eradication is ensured whenever $\Zt \ge \baz(\bnu)$: 
\begin{theorem}\label{thm_eradication_es}
Assume  \textbf{A.}1. For any admissible $\bnu$, there exists a  $\baz(\bnu):=\min\{\zt: h_v(\zt;\bnu)<0, h_i(\zt;\bnu)\le 0\}$. Then $P(\ce\mid\ \Zt = \zt)=1$  if and only if $\zt\ge\baz(\bnu)$.
\end{theorem}%
\noindent\textbf{Proof:} 
By \textbf{A.}1, there exists a $\baz(\bnu)\le M$ since, $h_i(M;\bnu)<h_v(M;\bnu)<0$. 
By monotonicity of $h_i$, $h_v$ and from 
\eqref{eqn_condi_prob_erad}-\eqref{eqn_eradication_condition}, for any $\zt \ge \baz(\bnu)$, we have,
 $P (\ce \mid
\Zt= \zt) =1$. For the converse, say $P (\ce \mid
\Zt = \zt) =1$ for some $(\zt,\bnu)$---then $h_v(\zt;\bnu)<0$ and $h_i(\zt;\bnu)\le 0$; by monotonicity and definition, this implies  $\zt \ge \baz(\bnu)$. \eop 

By virtue of the above theorem, the purpose of the leader (while using $\bnu$) would be to design a game among influencers that ensures the number of vaccinated influencers $\Zt$, at the corresponding NE,  is at least $\baz(\bnu)$ with  probability more than $(1-\delta)$. In the next section, we  analyse the stochastic game among the influencers, after concluding this section with a summary. 

\subsubsection*{Summary of the population game}
\begin{itemize}
    \item  In lemmas \ref{lem_self_erad} and \ref{lem_attractors}, we identified some important attractors of the population game for given disease characteristics and  strategy/outcome of top layers.
    \item In lemma \ref{lem_ess}, we identify and characterize  the attractors that result from an evolutionary stable (public) vaccine response. 
    \item In theorem \ref{thm_eradication_es}, we prove that it is possible to completely eradicate the disease if the number of vaccinated influencers is more than a certain $\baz(\bnu)$, depending upon dynamic  vaccine supply $\bnu.$
\end{itemize}

\section{Game among Influencers}\label{sec_inf_game}
By theorem \ref{thm_eradication_es}, the disease gets completely eradicated if at least $\baz(\bnu)$ (briefly referred to as $\baz$) influencers get vaccinated when the vaccines are supplied to public according to $\bnu$. Thus, the outcome of the bottom level game (ES-limit) depends upon    the distribution of the actual number of vaccinated influencers  $\Zt$. We  now derive the NE  of the influencers' game for any given strategy of the leader.

The leader's strategy, apart from $\bnu$, includes a second  component $\bg=(g_0,g_1,\dots,g_{\sbaz-1}),$ the vector of incentives. Any influencer upon vaccinating within the given $T$ days, can potentially receive incentive; this incentive  on any day $t$ is given by $g_z$, if \textit{the number of vaccinated influencers  by the end of the previous day}, $Z_{t-1} = z$. 
 Further, $g_z=0$ for $z\ge \baz$ as the disease will be eventually eradicated for such $z$ by theorem \ref{thm_eradication_es}, and it is not useful for the leader to provide any more incentive. Also, $\bg$ does not depend on $t$, the day of vaccination.
 Finally observe $\Zt = Z_T $.




For given  leader strategy $(\bg,\bnu)$, the incentive structure induces a {\it T-stage stochastic game} among the influencers (a.k.a. agents), with $T-1$ decision making epochs (days); the terminal reward depends upon the expected outcome of the bottom level game. At any epoch, any susceptible agent can either vaccinate or procrastinate the vaccination decision to the next epoch. \textit{Further, we assume that the agents do not get infected during the relatively short span of initial $T$ days, hence, they can either be in vaccinated or susceptible state};  only susceptible agents can make decisions.

The vaccination decision of the influencers  depends upon their own estimated costs of vaccine and infection, and the incentives announced by the leader.
The influencers are \textit{more informed}, thus their anticipated costs are different than that of the common public. 
For example, the vaccines are getting introduced  worldwide, and not just to the area under consideration; the influencers have access to regular updates of vaccine-side-effects from those areas, which common public may not have access to.
%
%

%
%


The initial  available information (on day one) on  the cost of vaccine-side-effects  is represented by $C_{\sd,1}$.
Every day, influencers receive  new (random) data about the side-effects. Say, some $N$  individuals (outside the area under consideration)  get vaccinated  every day. Then the estimate $C_{\sd, t}$ of the cost of vaccine-side-effects on 
any day $t \ge 2$ can be improved as below: 
\begin{eqnarray}\label{eqn_ck}
C_{\sd,t} 
&=& C_{\sd,t-1} + \frac{1}{t}\left(\bar{\xi}_{t}-C_{\sd, t-1}\right) \text{ with }\ 
\bar{\xi}_{t}=\frac{1}{N}\sum_{j= (t-1)N +1}^{tN} \xi_j,
\end{eqnarray}where $\xi_j$ is the cost of side-effects experienced by an individual $j$. We make the following assumption:
\begin{enumerate}[{\bf A.}2]
    \item  $\{\bar{\xi}_i\}_i$   are non-negative  i.i.d. random variables with  mean $\mvs>0$ and  density $f_{\bar{\xi}}$ supported on $[0,\infty)$. Further  $p_0:= P(\bar{\xi}_1=0)\ge 0$ indicates fraction of people with no side-effects. In all, $\bar{\xi}_i$ for any $i$ is distributed as $p_0 \dirac{0} (\cdot) + (1-p_0)f_{\bar{\xi}} (\cdot)dx$, where $\dirac{0} $ is Dirac measure at 0.
\end{enumerate}
\textit{All the expectations in the paper are conditioned on $C_{\sd,1}$.}
%
As in population game, the vaccination cost also  consists of an additional  fixed component $C_v$ (can be different from $c_v$, see table \ref{table_notations} for notations).

There are no incentives after $T$ days, hence we assume that the influencers  that did not vaccinate in $T$ days will not vaccinate thereafter, and remain susceptible throughout\footnote{ Firstly, the influencers do not get any incentives after $T$ days. Secondly, $T$ days provide sufficiently good estimates of vaccine-side-effects, as  people vaccinating  worldwide are in  large numbers. Thus, no extra information is revealed to the influencers to compel them to change their decision after $T$ days.}.
Thus, any susceptible influencer after $T$ days
incurs a cost of infection depending upon the public response (i.e., the ES-limit), which we describe next.  The anticipated probability of infection of such influencer is zero if the disease is eradicated at ES-limit (one which ensures complete eradication),  and one else. Recall, the ES-limit depends upon $\Zt=Z_T$ and from theorem \ref{thm_eradication_es}, the disease is completely eradicated if $\Zt\ge\baz$. When  $\Zt<\baz$, any susceptible influencer has cost of infection $C_i$ which can be different from $c_i$ because of various factors, for example medical facilities.

The above discussions are instrumental in developing the influencers' stochastic game, the
components of which are described below:

\noindent{\bf Decision epochs:} The decision epochs are days $t=1, \dots, T-1$.

\noindent{\bf State space ${\cal S}$:} All the agents at any decision epoch are aware of the vaccination status of other agents, they also have the estimate of cost of vaccine-side-effects $C_{\sd, t}$, as in \eqref{eqn_ck}. Let ${\cal X}_t^{j}\in\{\fv,\fs\}$ be the vaccination status of $j$-th influencer at epoch $t$, then the state of the system is given by $({\cal X}_t^{1},\dots, {\cal X}_t^{M}, C_{\sd, t})$.
The influencers are identical, thus it is sufficient to consider $X^j_t:=({\cal X}^j_t,Z_t,C_{\sd, t})$ as the state relevant for influencer $j$, where $Z_t:= \sum_{m} \indc{{\cal X}^m_t=\fv}$ is the number of vaccinated influencers by time $t$. Thus, a typical realisation of state relevant for agent $j$ is given by  $x^j_t=(\fs,z, c)$ or $(\fv,z, c)$ depending upon the vaccination status of the agent $j$; here
$c \in  \mathbbm{R}^+$ is a realization of $C_{\sd, t}$, 
and $z$ is a realization of $Z_t$.  
Observe that the state of all susceptible agents is the same and equals $(\fs,z, c)$ for some $(z,c)$, while that of the vaccinated agents (at  same $t$) equals $(\fv, z, c)$. 
We   represent the state space as below,

\vspace{-5mm}
$${\cal S}:=\{ (\scalx,z, c) : \scalx\in\{\fv,\fs\}, z \in \{0, \cdots,
M\} \mbox{ and } c \in \mathbbm{R}^+ \}.$$ 

\noindent{\bf Action space ${\cal A}$:} Any agent can choose an action based on its state. If the agent is vaccinated, then there is no action. If the agent is susceptible, then it can choose to either vaccinate i.e., $a= \fv$, or to remain susceptible, i.e., $a=  \fs$,  in that decision epoch. Thus the action space is  ${\cal A}_{x} = \{\fv,\fs\}$ when $x=(\fs,z,c)$ and ${\cal A}_{x} = \{\fv\}$  when $x=(\fv,z,c)$, \textit{where $\fv$ (same symbol chosen for uniformity) now represents a dummy action}. Then ${\cal A}:= \cup_x {\cal A}_x$. Further, the action of  agent $j$ at time $t$ is represented as $A^j_t$. 

\noindent{\bf Strategy:} It is known that the best response against  Markov strategies is also a Markov strategy (e.g., see \vcmnt{\cite{Puterman}}), thus we restrict our attention to the same. The strategy for  agent $j$ is the decision rule that prescribes an action (possibly randomized)  for every state of the agent and for every $t$: 
\begin{equation}\label{eqn_def_str}
    \pi^j:=(d_1^j,\dots, d^j_{T-1}), \ \mbox{ where } d^j_t: {\cal S} \to {\cal P}( {\cal A}) \mbox{ for } t= 1,\dots, T-1,
\end{equation}and ${\cal P}({\cal A})$ is the set of probability measures over ${\cal A}$ (to be more precise, $x \mapsto {\cal P}({\cal A}_x)$).
When more details are required, we also refer $d^j_t$ as $(d^j_t(x))$ and let  $d_t^j(x,\fv)  $ represent the probability of a susceptible agent $j$ choosing $a=\fv$ at epoch $t$ in state $x$.
Thus for any susceptible agent in state $x$, \textit{decision rule $d_t^j$ is completely specified by $d_t^j(x,\fv) \in [0,1]$} since $d_t^j(x,\fs)  = 1 -d_t^j(x,\fv)$; the same is trivially true for vaccinated agents as $\fv$ is the dummy action.
We briefly represent the strategy profile of all agents other than $j$ as $\pij$, and $\bpi:=(\pi^j,\pij)$.

\noindent{\bf Stage-wise cost:} We first define the terminal cost. Say agent $j$ is in state  $x=(\scalx,z,c)$ at epoch $T$ with $\scalx\in\{\fs,\fv\}$; any agent can anticipate the ES-limit of the population game based on $z$. Basically, if $z \ge \baz$, the disease will be eradicated form theorem \ref{thm_eradication_es}, and hence a susceptible agent $j$ expects to never get infected and the terminal cost is zero. However if $z < \baz$, the agent $j$ considers the worst case scenario, where the disease will not be eradicated, and anticipates to get infected with probability one in future. The expected cost of infection in this case is given by $C_i$. On the other hand, if $j$ is vaccinated, the agent incurs the cost $c$, the realization of $C_{\sd,T}$, which is sufficiently good estimate of vaccine-side-effects. Hence, the terminal cost for agent $j$ in state $x=(\scalx,z,c)$  is,
\begin{eqnarray}\label{eqn_ter_rew}
r_T(x)=\left\{ \begin{array}{cc}
   \CI \indc{z< \sbaz}  &  \mbox{ if } \scalx= \fs, \\
   c  &  \mbox{ if } \scalx=\fv.
 \end{array} \right.
\end{eqnarray}

We now describe the expected cost at intermediate epochs $t=1, \dots, T-1$. If an agent is vaccinated, there is no action, and the stage-wise cost for such agents is zero. If an agent is susceptible at $t$, then upon vaccination, it incurs a cost of vaccination $C_v$ and also receives incentive $g_z$. If it does not vaccinate, the cost is zero at that  $t$. Thus, the stage-wise cost for epoch $t$ equals:
\begin{eqnarray}\label{eqn_rew_action}
  r_t (x, a)=\left\{\begin{array}{lll}
   C_v - g_z  &  \mbox{ if } a=\fv,\  x=(\fs,z,c), \\
    0 & \mbox{ if } a=\fs,\  x=(\fs,z,c),\\
    0 & \mbox{ if } a=\fv, \ x= (\fv, z,c).
\end{array}\right.
\end{eqnarray}Hence, the expected reward at $t$ of agent $j$ using  decision rule $d^j_t$ is given by,
\begin{eqnarray}\label{eqn_rew}
r_t (x, d^j_t)=  E[r_t (x, a)\mid x,d^j_t]=
   d^j_t(x,\fv)(C_v - g_z) \indc{\sscalx = \sfs}, \ \  \forall\  x =(\scalx,z,c).\hspace{3mm}
\end{eqnarray}


\noindent{\bf Controlled transitions:}
%
%
A susceptible agent $j$ remains susceptible at $t+1$ with probability $(1-p_j)$, if its state at epoch $t$ is $x=(\fs,z,c)$ and action is $d^j_t(x,\fv)=p_j$.  A vaccinated agent remains vaccinated.  Further, the state component $Z_{t+1}$ of any agent in next epoch equals,
 \begin{equation}\label{eqn_z_evolve}
Z_{t+1}=z + Y_{t+1} +  \indc{{\cal X}^j_t\ne \sfv} \indc{{\cal X}^j_{t+1} = \sfv},
 \end{equation}
 where $Y_{t+1} $ is the (random) number of susceptible opponents that got vaccinated at $t$-th epoch. The distribution of $Y_{t+1}$ can be derived using $\pij$ (to be more precise, one requires only decision rules at time $t$, $d_t^m$ for each $m \ne j$). %
 The third component $C_{\sd,t+1}$ transits according to \eqref{eqn_ck} and is independent of the actions chosen by the agents.

\vspace{1mm}
\noindent{\bf Game formulation:} Let $E^{\pij}_{x}[\cdot]$ represent the expectation under opponent strategy profile $\pij$ and  conditioned on $x=(\fs,0,c_{\sd,1})$. Then, this game can be described as $G=\langle {\cal{M}}, \Pi, \uphi \rangle$, where ${\cal{M}}=\{ 1,\dots,M\}$ is the set of agents, $\Pi=\{\pi^j\}_j$ is the set of strategies, and the utility of agent $j$ is given by:
\begin{equation}\label{eqn_util_dyn}
\uphi (\pi^j,\pij )= \sum_{t=1}^{T-1} E_{x}^{\pij}[r(X^j_t,A_t^j) ] +  E_{x}^{\pij}[r(X^j_T)].
\end{equation} 

Our aim is to derive the Nash Equilibrium (NE) of this game. \textit{We restrict our attention to the symmetric NE}, and derive the best response and hence the NE using backward induction (dynamic programming).

\noindent{\bf Best response:}   
%
%
We begin with deriving the best response of agent $j$ against opponent strategy profile $\pij$. This best response  is the minimizer of \eqref{eqn_util_dyn}, and is clearly a Markov Decision Process problem (\cite{maitra}). This problem can be solved using Dynamic Programming (DP) equations (\cite{Puterman}), which are given below for any state $x$:
\begin{eqnarray}\label{eqn_dp}
u_T(\cs;\pij) &=& u_{T}(\cs) = r_T(\cs),  \mbox{ and for }   t \le T-1,\\
 u_t(\cs;\pij)&=& \min_{a \in {\cal A}_x} \{r_t(\cs,a)+ E^{\pij}_{x,a}[u_{t+1}(X^j_{t+1};\pij)]\}. \label{eqn_dp_2_1}
 \end{eqnarray}%
 In the above $E^{\pij}_{x,a}[\cdot]$ represents the expectation conditioned on current state $x$, opponent strategy profile $\pij$, and action $a$ chosen by $j$ (\textit{throughout dependencies are mentioned only when required}).
Recall, $X^j_{t+1}=({\cal X}^j_{t+1},Z_{t+1},C_{\sd,t+1})$ is the state relevant to agent $j$ at time $t+1$; also recall $Z_{t+1}$ depends upon opponent-strategy profile $\pij$ (see \eqref{eqn_z_evolve}).

Any strategy is completely specified by $\{d^j_t(x,\fv) \}_{\{x, t\}}$ and  we often observe randomized strategies at NE. Hence without loss of generality and for uniformity, we consider optimization in \eqref{eqn_dp_2} with respect to $p= d^j_t(x,\fv)$ (probability of vaccination in case of susceptible agents) as below with $r_t(\cs, p)$ as in \eqref{eqn_rew}: 
\begin{eqnarray}\label{eqn_dp_2}
     u_t(\cs;\pij)&=& \inf_{p \in [0,1]} \left \{r_t(\cs,p)+ E_{x,p}^{\pij}[u_{t+1}(X^j_{t+1};\pij)] \right  \}.\end{eqnarray}%
We refer to $u_t (\cdot; \cdot)$ as the value function as is usually done in MDP literature.     
\textit{The set of optimizers in \eqref{eqn_dp_2} is referred to  as \underline{$t$-stage-BR} against $\pij$}. Observe here that any strategy profile constructed by choosing one decision rule for each $x$ and $t$ from the corresponding $t$-stage-BR is a \textit{best-response strategy} against $\pij$ (\cite{Puterman}).

\noindent\textbf{Value function of vaccinated agents:} 
Define $
\Gamma_t (c) := E[ C_{\sd, T} \mid C_{\sd, t} = c]$ to be the conditional expectation of cost of vaccine-side-effects at time $T$, conditioned on that at $t\le T$.
We consider 
$\Gamma_t(c)$ as the eventual cost of vaccine-side-effects
 (when $C_{\sd,t}=c$) -- 
it is assumed that the agents get a sufficiently good estimate of the side-effects of the vaccine by time $T$,  as the number of people  vaccinating outside the considered area in \eqref{eqn_ck} is  large. 

One may expect  that the value-function $u_{t}((\fv,z,c);\pij)$ of vaccinated agents equals $\Gamma_t (c)$; in fact the same could be true for susceptible agents vaccinating at $t$. We indeed prove this along with others (proof in appendix~\ref{appen_a}):
\begin{lemma}\label{lem_vac_val}
The value function  for a vaccinated agent  $j$ in state $x=(\fv,z,c)$ equals $u_t(x;\pij)=\Gamma_t(c)$, for any   $z,c, \pij$ and  $t$. For the susceptible agents, we have, 

 \vspace{-3mm}
{\small
\begin{eqnarray}
\label{eqn_sk} 
 E[u_{t+1}(X^j_{t+1};\pij) \mid   X^j_{t} =(\fs,z,c), A^j_t=\fv ] = \Gamma_t (c). 
 \end{eqnarray}}%
 Further, 
 $\Gamma_{t}(c)=E[\Gamma_{t+1}(C_{\sd,t+1})\mid C_{\sd,t}=c] = \frac{t}{T}c+ \frac{T-t}{T}\mvs .$ \eop
\end{lemma}%
We now  derive the $t$-stage-BR set of a susceptible agent $j$, when more than $\baz$ influencers are vaccinated by $t$.  By theorem  \ref{thm_eradication_es}, the disease is guaranteed to get eradicated. Further, now the leader does  not offer any incentive. Hence, one can anticipate that agent $j$ would not vaccinate irrespective of $\pij$. This is indeed true, as shown below (proof   in appendix \ref{appen_a}).
\begin{lemma}\label{lem_opt_fbar}
At any $t$ and $x=(\fs,z,c)$ with $z\ge \baz$ and $c \in \mathbbm{R}^+$, the $t$-stage-BR given in \eqref{eqn_dp_2}, for a susceptible agent $j$ against any  $\pij$ equals
 $\{0\}$. Further the value function $u_t(\cs;\pij) = 0$. \eop
 \end{lemma}%
%

Thus,  none of the remaining influencers get vaccinated  once $Z_t > \baz$ for some $t$; but this $Z_t$ is already sufficient for complete eradication. The bigger question is whether $Z_t$ (or $Z_T$) under an appropriate equilibrium touches $\baz$ and with what probability. We precisely delve into it in the next. 

\subsection{Symmetric Nash Equilibrium}  
%
We prove in the following that one can have multiple symmetric NE and  summarise all such possible   equilibria.
Towards this, we consider the best response \eqref{eqn_dp}-\eqref{eqn_dp_2} against symmetric strategies of the opponents -- any symmetric strategy profile $\bpi^{(M)} := (\pi,\dots,\pi)$  is an NE if the best response against $\bpi^{(M-1)}$  is    again $\pi$. 
We begin with a class of special  strategies.

\noindent\textbf{Construction of a special strategy:} For the purpose of constructing a strategy that leads to a symmetric NE, we consider quantities that correspond to
situation where all the opponents vaccinate with same probability. Towards this, we define special notations for this subsection. 
Let $E^p_{x,a}[\cdot]$ represent the expectation of quantities related to \textit{next state} $X_{t+1}$ of the tagged susceptible agent (say $j$), when all opponents vaccinate with probability $p$ at that particular epoch, current state is $x= (\fs, z, c)$ and $j$ chooses action $a$. To be more precise, for any given function $f(\cdot)$, $E^p_{x,a}[f(X_{t+1})]$ is the expected value conditioned on $p$, $x$ and $a$; corresponding state transition is governed by \eqref{eqn_ck} and \eqref{eqn_z_evolve}  as described next. Consider the extra number of opponents that got vaccinated before the next epoch, i.e., $Y_{t+1}$ of \eqref{eqn_z_evolve}; we briefly represent it as $Y$. Then $Y\sim \bin(M-z-1,p)$ is a Binomial random variable with parameters $(M-z-1,p)$. In all, if agent $j$ vaccinates at $t$, then  its next state is $X_{t+1}=(\fv,z+Y+1,C_{\sd,t+1})$,  else $X_{t+1}=(\fs,z+Y,C_{\sd,t+1})$.

We now construct a special strategy $\pi=(d_1,\dots,d_{T-1})$ backward recursively in the following. 
We begin with the case when $\baz<M$.
For $t=T-1$, define the following set $\cp_{T-1}(\cs)$ for every $x=(\fs,z,c)$,
\begin{eqnarray}\label{eqn_cal_P_T_minus_1}
 \cp_{T-1}(\cs)&:=&\left\{\begin{array}{lll}
  \{1\}   & \mbox{ if } z< \baz \mbox{ and }C_v + \Gamma_{T-1}(c)- g_z \le 0,  \\ 
  \{0\}   &   \mbox{ if } z< \baz \mbox{ and } C_v + \Gamma_{T-1}(c)- g_z\ge \CI, \\
  \{p\} &  \mbox{ if } z< \baz \mbox{ and } 0<C_v + \Gamma_{T-1}(c)- g_z<\CI, \\
  \{0\} & \mbox{ else,} 
  \end{array} \right. 
  \end{eqnarray}
where $p \in (0,1)$ is the solution of following equation,  

\vspace{-3mm}
{\small
\begin{eqnarray}\label{eqn_tilde_p}
C_v - g_z + \Gamma_{T-1}(c)&=&E_x^{p}[\CI \indc{Y < \sbaz-z}] = C_i F_{M-1}(\baz-1;p) 
\mbox{ with }\\
F_{l}(m;p)&:= &\sum_{k=1}^{m} \binom{l}{k} p^k (1-p)^{l-k}, \mbox{ the distribution of }  \bin(l, p)  .\nonumber
\end{eqnarray}}Further, $p$ solving \eqref{eqn_tilde_p} \textit{always exists and is unique}, since $p \mapsto F_{M-1}(\baz-1;p)$ is  strictly decreasing (as $\bin(M-1,p_1)$ stochastically dominates $\bin(M-1,p_2)$  if $p_1>p_2$).
For every state $x$, choose any element of $\cp_{T-1}(x)$ to define the decision rule of the special strategy $\pi$ at $T-1$, i.e.,
%
\begin{eqnarray}\label{eqn_pi_T_minus_1}
d_{T-1}(x,\fv) &\in& \cp_{T-1}(x) \mbox{ and also define, }  \\
v_{T-1}(\cs)&:=& \min\{C_v + \Gamma_{T-1}(c)- g_z,\CI \}\indc{z<\sbaz}.
\label{eqn_u_t_minus_1}
\end{eqnarray}

 When all the opponents vaccinate with probability $p$ that solves \eqref{eqn_tilde_p} in state $x$, we will observe in the coming paragraphs that the value function \eqref{eqn_dp_2_1} of agent $j$ at $T-1$ and in state $x$  is achieved by both actions $a \in \{\fs,\fv\}$; observe  this is the usual characteristic of mixed strategy NE and hence the choice.

We will follow a similar procedure to define the decision rules for remaining  epochs backward recursively using   $\{v_t(x)\}$, as explained in the following.
For $t=T-2,\dots, 1$, define the  set $\cp_{t}(\cs)$ for every $x$ recursively as follows.
\begin{eqnarray}\label{eqn_thm_br_1}
\cp_{t}(\cs)&:=&
\left \{ \begin{array}{ll}
     \{0, 1\}  \cup  \cp_{(t)} & \mbox{ if }  z< \baz \mbox{ and }  C_v + \Gamma_{t}(c)- g_z\le 0, \\
    \{0\}  \cup  \cp_{(t)} &  \mbox{ if }z< \baz \mbox{ and }  C_v + \Gamma_{t}(c)- g_z > 0,\\
 \{0\} &    \mbox{ else, }
\end{array} \right .\\
\mbox{with } \cp_{(t)}&:=&\left \{p \in (0,1):C_v + \Gamma_{t}(c)- g_z= E^{p}_{x}[v_{t+1}(\tilde{X})] \right \}, \nonumber
\end{eqnarray}where $\tilde{X}$ is the next state when current state is $x$, opponents vaccinate with probability $p$ and agent $j$ does not vaccinate; basically $\tilde{X}=(\fs,z+Y,C_{\sd,t+1})$, $Y \sim \bin(M-z-1,p)$, and $C_{\sd,t+1}$  as in \eqref{eqn_ck}.
Again  choose any element of $\cp_t(x)$ to define the decision rule $d_t$ for stage $t$ and state $x$,
\begin{eqnarray}\label{eqn_vk}
d_{t}(x,\fv) &\in& \cp_{t}(x) \mbox{ and define, }  \nonumber\\
v_{t} (\cs)&:=& 
 \left\{\begin{array}{ll}
  E_x^{0}[v_{k+1}(\fs,z,C_{\sd,t+1})],    &  \mbox{if choosen } d_t(\cs,\fv) = 0,\\
    C_v + \Gamma_{t}(c)- g_z  &  \mbox{else,}
 \end{array}\right.
\end{eqnarray}where  $(\fs,z,C_{\sd,t+1})$ is the next state that results when state at time $t$, $x=(\fs,z,c)$,  and none of the agents (including $j$) vaccinate at $t$.

If $\baz=M$, the special strategy is constructed in the same manner, with set $\cp_t(x)$ for $t=1,\dots,T-1$ now defined as below (with $v_{t+1}(x)$ as in \eqref{eqn_u_t_minus_1},\eqref{eqn_vk}),
\begin{eqnarray}
\label{eqn_cal_P_T_minus_1_for_M}
    \cp_{T-1}(\cs)&:=&\left\{
  \begin{array}{lll}
  \{1\}   & \mbox{ if } z< \baz \mbox{ and }C_v + \Gamma_{T-1}(c)- g_z <\CI,   \\
  \{0\}   &   \mbox{ if } z< \baz \mbox{ and } C_v + \Gamma_{T-1}(c)- g_z>\CI,  \\
  {[0,1]} &  \mbox{ if } z< \baz \mbox{ and } C_v + \Gamma_{T-1}(c)- g_z = \CI, \\
  \{0\} & \mbox{ else, i.e., if } z \ge \baz,
\end{array}\right.\\
\mbox{and, }\cp_{t}(\cs)&:=&
\left \{ \begin{array}{ll}
    \{0,1\} \cup \{p\} & \mbox{if }  z< \baz \mbox{ and }   C_v + \Gamma_{t}(c)- g_z\le E^1_x[v_{t+1}(\bar{X})],\\
     \{0\} & \mbox{else. } 
\end{array} \right . \hspace{5mm}
\label{eqn_thm_br_1_zM}
\end{eqnarray}where  $\bar{X}:=(\fs,M-1,C_{\sd,t+1})$ is the next state, when current state is $x$, opponents vaccinate with probability 1, and agent $j$ does not vaccinate. Further, $p\in (0,1)$ is the unique solution of, 
\begin{equation}\label{eqn_p_for_bark_M}
C_v + \Gamma_{t}(c)- g_z= E^{p}_{x}[v_{t+1}(\tilde{X})] \ \ \mbox{ with } \tilde{X}:=(\fs,z+Y, C_{\sd,t+1}).
\end{equation}%
%
%
We summarise the procedure to construct this strategy in Algorithm \ref{NE_policy_algorithm}.
\begin{algorithm}\caption{Construction of a Special strategy $\pi=(d_1,\dots,d_{T-1})$} \label{NE_policy_algorithm}
\begin{enumerate}
\item Initialize $t=T-1$.
    \item For every $x$, choose $d_{T-1}(x,\fv) \in \cp_{T-1}(x)$ of \eqref{eqn_cal_P_T_minus_1} if $\baz<M$ or \eqref{eqn_cal_P_T_minus_1_for_M} else.
    \item Define $v_{T-1}(x)$ as in \eqref{eqn_u_t_minus_1}.
    \item Set $t=t-1$.
      \item For every $x$, choose $d_{t}(x,\fv) \in \cp_{t}(x)$ of \eqref{eqn_thm_br_1} if $\baz<M$ or \eqref{eqn_thm_br_1_zM} else.
    \item Define $v_{t}(x)$ as in \eqref{eqn_vk}.
    \item If $t=1$, stop, else go to step 4.
\end{enumerate}
\end{algorithm}%

Now, we have the following theorem that characterises all  the symmetric NE
(proof in appendix \ref{appen_a}, note $g_0$ is the incentive when zero influencers were vaccinated by previous day).
\begin{theorem}[{\bf Symmetric Nash Equilibrium}]\label{thm_eqlbm}
A strategy profile $\bpi^*:=(\pi^*, \dots, \pi^*)$ is a symmetric Nash Equilibrium  if and only if $\pi^*$ is constructed as in Algorithm \ref{NE_policy_algorithm}. Further, the equilibrium utility  of agent $j$, $\uphi(\pi^*,\bpi^{-j*}) = u_1(x;\bpi^{-j*})=v_{1} (x)$ and is upper-bounded by $ C_v + \Gamma_{1}(c_{\sd,1})- g_0$ when $x=(\fs,0,c_{\sd,1})$. \eop


%
\end{theorem}

\noindent\textbf{Remarks:} i) From \eqref{eqn_thm_br_1}-\eqref{eqn_thm_br_1_zM}, the set $\cp_t(x)$ has multiple choices, so  \textit{we have multiple symmetric NE} by the above theorem.

\noindent ii) Let $\pi^*_t:=(d^*_{t},\dots,d^*_{T-1})$ be the part of $\pi^*$ from stage $t$ onwards. Then $(\pi^*_t, \dots, \pi^*_t)$ is a symmetric Nash Equilibrium for $t$-sub-game,\footnote{{\bf $t$-sub-game: }
For each realization $(z,c)$ of $(Z_t, C_{\sd,t})$, the $t$-sub-game $G_t(z,c)$ is defined to be game among the influencers from epoch $t$ onwards. Thus, $G_t(z,c)$ is a $(T-t)$-stage game with strategies as in \eqref{eqn_def_str} but from epoch $t$ onwards, e.g., $
\pi^j_t := (d^j_t,\dots, d^j_{T-1})
$. With $\pit$  representing the strategy profile of opponents, and with $x^j$ representing the starting state, i.e., at stage $t$, 
 the utility of game $G_t(z,c)$ for agent $j$  is given by,
%
\begin{equation*}
\uphi_t (\pi^j_t,\pit ; x^j )= \sum_{k=t}^{T-1} E_{x^j}^{\pit}[r(X^j_k,A_k^j) ] +  E_{x^j}^{\pit}[r(X^j_T)]. \mbox{ for any } \pi^j_t \mbox{ and } \pit.
\end{equation*}
We would require that $\sum_j \indc{x^j =(\sfv,z,c)}=z$ and by symmetry, the game is the same for any $\{x^j\}$ satisfying this.} i.e., sub-game starting at stage $t$. 

\noindent\textbf{Outcome of game among influencers:} By theorem \ref{thm_eqlbm}, there are multiple symmetric NE. To understand the preferred  outcome of the game, we compare the utility of any agent/influencer  at various symmetric NE (utility is the same for all agents by symmetry at any fixed NE). \textit{We say an NE is preferred by the influencers if the utility at that NE is the minimum}. In the following corollary, we identify one such NE. 

Define  \textit{wait-and-watch strategy} $\pi^*_w:=(0,\dots,0,d_{T-1})$ with $d_{T-1}(x,\fv)\in \cp_{T-1}(x)$ as in \eqref{eqn_cal_P_T_minus_1}, \eqref{eqn_cal_P_T_minus_1_for_M}. From \eqref{eqn_thm_br_1}, \eqref{eqn_thm_br_1_zM},  $\cp_t(x)$  with $t < T-1$ always contains zero. Thus,  wait-and-watch strategy is a special strategy, and a part of symmetric NE by theorem \ref{thm_eqlbm}. We now have (proof in appendix \ref{appen_a}),
\begin{cor}[{\bf Preferred equilibrium of influencers}]\label{cor_wait_n_watch}
Among all the symmetric NE, wait-and-watch equilibrium,  $\bpi_w:=(\pi^*_w, \dots, \pi^*_w)$ has the minimum agent-wise utility, i.e., $\uphi(\pi^*_w,\bpi_w^{-j*}) \le \uphi(\pi, \pij)$ for any symmetric NE $\bpi \ne \bpi_{w}$. \eop
\end{cor}%

By definition,  $Z_{T-1}=0$ a.s. under wait-and-watch equilibrium.
Thus as the name suggests, all the influencers wait till the last epoch before making a vaccination decision, possibly in order to get the best  possible estimate of cost of  vaccine-side-effects. 
By above corollary,  \textit{this is the preferred equilibrium for influencers. Thus leader designs influencers game by anticipating wait-and-watch NE as the outcome.} 

When $\baz<M$, the set $\cp_{T-1}(x)$ in \eqref{eqn_cal_P_T_minus_1} is singleton for every $x$, and hence the wait-and-watch equilibrium is unique for any choice of $\bg$.

When
 $\baz=M$, from \eqref{eqn_cal_P_T_minus_1_for_M} the  set $\cp_{T-1}(x)$ is singleton for all $x=(\fs,0,c)$, except for the case when $C_v + \Gamma_{T-1}(c)-g_0 = \CI $. However, from \textbf{A.}2, $\Gamma_{T-1}(C_{\sd,T-1})$ has positive mass only at $\frac{c_{\sd,1} + \mvs}{T}$ and hence $P(C_v + \Gamma_{T-1}(C_{\sd,T-1})-g_0 = \CI)=0$ for all positive $g_0$ when $C_i>C_v + \frac{c_{\sd,1} + \mvs}{T}$. Otherwise, at $g_0 = C_v-C_i +  \frac{c_{\sd,1} + \mvs}{T}$ the above probability can be positive. For this  sole corner case, we again consider the worst case scenario, and  choose the wait-and-watch NE with $d_{T-1}(x,\fv)=0$ (for such $x,\bg$) as the outcome.

 In all, we summarise the outcome of stochastic game among influencers for any $(\bg,\bnu)$ as follows: a) none of the influencers vaccinate till $T-1$; b) at epoch $T-1$,  given $C_{\sd,T-1} = c$, all the  influencers vaccinate with probability, 
\begin{eqnarray}\label{eqn_outcome_vaccination}
p(g,c) &=& \indc{C_v - g + \Gamma_{T-1}(c)<\CI} \mbox{ when } \baz(\bnu)=M \mbox{ else for  } \baz(\bnu)<M , \\
    p(g,c) &=& \left\{ \begin{array}{ll}
 1   & \mbox{ if }\  \Gamma_{T-1}(c)-g \le - C_v,  \nonumber\\ 
  0   & \mbox{ if }  \   \Gamma_{T-1}(c)-g 
  \ge \CI - C_v, \\
  \tp & \mbox{ else}, 
    \end{array}  
    \right.\hspace{6mm} \\
   && \hspace{-2cm} \mbox{ where } \tp \mbox{
    is unique solution of
    }  F_{M-1}(\baz(\bnu)-1;\tp)= \frac{C_v + \Gamma_{T-1}(c)-g}{\CI}, \hspace{6mm} \label{eqn_p_outcome_F}
\end{eqnarray}%
where \textit{we represent $g_0$ briefly as $g$ henceforth}. Hence, for any given $(\bg,\bnu)$ (more precisely for any given $(g,\bnu)$), conditioned on $C_{\sd,T-1} = c$, the number of vaccinated influencers after $T$ days are  $\Zt= Z_T \sim \bin(M,p(g,c))$ for any $c$.

We conclude this section by summarising few interesting properties of NE vaccination probability function in \eqref{eqn_outcome_vaccination}. These are instrumental in deriving an optimal strategy for the leader on the top level (proof in appendix \ref{appen_a}). 
\begin{lemma}[{\bf NE vaccination probability}]\label{lem_monotonicity_of_p}
 For any $x =(\fs,0,c)$, the  function $p(g,c)$ in \eqref{eqn_outcome_vaccination} is  non-decreasing in $g$, and non-increasing in $c$. Further  if $\baz<M$, $p$ is a continuous function, and
 \begin{enumerate}[(i)]
     \item $g \mapsto p(g,c)$ is strictly increasing on $ \{p(g,c) \not \in \{0,1\}\}$.
     \item $c\mapsto p(g,c)$ is strictly decreasing on $ \{p(g,c) \not \in \{0,1\}\}$.  \eop
 \end{enumerate}
\end{lemma}

\subsubsection*{Summary of game among  influencers}
\begin{itemize}
    \item In theorem \ref{thm_eqlbm}, we identify all  the symmetric NE of influencers' game for any given leader strategy $(\bg,\bnu)$.
    \item We also identify the NE preferred by the influencers among the above in corollary~\ref{cor_wait_n_watch} and declare it as the outcome.
\end{itemize}

\section{Leader's Optimization} \label{sec_lead_game}

 The leader   aims to ensure complete eradication $\ce$ with required level of certainty. Towards this, it provides incentives $\bg$ to the influencers and also chooses VA rate strategy  $\bnu = (\nu_b,      \nu_{e} )$. Based on  $\bnu$ and $\baz(\bnu)$ (of theorem \ref{thm_eradication_es}),  
 the leader anticipates the ESSS for the bottom level. 
 It aims to bring the corresponding probability, $P(\ce) > (1-\delta)$, for some small $\delta>0$, \textit{while   optimizing the vaccinated fraction of public and the incentive scheme}.

Say the leader chooses an admissible $\bnu$, then by theorem \ref{thm_eradication_es},
$P (\ce) = P(\Zt \ge \baz(\bnu))$,  where $\Zt = Z_T$
is the   eventual number of vaccinated influencers.   Thus for all the vaccine strategies $\bnu$ that lead to the same value of $\baz(\bnu)$ (call it just $\baz$), $P (\ce)$ does not change.
Hence we begin with an optimal incentive scheme for any given $\baz$ (which implies for a range of $\bnu$), while the joint optimization problem (involving both vaccinated fraction and incentive scheme)  is  considered later in subsection \ref{sec_joint_opt}. Towards this, consider the following constrained optimization problem, with $\Np(\bg,\bnu) = \Np(\bg,\baz):=P (Z_T<\baz)$ representing the expected non-eradication probability under the given VA rate strategy $\bnu$:
\begin{eqnarray}\label{eqn_leader_opt}
\inf_{\bg}&& E\left[ \sum_{t=1}^{T-1} (Z_t - Z_{t-1}) \  g_{Z_{t-1}}   \right] \mbox{ subject to } \Np(\bg,\baz)  \le \delta.
\end{eqnarray}%
\textit{Thus the optimization of  incentive scheme  $\bg$  as in \eqref{eqn_leader_opt} depends only  on this number $\baz$  and not on   other details of VA strategy $\bnu$}.

In the previous section, we discussed all the possible symmetric equilibrium of the stochastic game among influencers. In corollary \ref{cor_wait_n_watch}, we further show that the ``wait-and-watch" is the most preferred NE; thus the leader anticipates it as the outcome.  At this NE, the total number of vaccinated influencers $Z_T$ depends only upon  the NE vaccination probability function (used at $T-1$) given in \eqref{eqn_outcome_vaccination}. This probability in turn depends upon the \textit{incentive $g$ with zero vaccinated influencers (recall $g$ represents $g_0$)} and the realisation of cost of vaccine-side-effects $C_{\sd,T-1}$ along with other parameters.

 
Under wait-and-watch, $Z_t=0$ a.s. for $t\le T-1$. Hence the objective function in \eqref{eqn_leader_opt} simplifies to $E[Z_T g ]$, where $Z_T \sim \bin(M,l)$ conditioned on $p(g,C_{\sd,T-1})=l$ (see \eqref{eqn_outcome_vaccination}).  
Thus the cost of the leader  due to incentives in \eqref{eqn_leader_opt}, and the non-eradication probability equals (depends only on $g, \baz$):
\begin{eqnarray}\label{eqn_lead_util_cons}
       \Ui(\bg,\bnu) &= \Ui(g,\baz) &= E[Z_T g ] =  M g E[p(g,C_{\sd,T-1}) ], \mbox{ and, }\\
     \Np(\bg,\bnu)  &= \Np(g,\baz) &= E[F_M(\baz-1;p(g,C_{\sd,T-1}))]. \nonumber
 \end{eqnarray}
In the above, the expectations are with respect to $C_{\sd,T-1}$ given in \eqref{eqn_ck}. Thus, the leader's optimization problem \eqref{eqn_leader_opt} simplifies to the following:
\begin{eqnarray}\label{eqn_leader_opt_wait}
&\inf_{g} \Ui(g,\baz) \mbox{ subject to } \Np(g,\baz) \le \delta.
\end{eqnarray}
If the   non-eradication probability without any incentive (i.e., with $g=0$) is  itself within the required limit, i.e., if $\Np(0,\baz)\le\delta$, then the optimal incentive scheme for the leader is $g^*=0$.  We  now consider the other case. 

Lemma \ref{lem_monotonicity_of_p} indicates that the leader's objective function $\Ui$ is  increasing in $g$. One can further show that the constraint $\Np$ is also monotone in $g$, and that a unique $g$ satisfying the constraint with equality is the optimizer, except for a corner case. We formally prove this  in the following (proof in appendix \ref{appendix_leaders_game}):


 \begin{theorem} \label{thm_lead_optimal}
Assume 
 {\bf A.}2,  consider $\delta\in(0,1)$ and and say $\Np(0,\baz)>\delta$. Then the optimizer  $g^*$ is the unique solution of $
     \Np(g,\baz) = \delta $.
 \end{theorem}%

Above theorem further simplifies the  problem \eqref{eqn_leader_opt} considered at wait-and-watch equilibrium, and is instrumental in further analysis. The optimizer of \eqref{eqn_leader_opt} for a given $\baz$ (referred as  optimal incentive scheme) for the leader is given by $\bg^*=(g^*,g_{1},\dots,g_{\sbaz})$ with $g^*$ as in theorem \ref{thm_lead_optimal}, and an arbitrary $g_{i} \in \mathbb{R}^+$ for each $i \ge 1$.

\noindent\textbf{Comparison across various $\baz$:}
  Using theorem \ref{thm_lead_optimal}, we derive one optimal incentive scheme for each $\baz$,  and compare these  schemes before moving on to 
the joint optimization    with respect to $(\bg, \bnu)$ in subsection \ref{sec_joint_opt}. 
Towards this, define ${\Ui}_k^*$ and $\og{k}$ to be optimal objective value and optimizer of \eqref{eqn_leader_opt_wait}
respectively when $\baz = k$ for any $k \le M$; for clarity, \textit{we would like to reaffirm that $\og{k}$ is optimal $g_0^*$ when $\baz=k$}. 

We \textit{begin with the case when perfect information on side-effects is available }
%
to the influencers from the start.  In other words, say
$C_{\sd,t}=c $ a.s., for any $t$, then  $\Gamma_{t} (c) = E[ C_{\sd, T} \mid C_{\sd, t} = c] $ of lemma \ref{lem_vac_val} equals $c$. We denote such a system as $\perf{c}$ -- one can derive it by setting  $\Bar{\xi}_i=c$ a.s. for all $i$ in \eqref{eqn_ck}.
Thus {\bf A.}2 is not satisfied; however one can directly derive the optimizers for any $\baz<M$ and $\epsilon$-optimizer for $\baz=M$. 
Once again the optimizer is $\og{\sbaz}=0$  when $\Np(0,\baz)\le\delta$, and otherwise we have the following with details in appendix~\ref{appendix_leaders_game}:

\underline{Optimal cost and optimizers for $\perf{c}$}:  {\it For any $\baz<M$,
\begin{eqnarray}\label{eqn_prf_info}
\Ui^*_{\sbaz} = M \og{\sbaz} p^*_{\sbaz}   
    \mbox{ and }  \og{\sbaz}= C_v + \Gamma - C_i F_{M-1}(\baz-1;p^*_{\sbaz}), \mbox{ with } \Gamma := c,\\
     \mbox{ where } p^*_{\sbaz} \mbox{ solves } F_M(\baz-1;p^*_{\sbaz}) = \delta. \nonumber
\end{eqnarray}
Fix $\epsilon>0$. Then the optimizer and the $\epsilon$-optimizer respectively for $\baz = 1$ and $M$, and the corresponding optimal costs are given by:  
\begin{eqnarray}
 \begin{array}{ll}
 \Ui^*_1 = M (C_v + \Gamma - \CI\delta^{\frac{M-1}{M}}) (1- \delta^{\frac{1}{M}}), &\hspace{5mm}  \og{1} = C_v + \Gamma - \CI \delta^{\frac{M-1}{M}}, \\
   \Ui^*_M = M(C_v +\Gamma - \CI),  &\hspace{5mm}   \og{M,\epsilon} = C_v +\Gamma -\CI + \epsilon.  
\end{array}  \label{eqn_perf_opt}
\end{eqnarray}}

Thus we have closed-form expressions for the corner cases with $\baz = 1$ or $M$ and further one of these values of $\baz$ are the best towards the cost of incentives as proved below (proof in appendix \ref{appendix_leaders_game}):
\begin{theorem} \label{thm_perf_info_var_z}
   In case of perfect information:  i) there exists a ${\bar \delta}_1>0$ such that for all $0<\delta \le {\bar \delta}_1$, 
$\Ui^*_1 < \Ui^*_k$ for all $k < M$;
ii) there exists another ${\bar \delta}_2 $ depending upon $M$ such that
$\Ui^*_M < \Ui^*_1$ for all $0<\delta \le {\bar \delta}_2$. \eop
\end{theorem}

\noindent{\bf Remarks:}  Thus for the case with perfect information, for any $\delta \le \Bar{\delta}_1$, one can compare $U_{\sbaz}^*$  for various $\baz$ that achieve the same level $(1-\delta)$ of non-eradication probability. \textit{The smallest incentive cost, $\min_{\sbaz \le M} U^*_{\sbaz}$  is  at one of the extreme points, i.e., either  at   $\baz=1$  or  at $\baz=M$ and not at any intermediate  $\baz$}. To be more precise, when $\delta \le \Bar{\delta}_2$ (defined in theorem \ref{thm_perf_info_var_z})   $\baz=M$ is optimal; while for larger $\delta \in (\Bar{\delta}_2,\Bar{\delta}_1]$ (if there is such a range), either $\baz=1$ or $\baz=M$ is optimal.  


\textit{We now consider the general case, where $\sigma^2$, the variance of $\{\bar{\xi}\}_i$ in \eqref{eqn_ck} is non-zero}, and  assume  {\bf A.}2. Recall, our aim is to   compare across various $\baz$, which depends upon $\sigma^2$ as well as $\delta$. Hence, we explicitly mention the   dependence on $\sigma^2$ and $\delta$ in this section.
The obvious range of $\delta$ important for this study is small values, and we consider the same in further analysis (proof in appendix~\ref{appendix_leaders_game}):
\begin{theorem} \label{thm_gM_smallest}
i) For any variance $\sigma^2 > 0$,  there exists an $\alpha>0$ and $\Bar{\delta}>0$ such that,%
\begin{eqnarray*}
    \og{M } (\sigma^2,\delta) + \alpha C_i & \le&  \og{k} (\sigma^2,\delta) \mbox { and }\\
    \Ui^*_M (\sigma^2,\delta)&<&\Ui^*_k(\sigma^2,\delta)   
    \mbox{ for all } k<M \mbox{ and } \delta < \Bar{\delta}.
\end{eqnarray*}%
\label{thm_cm_c1}
ii)There exists a $\Bar{\delta}>0$ such that  the following is true for all $\delta \le \Bar{\delta}$:  for any $\epsilon>0$, there exists a $\Bar{\sigma}^2>0$ (depending upon $\delta$ and $\epsilon$) such that 
\begin{eqnarray*}
\mid \og{\sbaz} (\sigma^2,\delta) - \og{\sbaz} \mid & \le & \epsilon \ \ \ \ \ \ \  \mbox{ and }  
\mid U^{*}_{\sbaz} (\sigma^2,\delta) - U^{*}_{\sbaz} \mid  \ \le  \ \epsilon 
\mbox{ for all  $\baz < M$,  while, }    \\
\og{M}(\sigma^2,\delta) &<& \og{M,\epsilon}  \ \ \mbox{and, } 
  U_M^*(\sigma^2,\delta)< U_M^{*} + M \epsilon \  \mbox{  for all $\sigma^2 \le \Bar{\sigma}^2$ },
\end{eqnarray*}where  $\og{M,\epsilon} $,  $\og{\sbaz}$,  $U^{*}_{\sbaz} $ and $U_M^{*}$ are as defined  in \eqref{eqn_prf_info}-\eqref{eqn_perf_opt} for $\perf{c_{\infty}}$ where $c_\infty:= E[\Gamma_{T-1}(C_{\sd,T-1})] = \nicefrac{(c_{\sd,1} + (T-1)\mvs)}{T}$.
 \eop
\end{theorem}

Thus the above theorem provides required comparison across various $\baz$ for small values of $\delta$. As in perfect information case, for all $\delta$ less than a threshold, $\baz = M$ achieves minimal $U^*_{\sbaz}$.  Furthermore, part (ii)  also shows that the optimal $\og{\sbaz}$ is close to the corresponding one in theorem \ref{thm_perf_info_var_z} for small  variance.


\subsubsection*{Numerical results}
We now reinforce the theory using numerical examples, that also illustrate new insights.  In all the examples, we consider $\{\bar{\xi}_i\}_i$ to be truncated Normal random variables, i.e., $\bar{\xi}_i=0$, if the corresponding random variable is negative.

\begin{figure}[h]
\vspace{-5mm}
    \centering
\begin{minipage}{0.48\textwidth}
\includegraphics[trim={4cm 8cm 3cm 7cm}, scale =0.4]{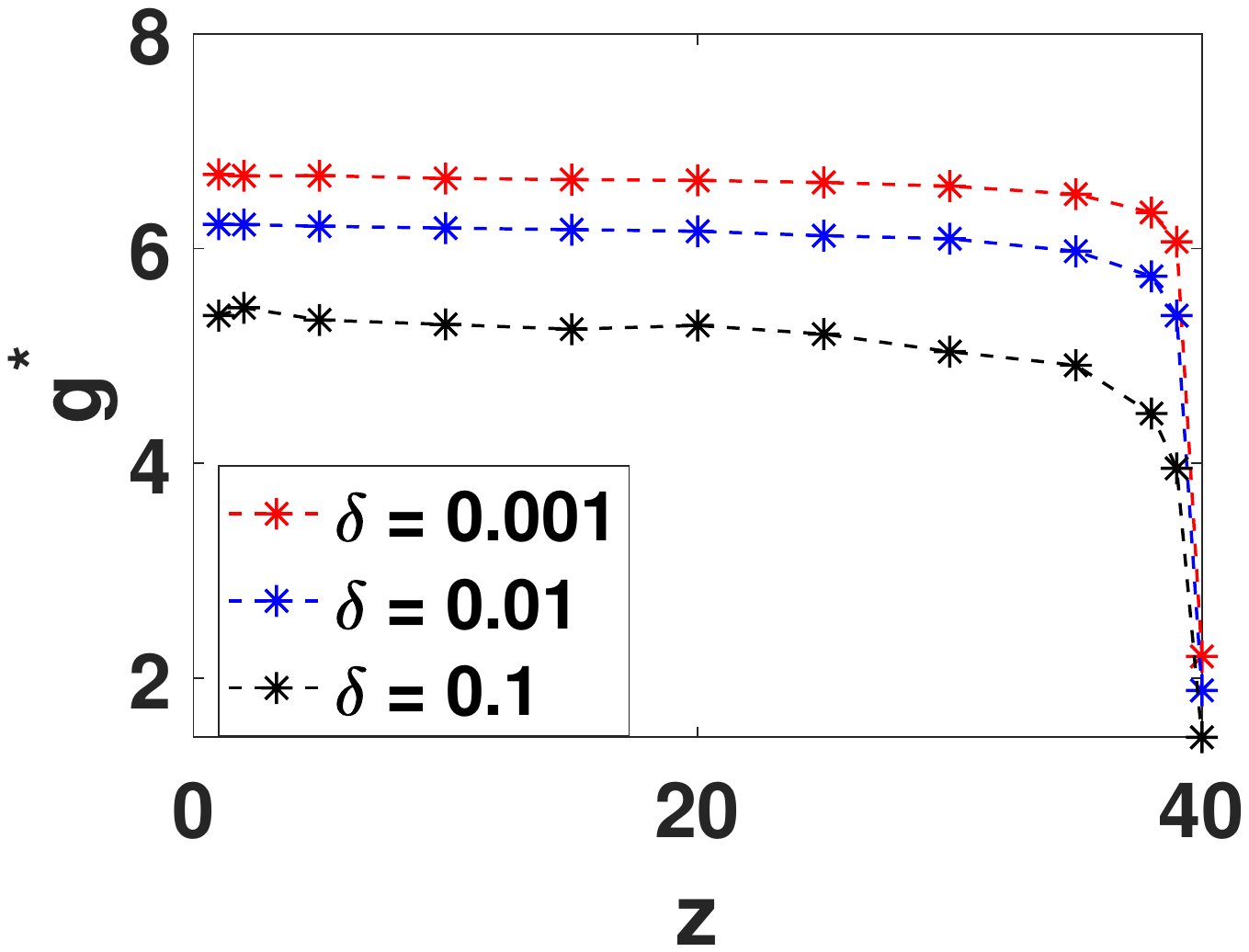}
\end{minipage}\hspace{2mm}
\begin{minipage}{0.48\textwidth}
\includegraphics[trim={4cm 8cm 3cm 7cm}, scale =0.4]{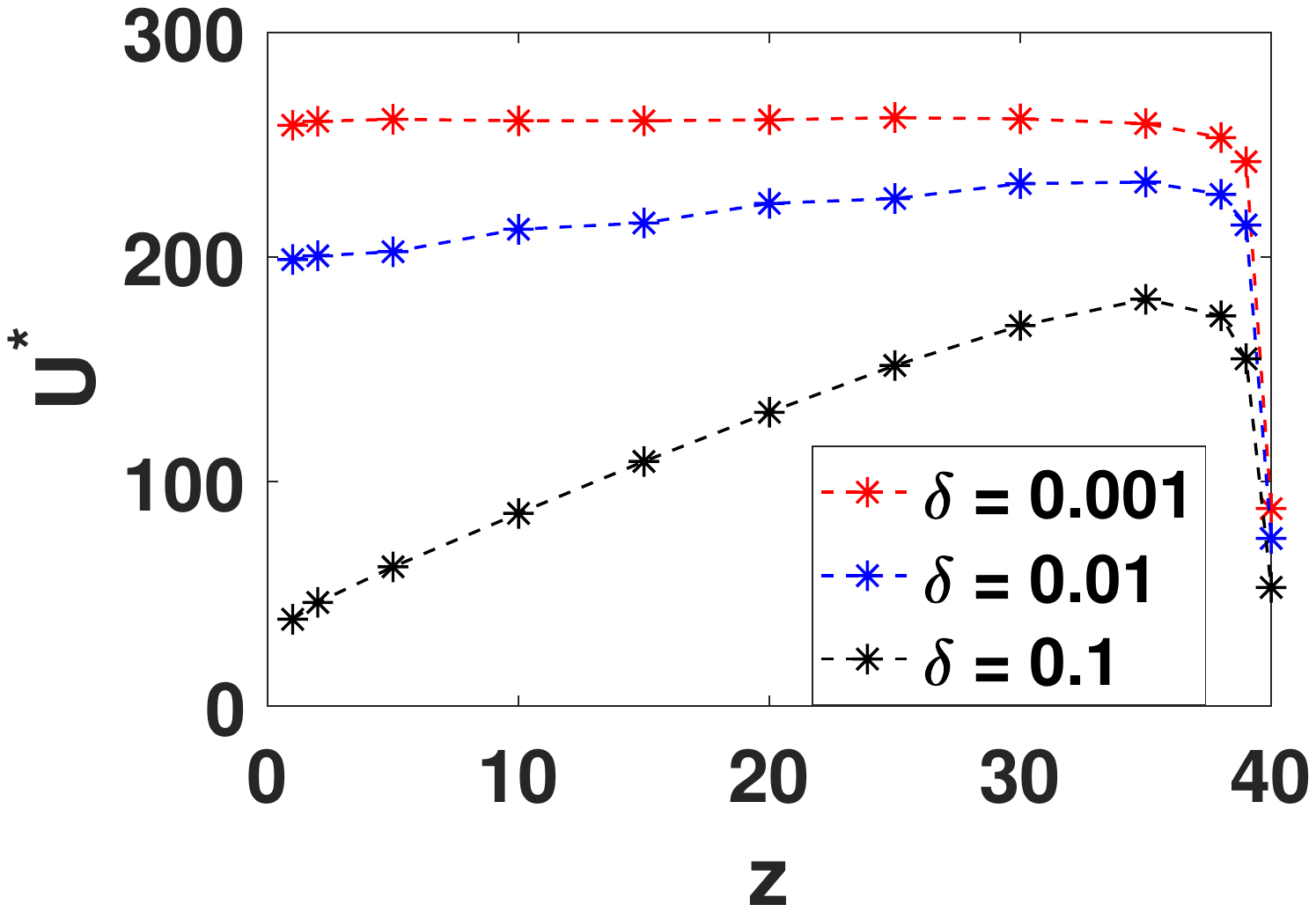}
\end{minipage}
 \caption{{\bf As a function of $\sbaz$}:  $M=40$,  $T=20$, $c=3$, $\mvs=5$, $\sigma^2=2$,  $C_v=1$, $\CI=5$ }
\label{fig:kbar_varies}
\end{figure}

We begin with an example in figure \ref{fig:kbar_varies}, where we plot the optimizers $\og{\sbaz}$ and optimal cost of incentives $\Ui_{\sbaz}^*$ (obtained using theorem \ref{thm_lead_optimal})  as a function of $\baz$ for various $\delta$.
Interestingly, $\og{1}$ is higher than $\og{\sbaz}$ for any $\baz \ge 2$, but $\Ui^*_1$ is smaller than $\Ui_k^*$ for all $k < M-1$.  \textit{Further $U_{\sbaz}^*$ with $\baz$ either  equal to $M$ or $1$ is the smallest among all possible  $\baz$, for all $\delta$.}


\begin{figure}
    \centering
    \begin{minipage}{0.48\textwidth}
    \includegraphics[trim={3cm 8cm 3cm 7.8cm}, scale=0.38]{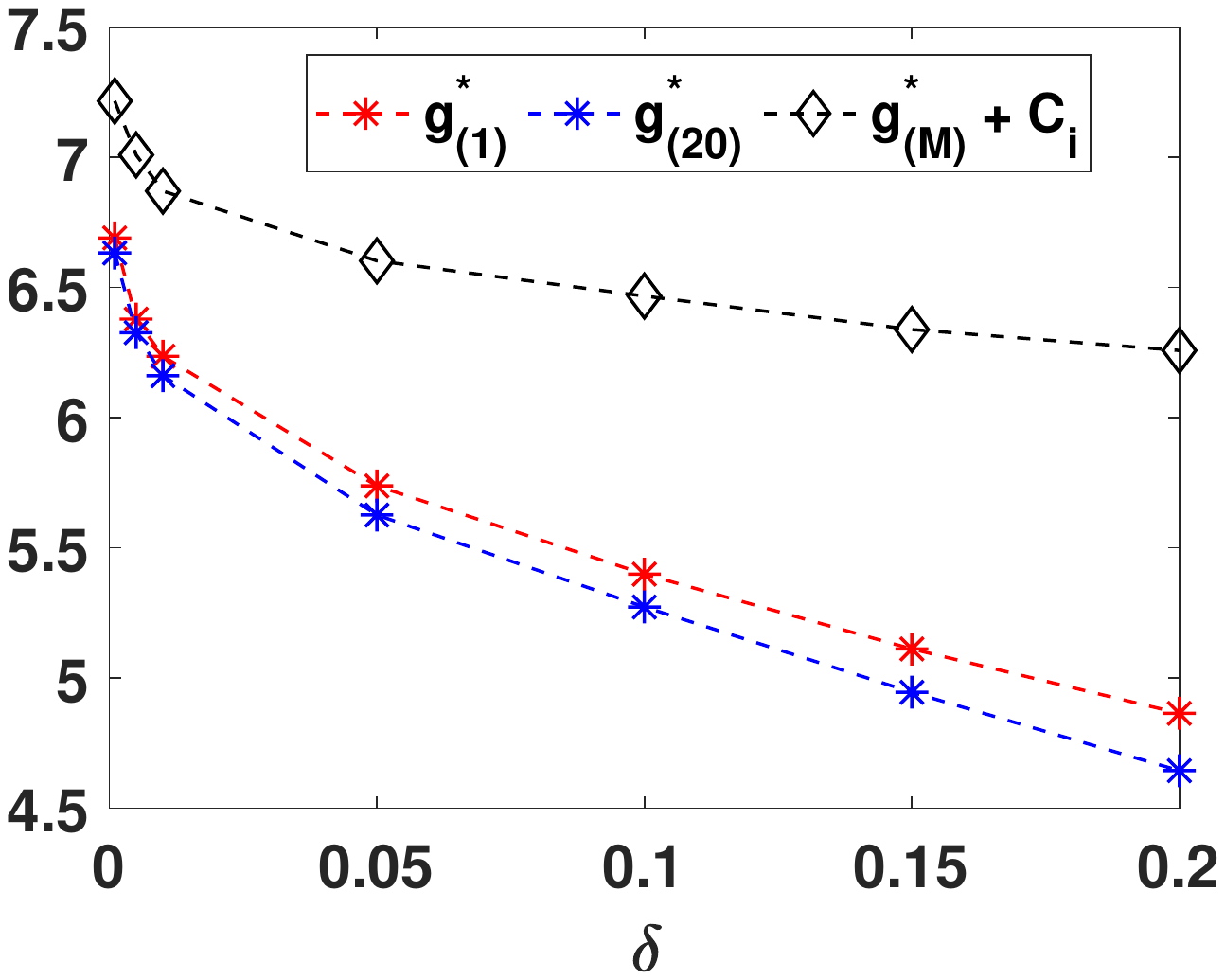}  
    %
    \end{minipage}
    \hspace{3mm}
    \begin{minipage}{0.48\textwidth}
    
    \includegraphics[trim={3cm 8cm 3cm 7.8cm}, scale=0.38]{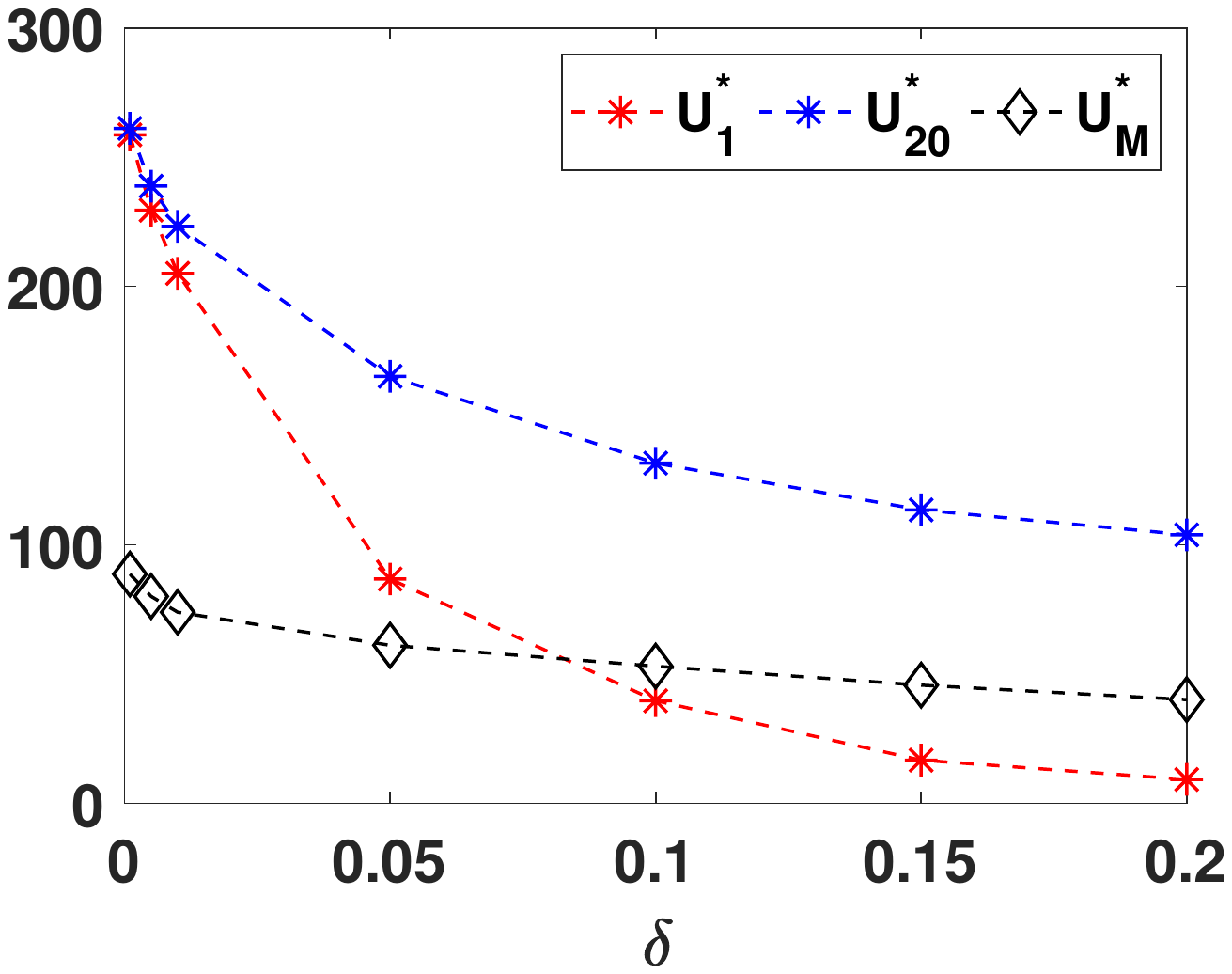}   
    \end{minipage}
      \caption{{\bf As a function of $\delta$:}   $M=40$,  $T=20$, $c=3$, $\mvs=5$, $\sigma^2=2$,  $C_v=1$, $\CI=5$ }
    \label{fig_gk_versus_delta} 
\end{figure}

\begin{figure}
    \centering
    \begin{minipage}{0.48\textwidth}
    \includegraphics[trim={1.5cm 6.5cm 5cm 7.3cm}, scale=0.29]{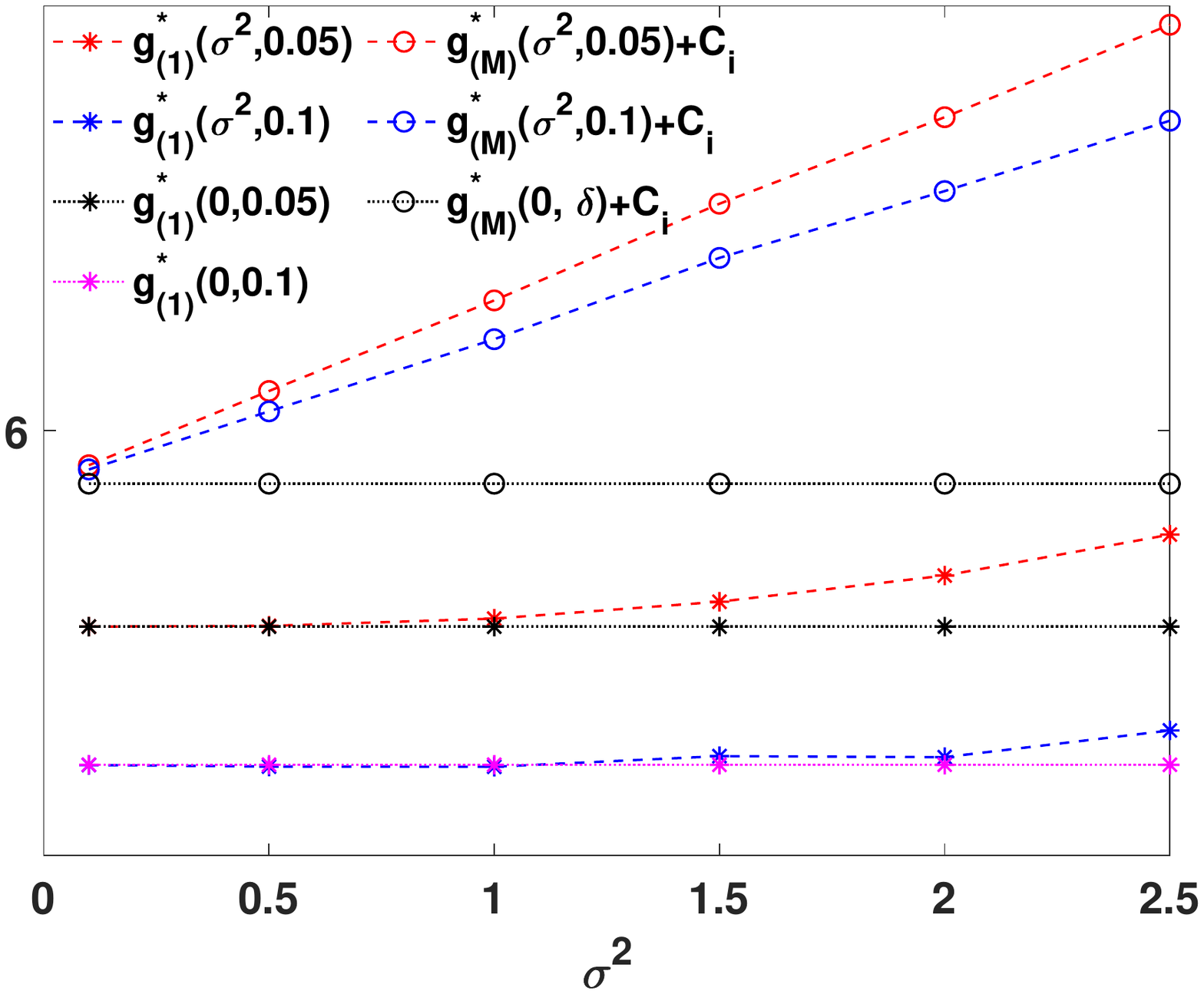}
    \end{minipage}
    \hspace{3mm}
    \begin{minipage}{0.48\textwidth}
    
    \includegraphics[trim={1.5cm 6.5cm 5cm 7.3cm}, scale=0.28]{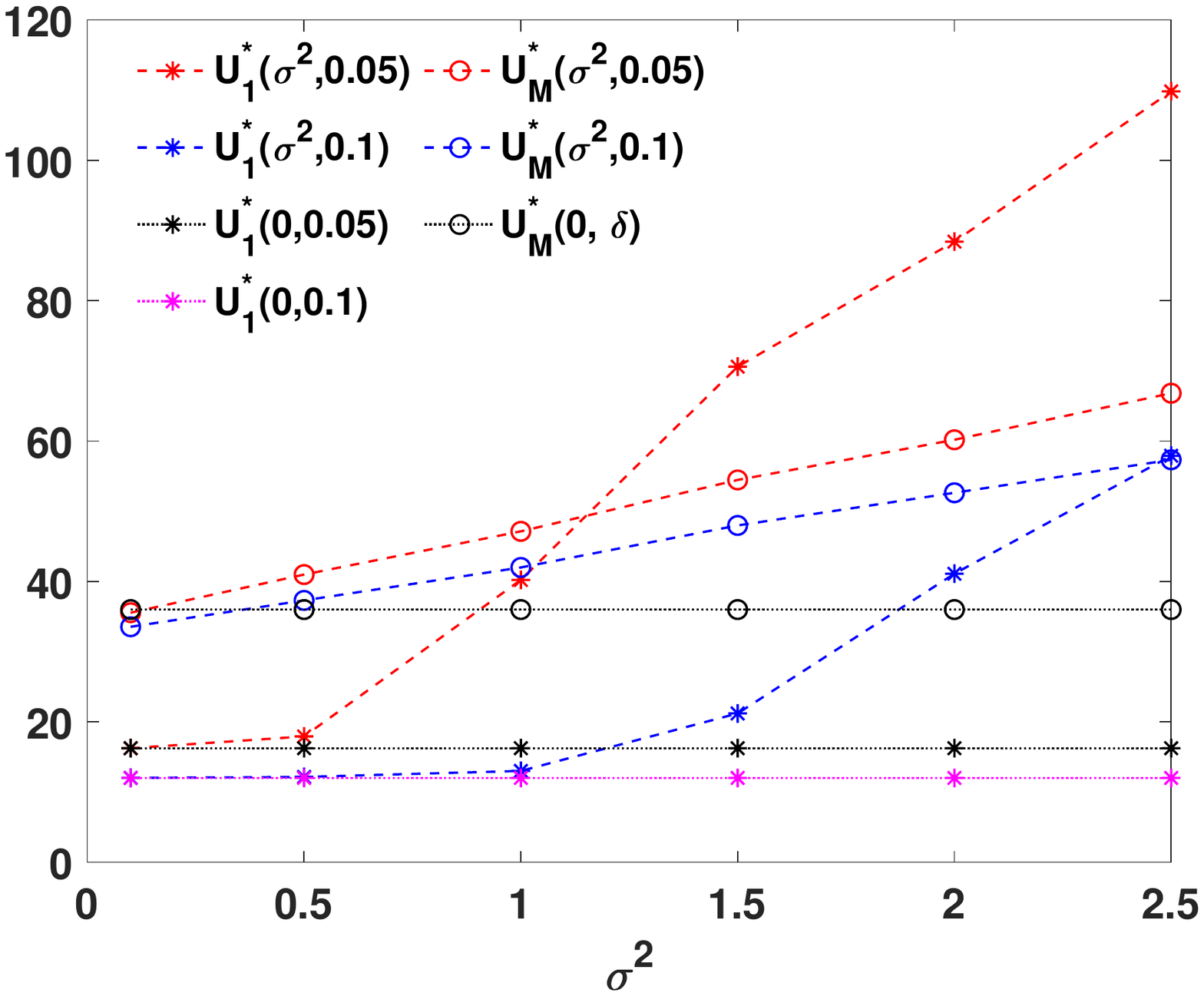}   
    \end{minipage}
    \caption{{\bf As a function of $\sigma^2$:}   $M=40$,  $T=20$, $c=3$, $\mvs=5$,  $C_v=1$, $\CI=5$. Observe from \eqref{eqn_perf_opt}, for zero variance $U_M^*(0, \delta)$ and $\og{M}(0, \delta)$ does not depend upon $\delta$. } \label{fig_versus_sigma_sqr}   
\end{figure}

 We have similar observation as above and more in 
 figure \ref{fig_gk_versus_delta};    we now  plot the convergence of various $\og{\sbaz}$ as a function of  $\delta$;  we  plot $\og{M}+C_i$ in place of $\og{M}$ for ease of comparison. For $\delta \le 0.07$, $U^*_M$ is the best, and when $\delta >0.07$, $U^*_1$ is the best.   
Thus the $\baz$  optimal for incentives is always among $ \{1, M\}$, with $\baz = M$ being optimal for smaller $\delta$.
Further the difference between $\og{\sbaz}$ for any $\baz <M$ and $\og{M}$ converges towards $C_i$ as $\delta \to 0$, reaffirming theorem \ref{thm_gM_smallest}.


In figure \ref{fig_versus_sigma_sqr}, we compare the optimal incentive schemes at  various  variances $\sigma^2$, for different $\delta$.  For clarity, we  provide the plot for $\baz=1$ and $M$, and plot $\og{M}(\sigma^2,\delta) + C_i$ instead of $\og{M}(\sigma^2,\delta)$.
It is clear that as variance reduces to zero, the optimizers $\og{1}(\sigma^2,\delta) \to \og{1}(0,\delta)$ and $ \og{M}(\sigma^2,\delta) \to \og{M,\epsilon}(0,\delta)$, where the limits are given in \eqref{eqn_perf_opt} (actual values represented by dashed lines  with markers, limits with dotted lines), further $U_k^*(\sigma^2,\delta) \to U_k^*(0,\delta)$.  This reaffirms the convergence given in theorem \ref{thm_gM_smallest}.(ii).  

From right sub-figure of figure \ref{fig_versus_sigma_sqr}, we make another interesting observation. The optimal incentive scheme is at  $\baz = M$ for a bigger range of $\delta$, when variance is large -- for $\sigma^2\ge 1$, $U_M^*$ is smaller than $U_1^*$ for all  $\delta \le 0.05$, while it is smaller even for $\delta = 0.1$ when $\sigma^2 = 2.5$.

 In all,   {\it it is optimal for the leader to aim either for $\baz = M$ (where free riding of influencers is completely curbed) or to aim for $\baz = 1$ (where it depends minimally on influencers).} It is also important to observe here that the leader will not achieve complete eradication  without influencers -- from \eqref{eqn_h_func}, $h_v(0;\bnu)>0$  (recall $c_f(0)  = 0$) and hence by lemma \ref{lem_ess}, eradicating attractor is not ESSS. 
Further more, $\baz =M$ provides optimal scheme for small $\delta$ (which is the obvious range of interest), and such  range of $\delta$ increases with variance, and thus we prefer $\baz = M$ in the joint optimization discussed  next.

\subsection{Joint Optimization for Leader} \label{sec_joint_opt}


As already mentioned, the leader's problem has  multiple objectives which we now describe in detail. We also discuss the notions of optimality and the optimizers in this subsection.

The main objective is to ensure the eradication of disease with probability at least $(1-\delta)$. This can be achieved by vaccinating sufficient fraction of the public, which leader would like to keep to the minimum. Thus, the second  objective is to minimize $\ps{e}$ (given  in table \ref{table_ess_candidate}), the fraction of vaccinated population at eradicating ESSS.
Apart from these, the leader would also be interested in an optimal incentive scheme that minimizes $\Ui$ of \eqref{eqn_leader_opt_wait}. 

Previously, we discussed  the optimal incentive scheme for any  VA rate strategy $\bnu$ (or $\baz(\bnu)$).
From  theorems \ref{thm_perf_info_var_z} and \ref{thm_gM_smallest}, we observe that  $\Ui^*_M<\Ui^*_k$ for all $k \le M-1$ for all sufficiently small $\delta$, which is the obvious range of interest. Thus, we have the following definition (recall $\bnu$ is admissible if it satisfies \eqref{eqn_eradication_condition}):
\begin{definition}[Incentive-optimal]
A strategy $(\bg,\bnu)$ is \textit{incentive-optimal} if $\bnu$ is admissible, $\baz(\bnu)=M$ and $g_0=g$ is the optimizer of \eqref{eqn_leader_opt_wait} with $\baz=M$. 
\end{definition}%
Lemma \ref{lem_min_psi} in appendix \ref{appendix_leaders_game} shows that for any admissible $\bnu$, $\ps{e}(\bnu) > \nvdf$  (here $\nvdf$ defined in table \ref{table_ess_candidate} is the fraction of infected population  at non-vaccinating ESSS).  Thus one can not achieve $\nvdf$ but can choose $(\bg,\bnu)$ to get arbitrarily close to it. Hence, we have the following,
\begin{definition}[$\epsilon$-vaccine-optimal]
For any $\epsilon>0$, leader's strategy $(\bg^\epsilon,\bnu^\epsilon)$ is called \textit{$\epsilon$-vaccine-optimal} if $\bnu^\epsilon$ is admissible,  $\ps{e}(\bnu^\epsilon) \le  \nvdf +\epsilon$, and  $g_0^\epsilon = g^\epsilon$ is the optimizer of \eqref{eqn_leader_opt_wait} with $\baz= \baz(\bnu^\epsilon)$.
\end{definition}
 Towards deriving the main result of the paper, for any $0\le k \le M$ define,
\begin{equation}\label{eqn_lk_defn}
    L_k :=\min\left\{-\min \left\{\Bar{c}_{v_2},\frac{c_{v_2}}{\nvdf}\right\}\frac{M-k}{M}, \ \frac{r+b}{r+2b}c_i   -  \bar{c}_{v_2}\frac{M-k}{M}\right\},
 \end{equation}and consider the following inequalities for  $1\le k \le M$ (recall $c_f(0):=0$):

 \vspace{-3mm}
{\small
\begin{eqnarray}
&&c_{v_1} - c_f(k) \le L_k, \mbox{ and, }\ c_{v_1} -c_f(k-1) > L_{k-1}, \mbox{ if } \Bar{c}_{v_2}>\frac{c_{v_2}}{\nvdf}, \label{eqn_k_eps_opt2} \\
&&c_{v_1} - c_f(k) < L_k, \mbox{ and, }\ c_{v_1} -c_f(k-1) \ge L_{k-1},   \mbox{ else. }\label{eqn_k_eps_opt} 
\end{eqnarray}}%
With all the definitions in place, we present our final result. This result says is not always possible to achieve  incentive optimality. However, one can always achieve $\epsilon$-vaccine-optimality, which may have higher importance.

\begin{theorem}\label{thm_iff_conditionn_for_optimality}
i) \textbf{Existence of $\epsilon$-vaccine-optimal policy:} There exists a unique $k \in \{1, \cdots, M\}$ that satisfies \eqref{eqn_k_eps_opt2} or \eqref{eqn_k_eps_opt}, whichever is applicable. 
Let $\bg$ be the optimal incentive scheme with $\baz=k$. Then  for any $\epsilon >0 $ there exist $\bnu^\epsilon$  such that,   $k=\baz(\bnu^\epsilon)$ and $(\bg,\bnu^\epsilon)$ is $\epsilon$-vaccine-optimal; we further  have, 
$$\nu_b^\epsilon \in ( b\rho \nvdf-\epsilon, b\rho \nvdf) \ \mbox{ and } \  \nu_e^\epsilon = b \rho -\frac{\nu_b}{\nvdf}+o(\epsilon).$$
ii) \textbf{Incentive optimality:} There exists an incentive optimal policy  if and only if,

\vspace{-3mm}
{\small
\begin{equation}\label{eqn_inc_opt_cond}
   c_{v_1}-c_f(M-1) > -\frac{\bar{c}_{v_2}}{M} \mbox{ when } \Bar{c}_{v_2}>\frac{c_{v_2}}{\nvdf}, \mbox{ and } c_{v_1}-c_f(M-1) \ge -\frac{\bar{c}_{v_2}}{M} \mbox{ else.}  \hspace{4mm} 
\end{equation}}\vspace{-1mm}\eop
\end{theorem}%
\vspace{-5mm}
\noindent\textbf{Remarks:}  There exists  a policy $(\bg, \bnu)$ that is $\epsilon$-vaccine-optimal as well as incentive optimal for any $\epsilon>0$ if $k$ satisfying \eqref{eqn_k_eps_opt2}-\eqref{eqn_k_eps_opt} equals $M$.

At incentive optimality, $\ps{e}$ can be far away from $\nvdf$. At vaccine-optimality, the cost of incentives is not good, it is higher as also shown in numerical section.

For a given disease characteristics $(\lambda, r,b)$ and public behaviour, if the influence is smaller, i.e, $c_{f}(M-1)$ is less, then one can get both optimality. On the other hand, for the same scenario 
if the influence is high, then leader needs to pay more to the influencers as incentive optimality is not achieved.

If disease is less infectious, i.e., $\nvdf$ is small, then inequality \eqref{eqn_k_eps_opt2} is applicable, and further aspects do not depend upon $\nvdf$. Thus below a threshold on $\nvdf$, the possibility of both optimality or incentive optimality remains the same. Nonetheless, it is obvious that the `best' vaccinated fraction of public, $\ps{e}$, also decreases with $\nvdf$. 

With increase in $\nvdf$, from   \eqref{eqn_lk_defn} and \eqref{eqn_k_eps_opt}, $\epsilon$-optimality is achieved with smaller values of $k$. Thus, leader requires to motivate lesser number of influencers if disease is more infectious; but the irony is that, it is less expensive for the leader to motivate all of the influencers, i.e. $\baz=M$. Another obvious observation is that the fraction of vaccinated public is also higher. As seen from \eqref{eqn_lk_defn}, the severity of disease $c_i$ has similar effect as $\nvdf$, however the vaccinated fraction of public does not alter.

\subsubsection*{Summary of Leader's Optimization and the Stackelberg game}
\begin{itemize}

\item All the policies of this section  ensure non-eradication  probability is below the given threshold $\delta$. We    
    mainly focus on small values of $\delta$ (the range important for this study). 

    \item We first derive the optimal incentive scheme for any given dynamic vaccine supply rate policy  $\bnu$ or equivalently for given $\baz(\bnu)$.


\item For perfect information case, some optimal incentive schemes are  in \eqref{eqn_prf_info}-\eqref{eqn_perf_opt},
while for the general case,  theorem \ref{thm_gM_smallest} does the same. 

\item We also compare the  above incentive schemes with distinct values of $\baz(\bnu)$ and identify that the cost of incentives is minimum either when
\begin{itemize}
    \item $\baz(\bnu) = M$, all  the influencers need to vaccinate and this curbs the  free riding among  the influencers;
    \item or when $\baz(\bnu) = 1$,  minimizing the dependency on influencers. 
\end{itemize}

\item We define two notions of optimality, vaccine-optimality (minimal fraction of public vaccinated) and incentive-optimality (incentives are minimized)

\item In theorem \ref{thm_iff_conditionn_for_optimality},  we establish that vaccine optimal scheme always exists; we also identify the conditions for  existence of incentive optimal scheme.   
    
\end{itemize}

\section{Numerical observations} \label{sec_num_res}
The influence of the influencers on 
the public is captured via a linear curve:
$
c_f (\baz) =  s \baz
$, 
where $s$ is the sensitivity parameter. 
We first depict the variation in leaders strategy as a function of $\nvdf$ and sensitivity parameter $s$. 

\begin{figure}[h]
\begin{minipage}{0.48\textwidth}
\includegraphics[trim={4.5cm 8cm 1cm 7cm}, scale =0.41]{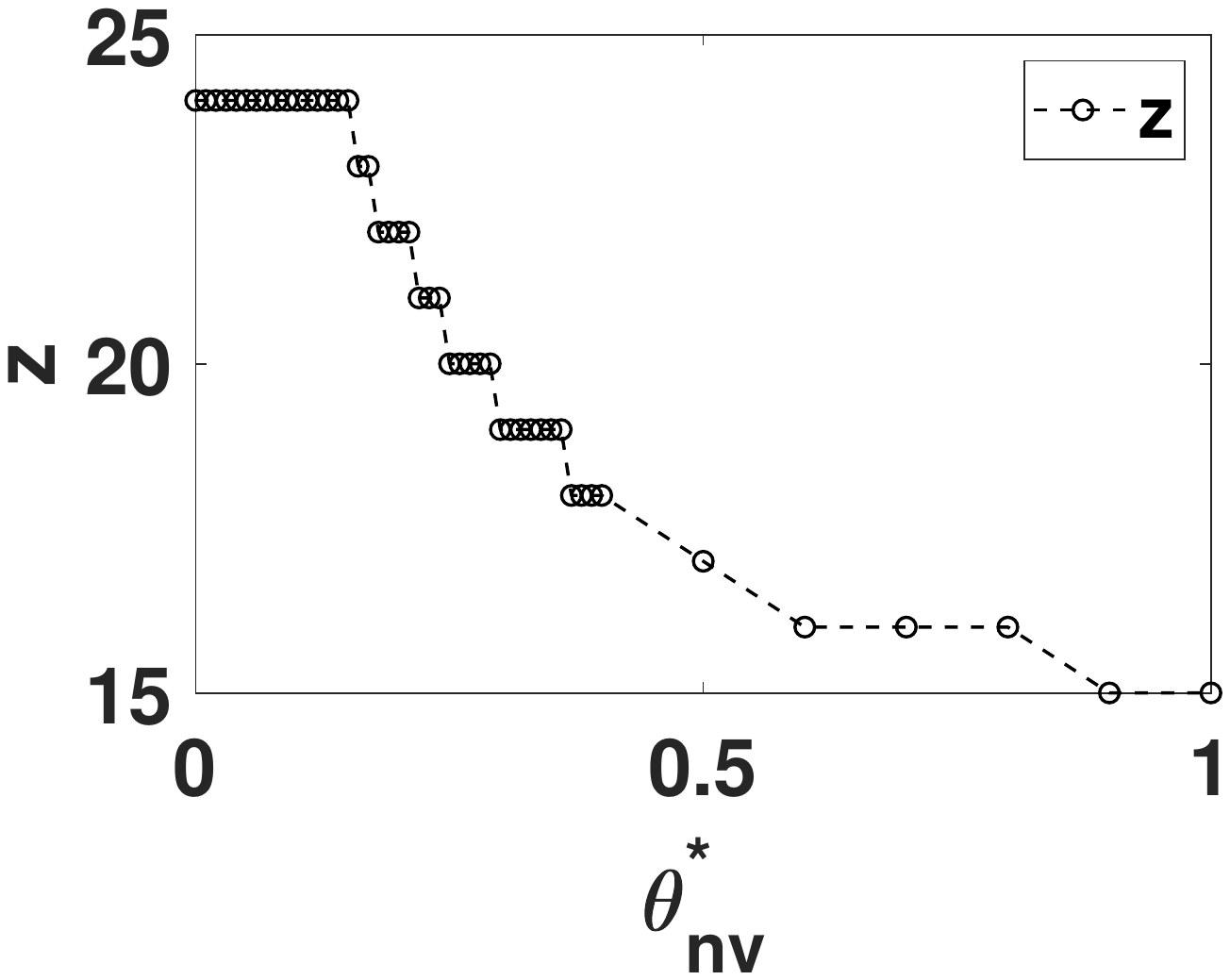}
\caption{\centering{\bf $\sbaz$ versus $\nvdf$}: $M=40$, $c_{v_1}=6$, $c_{v_2}=2$, $\Bar{c}_{v_2}=15$, $c_i=50$, $c_f(k)=0.5k$, $r=5$,
 $b = 2$}
 \label{fig_theta_vs_zbar}
\end{minipage}%
\hspace{4mm}%
\begin{minipage}{0.47\textwidth}
\vspace{-2mm}
\includegraphics[trim={2cm 8cm 1cm 7cm}, scale =0.39]{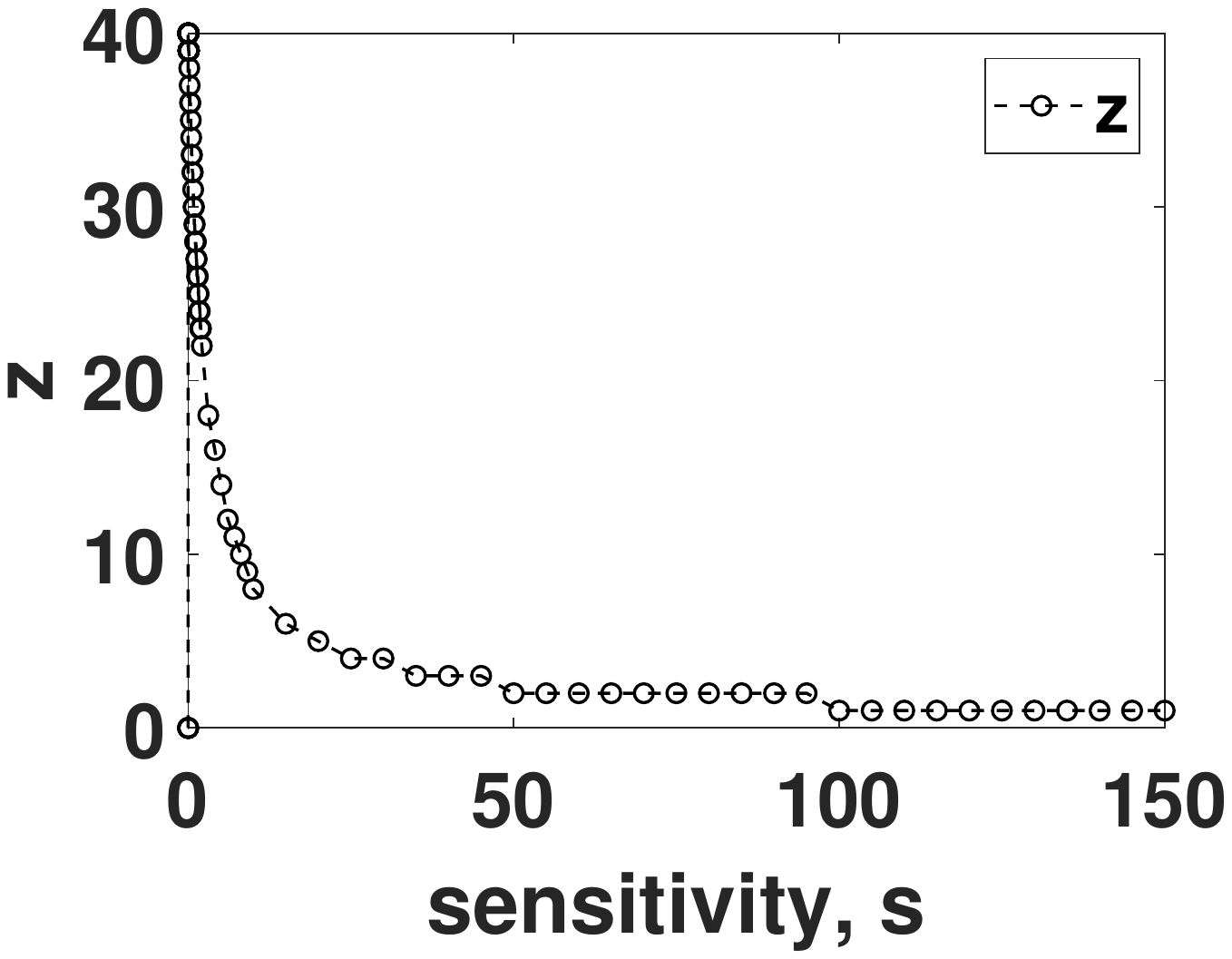}
\vspace{2mm}
\caption{\centering{\bf $\sbaz$ versus $s$}: parameters as in figure \ref{fig_vaccine_optimal_zbar}}
  \label{fig_s_vs_zbar}
\end{minipage}
\end{figure}

From theorem \ref{thm_iff_conditionn_for_optimality}, an $\epsilon$-vaccine-optimal strategy $(\bg,\bnu)$ has  $\baz(\bnu)=\baz$ where $\baz$ uniquely satisfies \eqref{eqn_k_eps_opt2}-\eqref{eqn_k_eps_opt}. Recall, such a $
\baz$ is the target number of influencers that the leader aims to get vaccinated to ensure complete eradication. We plot such $\baz$  as a function of $\nvdf$ and  $s$ in figures \ref{fig_theta_vs_zbar}-\ref{fig_s_vs_zbar}; only one parameter is varied, while the rest are fixed and are given in the respective captions. In figure \ref{fig_theta_vs_zbar}, we observe $\baz$ varies only in a restrictive range with respect to $\nvdf$. Such a variation is evident from \eqref{eqn_k_eps_opt2}-\eqref{eqn_k_eps_opt}, as dependence upon $\nvdf$ comes through $\min\{\Bar{c}_{v_2},\frac{c_{v_2}}{\nvdf}\}$ and the outer minimum in \eqref{eqn_lk_defn}. Thus interestingly, below a certain threshold and above a certain threshold on $\nvdf$, the leaders strategy does not change further with $\nvdf$ (which is indicative of the infection level).
In all, as the disease becomes more infectious, i.e., higher $\nvdf$, we have $\baz<M$, and thus it becomes more expensive for the leader to encourage the influencers towards vaccination to ensure complete eradication. 

\begin{figure}[h]
    \centering
\begin{minipage}{0.48\textwidth}
\includegraphics[trim={2cm 8cm 1cm 7cm}, scale =0.31]{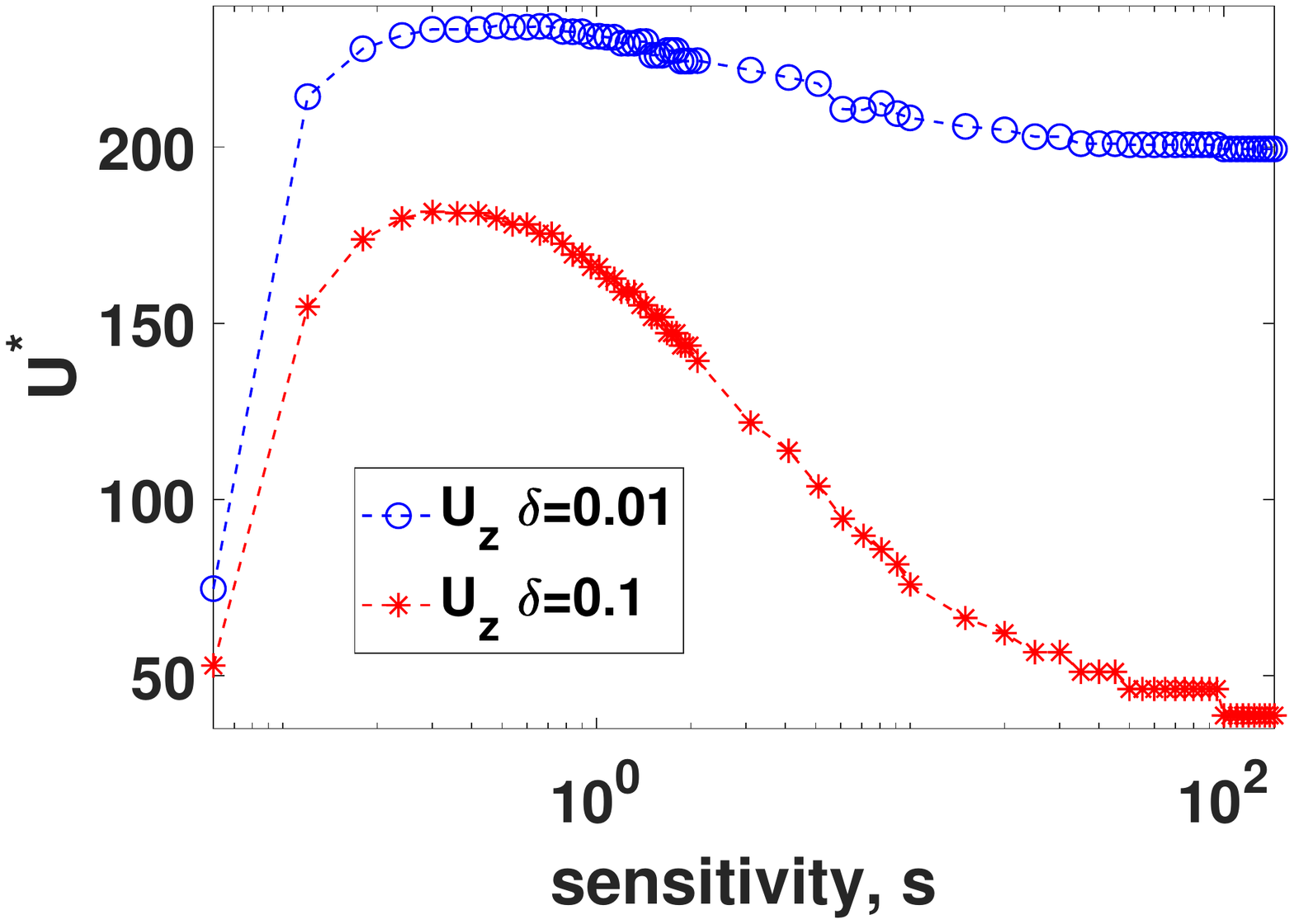}
\end{minipage}\hspace{2mm}
\begin{minipage}{0.48\textwidth}
\includegraphics[trim={2cm 8cm 1cm 7cm}, scale =0.4]{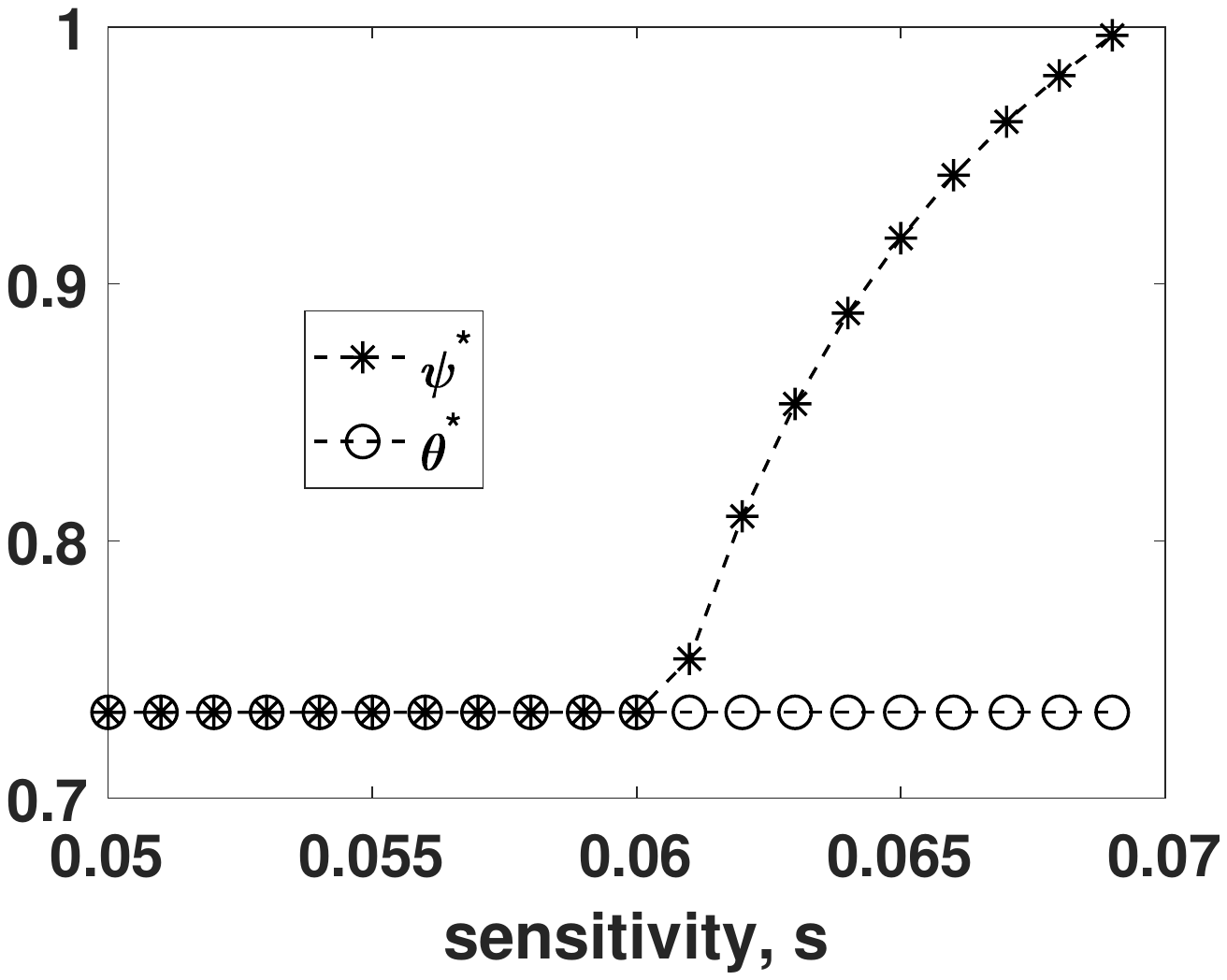}
\end{minipage}
\vspace{3mm}
 \caption{{\bf $\epsilon$-vaccine-optimal $U(\bg,\bnu)$ and incentive-optimal $\ps{e}$ versus sensitivity parameter ($s$)}: $M=40$, $c_{v_1}=0.2$, $c_{v_2}=0.05$, $\Bar{c}_{v_2}=100$, $c_i=0.5$, $c_f(k)=sk$, $r=2$, $\lambda=15$,
 $b = 2$, $T=20$, $c=3$, $\mvs=5$, $\sigma^2=2$,  $C_v=1$, $\CI=5$ }
\label{fig_vaccine_optimal_zbar}
\end{figure}

In contrast to variations in figure \ref{fig_theta_vs_zbar}, $\baz$  varies from $M$ to 1, as the sensitivity parameter $s$ changes in figure \ref{fig_s_vs_zbar}. 
We also plot the corresponding $U^*$,  the cost of incentives at the equilibrium, under  $\epsilon$-vaccine-optimal strategies in the left sub-figure of figure \ref{fig_vaccine_optimal_zbar}. In the right sub-figure,
we   plot   $\ps{e}$, the fraction of vaccinated public at the equilibrium, under   incentive optimal strategies.

We observe, $\baz$ is  decreasing with $s$, however the cost of the leader also depends upon $\delta$ (increases as $\delta$ decreases). More interesting is the comparison with same $\delta$ and for different values of $s$.  When $\delta$ is small ($\delta=0.01$, blue curve) and the population is highly sensitive towards the influencers (higher $s$), it is 
 expensive for the leader to ensure complete eradication as seen in the left sub-figure. 
On the other hand, with larger $\delta$ ($\delta=0.1$, red curve),  it is relatively less expensive for the leader to encourage influencers towards complete eradication when the sensitivity is at the two extremes. Recall with larger $\delta$, $\baz=1$ is the optimizer, and the cost with $\baz=M$ is not too high. 
When sensitivity is moderate, the leader has to pay higher cost for all $\delta$, possibly  because of free-riding behaviour of the influencers.

In the right sub-figure, we plot the regime where incentive optimality is achievable, i.e., \eqref{eqn_inc_opt_cond} is satisfied. In this regime, the incentive cost of the leader, $U^*$ is optimal (for small $\delta$) upon choosing an incentive optimal policy, however, the vaccinated fraction of population can be different from $\nvdf$. When $s \le 0.06$, we have both incentive and vaccine-optimality, but as $s$ increases, we only have incentive optimality, and the fraction of vaccinated population becomes significantly higher than $\nvdf$.

\section*{Conclusions}
It is understood in the literature that voluntary vaccination without incentives may not lead to the eradication  of an infectious disease. In this paper, we explore the possibility of eradicating the same  with the help of influencers.  

There is a minimum  fraction of the public that needs to be vaccinated to ensure complete eradication -- this fraction depends only on disease characteristics. We show that the vaccination of this fraction can be achieved, and such a state is evolutionary stable when a certain number of influencers vaccinate. We further show that there always exists a dynamic vaccine supply strategy and an incentive scheme (optimal for that vaccine supply strategy) for the leader, that motivates this required number of influencers; we call such strategies as vaccine-optimal strategies.

The number of influencers required to ensure complete eradication differs for different vaccine supply strategy choices. Thus the leader can control the required number of influencers for eradication by controlling the vaccine supply rate.
Interestingly, the cost of incentives is minimum when all the influencers are required to vaccinate. It is not so surprising when one observes that such a strategy actually curbs the free-riding behavior among influencers. This is the case when the leader wants to ensure complete eradication with a high level of certainty. 
Suppose the leader only wants to ensure eradication with a relatively smaller level of certainty. In that case, the cost of incentives is minimum when the leader needs to motivate only one influencer towards vaccination. 
In contrast to vaccine-optimal strategies, the existence of such incentive-optimal strategy depends upon disease and public characteristics. 

Interestingly, when the sensitivity of the public towards influencers is either high or low, there exists a strategy that is both optimal. With moderate sensitivity, the leader must choose between incentive and vaccine-optimal policy when the former exists.
 With an  incentive-optimal policy, the incentive cost is minimal, but the vaccinated fraction  can be significantly higher than required and vice-versa for the vaccine-optimal policy. 
 
As the level of infection increases, we observe that the leader needs to motivate fewer influencers. However,  motivating a smaller number is more expensive. It might be possible for the leader to create an illusion that all the influencers are required to vaccinate to ensure complete eradication; if that is possible, the leader can curb the free-riding behavior of the public and influencers more effectively when the infection  rate is high. The study of such a possibility will be an interesting direction for future work.


\begin{appendices}
\section{Proofs for Population Dynamics}\label{appen_ode}

\subsection*{ODE approximation}
The population  dynamics are exactly similar to that in \cite{vaccination} except for the vaccination rates, as explained in section \ref{sec_pop_dyn}; this modification does not alter  the analysis provided there.  We reproduce the results here for completeness.

Our aim is to use the attractors of the ODE \eqref{eqn_ODE} as the limit proportions for the stochastic system for determining the ESSS and further analysis.
In \cite{vaccination}, we along with other coauthors derived a result which provides the justification towards the same in two ways: i) the first part shows that the solutions of the ODE  almost surely approximate the stochastic trajectory over any finite time horizon, and ii) the second part shows that the stochastic trajectory converges towards the attractors of the ODE under an additional assumption {\bf A}.3.  We did not make it very clear in \cite{vaccination}, but  {\bf A}.3 is only required for the second part and the first part is true even without it.  We now state the precise assumption and the theorem statement with $\up_n := (\theta (\tau_n),\psi (\tau_n),\eta (\tau_n))$, the embedded chain, where $\tau_n$ are the transition epochs of the jump process.

\begin{enumerate}[{\bf A.}3]
    \item Let the set $A$ be locally asymptotically stable in the sense of Lyapunov for the ODE \eqref{eqn_ODE}.  Assume that $\{\up_n\}$ visits a compact set  in the domain of attraction  of $A$ infinitely often with probability  $p > 0$.
\end{enumerate}

\begin{theorem}\label{thm_ode_conv}
 i) For every $T>0$, almost surely there exists a sub-sequence $(k_m)$ such that: ($t_k:=\sum_{i=1}^{k} \frac{1}{1+i}$)
            $$
            \sup_{k: t_k \in [t_{k_m}, t_{k_m} + T]} d(\up_k, \up_{*}(t_k - t_{k_m})) \to 0,  \mbox{ as } m \to \infty, \mbox{ where, }
            $$
        $\up_*(\cdot) = (\theta_*(\cdot), \psi_*(\cdot) ,\eta_*(\cdot)  ) $ is the  solution of ODE \eqref{eqn_ODE} with initial condition  $\up_*(0) = \lim_{k_m} \up_{k_m}$. 
        
        ii) Additionally assume  \textbf{A.}3. Then the sequence converges, $\up_n \to A$ as $n \to \infty$ with probability at least $p$. \\
    {\normalfont \textbf{Proof} follows exactly as in \cite{vaccination}.} \eop
\end{theorem}%

Part (i) provides only a partial justification to use the ODE attractors; \textbf{A}.3 is required for part (ii) and hence for complete justification. In the following we would show that the attractors are locally stable in the sense of Lyapunov, but can not comment on the domain of attraction.

\textbf{Common arguments for Lemma \ref{lem_self_erad} and Lemma \ref{lem_attractors}:}   From the structure of ODE \eqref{eqn_ODE}, it is clear that $(\theta,\psi)$ is always contained  in $[0,1]^2$, in fact, $0 \le\theta+\psi \le 1$ if the initial point $(\theta_0,\psi_0) \in [0,1]$ and  $\theta_0+\psi_0 \le 1$. To see this, observe as $\theta$ approaches 1, $\dot{\theta}<0$, and  $\theta$ can not be negative as 0 is an equilibrium point.  Similar arguments hold for $\psi$. Further, as $\theta+\psi$ approaches 1, $\dot{\theta}<0$, and $\dot{\psi}<0$, and hence the claim is true (such arguments are standard in ODE literature \cite{invariance}).  

In all the proofs related to lemmas \ref{lem_self_erad} and \ref{lem_attractors}, we construct appropriate Lyapunov functions that  satisfy the conditions required in \cite[pp 195]{chaos} to show corresponding equilibrium points are attractors. The proof that the functions indeed satisfy the required conditions follows almost as in \cite{vaccination} \TR{and hence the relevant steps are skipped here. However, there are slight changes required, and all the details are provided in arxiv version of this paper \cite{arxive}}{and is provided in the following}.

\noindent\textbf{Proof of Lemma \ref{lem_self_erad}:} 
Let $\up :=(\theta,\psi,\eta)$ and $\up^*=(0,0,\eta^*)$.  \TR{To prove $\up^*$ is an attractor, consider the following Lyapunov  function,
\vspace{-3mm}
\begin{eqnarray*}
     V(\up) &:=&   
       \theta A_\theta^2 +  \psi B_\psi^2 + C_\eta(\eta^* - \eta),\ \  \mbox{ where } \ A_\theta:=1-\theta - \frac{1}{\rho},\\
        B_\psi&:=&(1-\psi)\beta\nu_b-b, \ \mbox{ and } \  \ C_\eta:= \eta^* - \eta.
       \end{eqnarray*}
       }{Let $A:=1-\theta -\psi - \frac{1}{\rho}$, $B:=(1-\theta -\psi)\beta(\nu_b+\nu_e\psi)-b$ and $C:= \frac{b-d }{\varrho} - \eta$. 
Observe that one can re-write the ODE \eqref{eqn_ODE} in a neighbourhood of $\up^*$ as below:

\vspace{-3mm}
{\small
\begin{eqnarray}\label{eqn_ode_0c1_2}
 \dot{\up} = g(\up), \mbox{\normalsize where }   g(\up)= \left(\frac{\theta \lambda}{ \eta\varrho}A,\ \frac{B\psi}{\eta\varrho }, \ 
  C \right).
\end{eqnarray}}%
Now define the  function $
 V(\up) :=   
       \theta A_\theta^2 +  \psi B_\psi^2 + C_\eta(\eta^* - \eta) $, where   $A_\theta := A(\theta, \psi^*, \eta^*)$,  $B_\psi :=(1-\theta^* -\psi)\beta(\nu_b+\nu_e\psi^*)-b$ and $C_\eta := C(\theta^*,\psi^*,\eta)$.} %
 \TR{As mentioned before, the proof is as in \cite[Theorem 2]{vaccination}, with details   in \cite{arxive}.}
{It is clear that $V$ is continuously differentiable, $V(\up^*)=0$ from table \ref{table_ess_candidate}, and $V(\up)>0$ for all $\up \ne \up^*$ in an appropriate neighbourhood of $\up^*$.  Now we  show $\dot{V}(\up)<0$ for all $\up \ne \up^*$ in an appropriate neighbourhood to prove $V$ is a Lyapunov function.  
From \eqref{eqn_ode_0c1_2},

\vspace{-3mm}
{\small
\begin{eqnarray*}
    \dot{V}(\up)= \langle \nabla V, g\rangle =\frac{1}{\eta\varrho}  \left (\lambda \theta A(A_\theta^2 - 2 \theta A_\theta) + B \psi (B_\psi^2 - 2 \psi B_\psi\beta\nu_b) \right) - C (C_\eta  + \eta^*-\eta ),
\end{eqnarray*}}%
Since, $A(\up^*)=A_{\theta^*} < 0$, when $\rho<1$, from the continuity of $A, A_\theta$ one can choose an appropriate neighbourhood of $\up$ such that $\frac{1}{\eta \varrho}\lambda \theta A(A_\theta^2 - 2 \theta A_\theta)<-\epsilon_1$ for $\epsilon_1>0$. If $\rho = 1$ then $A(\up^*)=A_{\theta^*} = 0$, but we again have $A, A_\theta<0$ in a neighbourhood of $\up^*$ since $\theta>0,\psi>0$.

Similarly, $B(\up^*)=B_{\psi^*}<0$ when $\beta<\nicefrac{b}{\nu_b}$, from the continuity of $B, B_\psi$ one can choose an (smaller if required) appropriate neighbourhood  of $\up$ such that $\frac{1}{\eta \varrho} B \psi (B_\psi^2 - 2 \psi B_\psi\beta\nu_b)<-\epsilon_2$ for $\epsilon_2>0$. Towards the last component, observe for any $\up$,

\vspace{-3mm}
{\small
\begin{eqnarray*}
   C (C_\eta  + \eta^*-\eta )= 2C(\eta^*-\eta) = 2 \left(\frac{b-d}{\varrho} -\eta \right)\left(\frac{b-d}{\varrho^*} -\eta \right).
\end{eqnarray*}}Thus, there exists a neighbourhood of $\up^*$ such that, for any $\up \ne \up^*$, if $\frac{b-d}{\varrho^*}-\eta<0$, then $\frac{b-d}{\varrho}-\eta<0$ and if $\frac{b-d}{\varrho^*}-\eta>0$, then $\frac{b-d}{\varrho}-\eta>0$ (possible by continuity of $\varrho$ in $(\theta,\psi)$). Hence in an appropriate neighbourhood, $\dot{V}(\up)<0$. Thus,  $\up^*$ is an attractor.}

Towards the converse, 
say  $\rho > 1$. Consider a neighbourhood such that $\psi \to 0$, if there is no such neighbourhood then we have nothing to prove, $\up^*$ is not an attractor. Thus when started in that neighbourhood, for any $\epsilon>0$ there exist a $\tau_\epsilon< \infty$ such that $\psi< \epsilon$ for all $t \ge \tau_\epsilon$. Hence, from \eqref{eqn_ODE}, there exists some $\epsilon>0$ such that $\dot{\theta}(t)>0$ whenever $\theta(t)<\epsilon$, and $t > \tau_\epsilon$. Thus if  $\theta \to 0$, $\theta(t)< \epsilon$ for all $t \ge \tau'$ for some bigger $\tau'$, and then $\dot{\theta}(t)>0$ for all $t \ge \tau'$, contradicting the convergence to 0.
\eop

\noindent\textbf{Proof of Lemma \ref{lem_attractors}:} We prove that the equilibrium points of ODE \eqref{eqn_ODE} listed in table \ref{table_ess_candidate} are the attractors satisying conditions of footnote \ref{foot_remains1} if and only if the  conditions in the first column of the table \ref{table_ess_candidate} are satisfied. Let $\up:=(\theta,\psi,\eta)$ and the equilibrium point $\up^*:=(\theta^*,\psi^*,\eta^*)$.

\TR{}{\textbf{Case 1:}} Consider the first row of table \ref{table_ess_candidate} and i.e.,  equilibrium point $\up^* = (\nvdf,0,\eta^*)$. \TR{The   Lyapunov function to prove it an attractor is (see \cite{arxive} for proof  details) 
\vspace{-3mm}
\begin{eqnarray*}
    V(\up) &:=&   
    A_\theta^2 +  \psi B_\psi^2 + C_\eta (\eta^* - \eta), \ \mbox{ where } \ A_\theta:=1-\theta  - \frac{1}{\rho},\\
    B_\psi&:= &(1-\nvdf -\psi)\beta\nu_b-b \  \mbox{ and } \ C_\eta:= \eta^* - \eta.
\end{eqnarray*}}{ Let $A:=1-\theta -\psi - \frac{1}{\rho}$, $B:=(1-\theta -\psi)\beta(\nu_b+\nu_e\psi)-b$ and $C:= \frac{b-d }{\varrho} - \eta$. 
Observe that one can re-write the ODE \eqref{eqn_ODE} in a neighbourhood of $\up^*$ as below:

\vspace{-3mm}
{\small
\begin{eqnarray}\label{eqn_ode_c1_2}
 \dot{\up} = g(\up), \mbox{\normalsize where }   g(\up)= \left(\frac{\theta \lambda}{ \eta\varrho}A,\ \frac{B\psi}{\eta\varrho }, \ 
  C \right).
\end{eqnarray}}%
Now define the function, $
 V(\up) :=   
    A_\theta^2 +  \psi B_\psi^2 + C_\eta (\eta^* - \eta)$ \TR{.}{, where   $A_\theta := A(\theta, \psi^*, \eta^*)$,  $B_\psi :=(1-\theta^* -\psi)\beta(\nu_b+\nu_e\psi^*)-b$ and $C_\eta := C(\theta^*,\psi^*,\eta)$.}} %
 \TR{}
{It is clear that $V$ is continuously differentiable, $V(\up^*)=0$ from table \ref{table_ess_candidate}, and $V(\up)>0$ for all $\up \ne \up^*$ in an appropriate neighbourhood of $\up^*$.  Now we  show $\dot{V}(\up)<0$ for all $\up \ne \up^*$ in an appropriate neighbourhood to prove $V$ is a Lyapunov function.  
From \eqref{eqn_ode_c1_2}

\vspace{-3mm}
{\small
\begin{eqnarray*}
    \dot{V}(\up)
    &&= \frac{1}{\eta\varrho}  \left (-2 \lambda A_\theta A + B \psi (B_\psi^2 - 2 \psi B_\psi\beta\nu_b) \right) - C (C_\eta  + \eta^*-\eta ),\\
    && \hspace{-1cm}=  \frac{\psi}{\eta\varrho}  \left ( 2 \lambda\psi (\theta^* - \theta) + B B_\psi^2 \right) - \frac{2}{\eta\varrho}  \left (\lambda A_\theta^2 + B  B_\psi \psi^2\beta\nu_b) \right) - C (C_\eta  + \eta^*-\eta ),\end{eqnarray*}}since $AA_\theta = A_\theta^2 - \psi(\theta^* - \theta)$.
Observe that, $B(\up^*)=B_{\psi^*}<0$ when $\beta<\nicefrac{b \rho}{\nu_b}$, from the continuity of $B, B_\psi$, one can choose an appropriate neighbourhood  (smaller if required)   of $\up^*$ such that following holds (for some $\epsilon_1>0$):

\vspace{-3mm}
{\small
\begin{eqnarray*}
    B<0,\ B_\psi <0, \mbox{and, }  \frac{1}{\eta\varrho}  B  B_\psi^2 <- \epsilon_1.
\end{eqnarray*}}%
Then, arguments for last component follow as  in proof of lemma \ref{lem_self_erad},  we have $\dot{V(\up)}<0$. Thus, by \cite[pp 195]{chaos}, $\up^*$ is an attractor.}

Towards the converse, say $\up^*= (\nvdf,0,\eta^*)$ is the attractor. Now say $\beta >  \frac{b\rho}{\nu_b}$. Consider a neighbourhood of $\up^*$ such that $\theta \to \nvdf$; if no such neighbourhood exists, we have nothing to prove. When started in this neighbourhood,  for every $\epsilon>0$ there exists $\tau_\epsilon$ such that $\theta(t) \in (\nvdf-\epsilon,\nvdf+\epsilon)$ for all $t \ge \tau_\epsilon$. Hence, from \eqref{eqn_ODE}, there exists some $\epsilon>0$ such that $\dot{\psi}>0$  for any $t > \tau_\epsilon$ if further $\psi(t)<\epsilon$. Thus, if $\psi \to 0$,  $\psi(t)< \epsilon$ for all $t \ge \tau'$ for some bigger $\tau'$, and then $\dot{\psi}(t)>0$ for all $t \ge \tau'$, contradicting the convergence to 0.

Now, consider the case when $\beta\psi^*>1$ (second and third row of table \ref{table_ess_candidate}), here $\varpi(\psi^*)=1$.  \TR{To prove $\up^*$ is an attractor, define the following Lyapunov function (proof details provided in \cite{arxive}):
{
\begin{eqnarray*}
 V(\up) &:=&    \left\{\begin{array}{ll}
       \theta A_\theta^2 +   B_\psi^2 + C_\eta (\eta^*- \eta), &\mbox{ if }  \nu_e> b \rho -\frac{\nu_b}{\nvdf},\\
   w_1 (\htheta - \theta) A_\theta + w_2 (\hpsi - \psi) B_\psi + C_\eta(\eta^* - \eta),  & \mbox{ if } \nu_e < b \rho -\frac{\nu_b}{\nvdf}, \\
 \end{array}\right. \\   
 \mbox{ where}&& A_\theta := 1-\theta-\psi^* -\frac{1}{\rho}, \ \  B_\psi :=(1-\theta^* -\psi)(\nu_b+\nu_e\psi^*)-b, \\
 C_\eta &:=& \eta-\eta^* \  \mbox{ and $w_1$, $w_2$ are appropriate constants.}
\end{eqnarray*}}}{Then one can re-write the ODE \eqref{eqn_ODE} as below in an appropriate neighbourhood of $\up^*$:

\vspace{-3mm}
 {\small
\begin{eqnarray}\label{eqn_ode_c3_4}
 \dot{\up} = g(\up), \mbox{ where }   g(\up)= \left(\frac{\theta \lambda}{ \eta\varrho}A,\ \frac{B}{\eta\varrho }, \ 
  C \right),
\end{eqnarray}}with $A,C$ as in \eqref{eqn_ode_c1_2} and $B:=(1-\theta -\psi)(\nu_b+\nu_e\psi)-b\psi$.} 

\TR{}
{\textbf{Case 2:} Say $\nu_e> b \rho -\frac{\nu_b}{\nvdf}$; here $\up^*=(0,\ps{e},\frac{b-d}{\varrho^*})$.  Define the following Lyapunov function:
$$V  = \theta A_\theta^2 +   B_\psi^2 + C_\eta (\eta^*- \eta),$$
where   $A_\theta := A(\theta, \psi^*, \eta^*)$,  $B_\psi :=(1-\theta^* -\psi)(\nu_b+\nu_e\psi^*)-b\psi^*$ and $C_\eta := C(\theta^*,\psi^*,\eta)$. Now, from \eqref{eqn_ode_c3_4},

\vspace{-3mm}
{\small
\begin{eqnarray*}
    \dot{V}(\up) = \frac{1}{\eta\varrho}  \left (\lambda \theta A(A_\theta^2 - 2 \theta A_\theta)- 2B B_\psi(\nu_b+\nu_e\psi^*)\right) - C (C_\eta  + \eta^*-\eta ).
\end{eqnarray*}}Now, since $\nu_e > b \rho - \frac{\nu_b}{\nvdf}$, $A(\up^*) =A_{\theta^*}<0$. Choose an appropriate neighbourhood of $\up^*$ such that, 

\vspace{-3mm}
{\small
\begin{eqnarray*}
    A<0,\ A_\theta <0,\ \ \frac{\lambda \theta}{\eta\varrho}  A A_\theta^2 <- \epsilon_1, \mbox{ and } 2B B_\psi(\nu_b+\nu_e\psi^*) < \frac{\epsilon_1}{2},
\end{eqnarray*}}then first component is negative. Arguments for last component follow as before and we have $\dot{V(\up)}<0$.

\textbf{Case 3:} Now we are only left with the case when $\nu_e < b \rho - \frac{\nu_b}{\nvdf}$ and $\beta \psi^* >1$; here $\up^*=(\tht{o},\ps{o},\frac{b-d}{\varrho^*})$.  Define the following function which we will prove to be Lyapunov:
\begin{align}
    V(\up) := w_1 (\htheta - \theta) A_\theta + w_2 (\hpsi - \psi) B_\psi + C_\eta(\eta^* - \eta), \mbox{ where }
\end{align} $w_1$, $w_2$ are appropriate constants will be chosen later,
$A_\theta$, $B_\psi$ and $C_\eta$ are as  in case 3; observe $A_{\theta^*}= B_{\psi^*}=0$. 
Clearly,  $V$ is continuously differentiable, 
$V(\Upsilon^*) = 0$ and   $V(\Upsilon) > 0$ for all $\Upsilon \ne \Upsilon^*$.

For simplicity in notations, call $\ttheta := \htheta - \theta$ and $\tpsi := \hpsi - \psi$.    
The derivative of $V(\up(t))$ with respect to time is:

\vspace{-2mm}
{\small    \begin{align}\label{eqn_derv_V_lemma}
     \dot{V} 
     &= -\left(A_\theta + \ttheta \right)  \frac{A\theta \lambda w_1}{\eta \varrho}  - \left(B_\psi + \tpsi(\nu_b + \nu_e \psi^*) \right)  \frac{B w_2}{\eta \varrho} - (C_\eta + \eta^*- \eta)C .
    \end{align}}
    
One can prove that the  last component,  $- (C + \eta^*- \eta) C$ is strictly negative in an appropriate neighborhood of $\up^*$ as in other cases.  Now, we proceed to prove that other terms in $\dot{V}$   are also strictly negative in a neighborhood of $\up^*$.

Consider the term\footnote{Observe that $A_{\theta} A = (A_{\theta^*} + \ttheta ) A$, and $A_{\theta^*} = 0$.} $\left(A_{\theta} + \ttheta \right)A$, call it $A_1$:

\vspace{-2mm}
{\small
\begin{align}\label{A_term}
    \begin{aligned}
A_1 & = 2\ttheta A = 2\left( \ttheta^2 + \ttheta \tpsi \right) = 2\left( \left( \ttheta c_1 + \frac{1}{2c_1}\tpsi\right)^2 + (1-c_1^2)\ttheta^2 - \frac{1}{4c_1^2} \tpsi^2  \right),
\end{aligned}
\end{align}}where $c_1$ will be chosen appropriately in later part of proof. Further,  let $\nnu :=\nu_b + \nu_e \psi$ and $\nnus := \nu_b + \nu_e \psi^*$, and define $B_1 := \left(B_\psi + \tpsi \nnus \right)B$. Then 
\begin{align}\label{B_term}
    \begin{aligned}
    B_1 &= 2\tpsi \nnus B= 2\tpsi \nnus  (1-\theta-\psi) (\nnu  -\nnus )  + 2\tpsi^2 b \nnus  + 2\nnus (\ttheta \tpsi + \tpsi^2).
\end{aligned}
\end{align}
Note that the first term in \eqref{B_term} can be written as:
\begin{eqnarray*}
   B_2 & = &  2 \tpsi  \nnus  (\nnu  -\nnus )   (1 + \ttheta +\tpsi - \htheta - \hpsi)  =  2 \tpsi  \nnus  (\nnu  -\nnus )  (1 + \ttheta + \tpsi - \nvdf),\\
   & =&   2   \nnus  (\nnu  -\nnus ) [  \tpsi(1 + \tpsi - \nvdf) + \ttheta\tpsi],\\
   &\stackrel{a}{=}&  2   \nnus  [-\nu_e  (1 + \tpsi - \nvdf)\tpsi^2 - (\nu_e  \tpsi) \ttheta\tpsi] =  2 \nnus  [ p_1(\up) \tpsi^2 - p_2(\up)\tpsi \ttheta ],
\end{eqnarray*}%
with $p_1(\up):= -\nu_e (1-\nvdf + \tpsi)$ and $p_2 := \nu_e \tpsi$; $a$ follows using $\nnu-\nnus= - \nu_e \tpsi$. Observe that, in an $\epsilon$-neighbourhood of $\up$ such that $-\epsilon < \tpsi < \epsilon$, we have  $-\nu_e (1-\nvdf) - \nu_e \epsilon \le p_1(\up) \le -\nu_e (1-\nvdf)  
 + \nu_e \epsilon$ and $ - \nu_e \epsilon \le p_2(\up) \le 
  \nu_e \epsilon$. 
Now using $B_2$, one can re-write $B_1$  as:

\vspace{-2mm}
{\small 
\begin{align}\label{eqn_B1_term}
    \begin{aligned}
    B_1 &=2 \nnus  \left [ \left ( b + 1 +  p_1(\up) \right ) \tpsi^2 + 
    \left ( 1 - p_2(\up) \right ) \tpsi\ttheta
    \right  ]
    \\
    &= 2\nnus\left[ \tpsi^2\Big( b + p_1(\up) + 1-c_2^2 \Big)  -  \frac{1}{4c_2^2} \left(1 - p_2(\up) \right)^2 \ttheta^2 \right]\\
&\hspace{1cm}+ 2\nnus  \left(c_2 \tpsi + \frac{1}{2c_2} \left(1 - p_2(\up) \right) \ttheta \right)^2.
\end{aligned}
\end{align}}%
Thus, we get: $\dot{V} 
     <  -A_1 \frac{  \theta \lambda w_1}{\eta \varrho}  - B_1 \frac{w_2}{\eta \varrho} $ (recall the last component in \eqref{eqn_derv_V_lemma} is strictly negative). Now, for $\dot{V}$ to be negative, we need (using terms, in \eqref{A_term} and \eqref{eqn_B1_term}, corresponding to $\ttheta^2, \tpsi^2$):

\vspace{-2mm}
{\small 
\begin{align*}
    \nnus   &\leq \frac{4 c_2^2   }{ w_2 \left( 1- p_2(\up) \right)^2 }  \lambda w_1 (1-c_1^2) \theta, \mbox{ and }\\
     \frac{1}{2c_1^2}\theta \lambda w_1 &\leq  2 \nnus  w_2   \Big(b+  p_1(\up) +1-c_2^2\Big).
\end{align*}}To this end, it is sufficient to choose a further smaller neighborhood of $\up^*$ (call it again $\epsilon$-neighborhood, and such that $\htheta - \epsilon >0$) and some of the parameters $c_1, c_2, w_1$ and $w_2$  such that $1 - c_1^2, 1 - c_2^2 > 0$ and following is true: 
\begin{align}\label{eqn_bound}
 \nnus   &=\frac{4 c_2^2   }{ w_2 \left( 1 +\nu_e \epsilon\right)^2 }  \lambda w_1 (1-c_1^2) (\theta^* - \epsilon), \mbox{ and } \nonumber \\
    \frac{1}{2c_1^2}(\htheta + \delta) \lambda w_1 &\leq  2 \nnus  w_2 \Big( b + 1-c_2^2  -\nu_e (1-\nvdf +\epsilon) \Big).
\end{align}

Let $c_1^2 := \Omega_1 w_1$, and $c_2^2 := \Omega_2 w_2$, for appropriate $\Omega_1, \Omega_2$.
Then, by substituting the first equality in the second inequality of  \eqref{eqn_bound}: 

\vspace{-2mm}
{\small
\begin{align*}
    (\htheta +\epsilon)  &\leq \frac{ 2c_1^2    4 c_2^2}{ w_2 \left( 1+\nu_e \epsilon\right)^2 }   (1-c_1^2) (\htheta - \epsilon) 
    2 w_2    \Big( b  + 1-c_2^2 -\nu_e (1-\nvdf +\epsilon)  \Big)
    \\
    &=     (\htheta - \epsilon) \Omega_1 w_1 (1- \Omega_1 w_1) \Omega_2 w_2 \Big( 1   + b -\nu_e (1-\nvdf +\epsilon)  -\Omega_2 w_2 \Big)\left( \frac{ 4}{1 +\nu_e \epsilon} \right)^2 \\
    &\leq    (\htheta - \epsilon)  \frac{1}{4}  
    \frac{1}{4}\Big( 1+ b -\nu_e (1-\nvdf +\epsilon)    \Big)^2 \left( \frac{  4 }{1 +\nu_e \epsilon} \right)^2 \\
    &=     (\htheta - \epsilon)   \left ( \frac{ 1+ b -\nu_e (1-\nvdf +\epsilon) }{1 +\nu_e \epsilon} \right )^2 , 
\end{align*}}
The last step is obtained by choosing   $\Omega_1w_1=1/2$ that provides the maximum value $1/4$ (for product  $\Omega_1w_1 (1-\Omega_1w_1)$) and $\Omega_2 w_2$ that provides the maximum value $( 1+ b -\nu_e (1-\nvdf +\epsilon)  )^2/4 $ for term $\Omega_2 w_2 \Big( 1+ b -\nu_e (1-\nvdf +\epsilon)  -\Omega_2 w_2 \Big)$. 

Now observe that, 
\begin{eqnarray}\label{eqn_cond_lemma} 1+ b -\nu_e (1-\nvdf)  = \frac{\rho +b\rho - \nu_e}{\rho} >1 \ \ \ (\mbox{as } \nu_e <  b\rho - \frac{\nu_b}{\nvdf})
\end{eqnarray}%
Thus \eqref{eqn_bound} holds in an appropriate $\epsilon$-neighbourhood, and we have the result.}

Towards the converse, $\beta\psi^*>1$ is true since $\up^*$ satisfies condition of footnote \ref{foot_remains1}. 
Now, say $(0,\ps{e},\eta^*)$ is the attractor  and $\nu_e < b \rho -\frac{\nu_b}{\nvdf}$. This implies $1-\ps{e} > \frac{1}{\rho}$. Consider a neighbourhood of $\up^*$ such that $\psi \to \ps{e}$ and $\theta \to 0$, if no such neighbourhood exists, we have nothing to prove. When started in this neighbourhood, for any $\epsilon>0$, there exists a $\tau_\epsilon$ such that $\psi(t) \in (\ps{e}-\epsilon, \ps{e}+\epsilon)$ for all $t \ge \tau_\epsilon$.  Thus, from \eqref{eqn_ODE},  we have $\dot{\theta}(t)>0$ when $\theta(t)<\epsilon$ and $t \ge \tau_\epsilon$ for some appropriate $\epsilon$. Thus the contradiction to convergence of $\theta$ to 0.

Lastly, say $(\tht{o},\ps{o},\eta^*)$ which also satisfies the conditions (e.g., $\tht{o}\ne 0$) of footnote \ref{foot_remains1}  is the attractor. This implies $\tht{o}>0$, which implies $\nu_e < b \rho -\frac{\nu_b}{\nvdf}$. Thus the proof.\eop

\begin{lemma}\label{lem_cont_of_attract}
Let $\up^*:=(\theta^*,\psi^*,\eta^*)$ be the attractor of the  ODE   \eqref{eqn_ODE} under population response $\varpi(\psi)$, that satisfies $q^*\in\{0,1\}$ and conditions of table \ref{table_ess_candidate}. Let $\up_\epsilon^*:=(\theta^*_\epsilon,\psi^*_\epsilon,\eta^*_\epsilon)$ be  attractor of ODE corresponding to $\epsilon$-mutant of this population response,  $\varpi_{\epsilon}(q)$ for some $q\in(0,1)$. Then, \\
i) there exists an ${\bar \epsilon} (q) > 0$ such that the
attractor $\up_\epsilon^*$ is unique  and  is a continuous function of $\epsilon$ 
for all $\epsilon \le {\bar \epsilon}$ with $\up^*_0 = \up^*$. 

\noindent ii) Further  $\bar{\epsilon}$ could be chosen such that the sign of $h(\theta^*_\epsilon,\psi^*_\epsilon)$ of \eqref{eqn_utility} remains the same as that of $h(\theta^*,\psi^*)$ for all $\epsilon \le {\bar \epsilon}$, when the latter is not zero. 
\end{lemma}
\noindent {\bf Proof:} Let  $b_\epsilon := \frac{b}{\epsilon q + (1-\epsilon)}$. 
Say $\rho>1$, $0\le \nu_e< b\rho -\frac{\nu_b}{\nvdf}$. From table \ref{table_ess_candidate}, $\up^* =(\tht{o},\ps{o}, \eta^*)$ for all $\beta>\frac{1}{\ps{o}}$ and $q^*=1$. Now consider the ODE corresponding to  $\varpi_\epsilon(q)$ (obtained by replacing $\varpi(\psi)$ by $\varpi_\epsilon(q)$ in \eqref{eqn_ODE}) in a neighbourhood of $\up^*$.  One of the equilibrium points (call it $\up_\epsilon^*$) of this ODE is given by  zero of the following:

 \vspace{-3mm}
{\small
\begin{align}\label{eqn_perturbed_ODE_lemma}
1-\theta -\psi - \frac{1}{\rho}   ,\ \ (1-\theta -\psi)(1-\epsilon + \epsilon q)(\nu_b+\nu_e\psi)-b\psi,\ \ \frac{b-d }{\varrho} - \eta. 
\end{align}}%
By direct computation, $\up_\epsilon^*=(\theta^*_\epsilon,\psi^*_\epsilon, \eta^*_\epsilon)$ with $\theta^*_\epsilon = \nvdf - \psi^*_\epsilon$,   $\psi^*_\epsilon  = \frac{\nu_b}{ b_\epsilon \rho - \nu_e}$,  and $\eta^*_\epsilon = \frac{b-d}{\varrho_\epsilon}$ where $\varrho_\epsilon$ is $\varrho$ of \eqref{eqn_ODE} computed at $(\theta^*_\epsilon,\psi^*_\epsilon)$. Observe $\up_\epsilon^*$ is unique, continuous in $\epsilon$ (in some ${\bar \epsilon}$-neighbourhood) and coincides with $\up^*$ at $\epsilon =0$.
Further, arguing as in  proof of lemma \ref{lem_attractors}, (using corresponding Lyapunov function  with obvious modifications), $\up_\epsilon^*$ is the attractor.

Now consider the case when $\nu_e> b\rho -\frac{\nu_b}{\nvdf}$. From table \ref{table_ess_candidate}, $\up^* = ( 0,\ps{e},\eta^*)$ for all $\beta>\frac{1}{\ps{e}}$. The corresponding zero $\up^*_\epsilon$  under mutation policy $\varpi_\epsilon(q)$ equals $\theta^*_\epsilon=0$,

\vspace{-3mm}
{\small
$$\psi^*_\epsilon= \frac{-(b_\epsilon+\nu_b-\nu_e)+\sqrt{(b_\epsilon +\nu_b-\nu_e)^2 + 4 \nu_e \nu_b}}{2 \nu_e } \indc{\nu_e>0}  + \frac{\nu_b}{b_\epsilon + \nu_e},$$}%
 and $\eta^*_\epsilon=\frac{b-d}{\varrho_\epsilon}$. Once again,  $\up^*_\epsilon$ is unique, continuous in $\epsilon$ and coincides with $\up^*$ at $\epsilon=0$. Further arguing as in  proof of lemma \ref{lem_attractors},  $\up_\epsilon^*$ is also the attractor.

Similarly, when $\beta<\frac{b\rho}{\nu_b}$, from table \ref{table_ess_candidate} $\up^* = (\nvdf, 0,\eta^*)$.  The corresponding zero $\up_\epsilon^*$ equals $\up^*$ in a neighbourhood. Once again the Lyapunov arguments go through as before. 

Lastly, say $\rho\le 1$, from table \ref{table_ess_candidate} $\up^*=(0,0,\eta^*)$ for all $\beta<\frac{b}{\nu_b}$. For an appropriate $\bar{\epsilon}$, $\up^*_\epsilon=\up^*=(0,0,\eta^*)$ for all $\epsilon<\bar{\epsilon}$ and Lyapunov arguments go through as before.  This proves part (i).

The second part follows by continuity of $h$ function \eqref{eqn_utility}.  \eop

\noindent\textbf{Proof of Lemma \ref{lem_ess}:} The proof follows in  exact similar lines as in \cite{vaccination}. We discuss the same in the following for the sake of completion. In the following, we construct appropriate population response   $\varpi = \min \{1,\beta\psi\}$ that are ESSS for different cases.

Say $\rho\le1$ and consider the population response  $\varpi $ with $\beta<\frac{b}{\nu_b}$.  From lemma  \ref{lem_self_erad}, the attractor of ODE \eqref{eqn_ODE} equals $(0,0,\eta^*)$. Then $h(\theta^*,\psi^*)$ in \eqref{eqn_utility} equals $c_{v_1} +\left(1-\frac{\zt}{M}\right){\bar c}_{v_2} >0$. Thus $q^*=\varpi(\psi^*) =0 $, and best response set $ {\cal B} (\varpi)$ in \eqref{eqn_BR_ess} equals $\{0\}$, satisfying first condition for ESSS (see definition \ref{def_ess}). 
Further, under the mutation policy $\varpi_\epsilon(q)$, from lemma \ref{lem_cont_of_attract}, $h(\theta^*_\epsilon,\psi^*_\epsilon)>0$, hence $ {\cal B} (\varpi_\epsilon(q)  )  =\{0\}$ satisfying second condition. Thus $\varpi$ is an ESSS. Now say $\rho>1$.

When  $h_i(\zt,\bnu)>0$,  consider   $\varpi$ with $\beta<\frac{b\rho}{\nu_b}$.
From table \ref{table_ess_candidate},  the attractor $\up^* = (\nvdf,0,\frac{b-d}{\varrho^*})$ and  $q^*=0$. Thus $h(\nvdf,0) =h_i(\zt,\bnu)>0 $, and hence ${\cal B} (\varpi ) = \{0\}$.
Further, under the mutation policy $\varpi_\epsilon(q)$, from lemma \ref{lem_cont_of_attract}, $h(\theta^*_\epsilon,\psi^*_\epsilon)>0$, hence ${\cal B} (\varpi_\epsilon(q)  )  =\{0\}$. Thus  $\varpi $ is an ESSS.

When $ h_v(\zt;\bnu)<0$ and $\nu_e>b\rho - \frac{\nu_b}{\nvdf}$, consider $\varpi$ with $\beta > \frac{1}{ \ps{e}}$, from table \ref{table_ess_candidate}, $\up^* = (0,\ps{e},\eta^*)$ and $q^*=1$. Thus, the best response set ${\cal B} (\varpi ) = \{1\}$ as $h(0,\ps{e}) = h_v(\zt;\bnu)<0$. Further, from lemma \ref{lem_cont_of_attract}, $h(\theta^*_\epsilon,\psi^*_\epsilon)<0$, hence ${\cal B} (\varpi_\epsilon(q)  )  =\{1\}$. Thus  $\varpi $ is an ESSS.

Finally when $\nu_e< b\rho - \frac{\nu_b}{\nvdf}$ and $h_v^o(\zt;\bnu)<0$, consider $\varpi$ with $\beta > \frac{1}{\ps{o} }$. From table \ref{table_ess_candidate}, $\up^* = (\tht{o},\ps{o},\eta^*)$ and thus $q^*=1$. Thus    the best response set ${\cal B} (\varpi ) = \{1\}$ as $h(\tht{o},\ps{o}) = h_v^o(\zt;\bnu)<0$. 
From lemma \ref{lem_cont_of_attract}, $h(\theta^*_\epsilon,\psi^*_\epsilon)<0$, hence ${\cal B} (\varpi_\epsilon(q)  )  =\{1\}$. Thus  $\varpi $ is an ESSS.

If none of the conditions are satisfied, then there is no  $\varpi$ such that $\{q^*_{\varpi} \} = {\cal B} (\varpi)$, and hence no ESSS. \eop


\section{Proof of influencers game}

\label{appen_a}
{\bf Proof of lemma \ref{lem_vac_val}:} The proof follows using  backward induction. At $T$ from \eqref{eqn_dp} and \eqref{eqn_ter_rew}, if  state $X^j_T = x=(\fv,z,c)$, then the value function $u_T(x ; \pij )= c= \Gamma_T(c)$ trivially for any $\pij$. At $T-1$, the value function of a vaccinated agent in state $x = (\fv,z,c)$ from \eqref{eqn_dp_2}  is, 

\vspace{-3mm}
{\small\begin{eqnarray*}
    u_{T-1}(x;\pij) &=& 0 + E_x[u_{T}(X^j_{T};\pij)] = E[ C_{\sd, T}\mid C_{\sd, T-1}=c] = \Gamma_{T-1}(c).\\
 \mbox{ Also, from \eqref{eqn_ck},} && \Gamma_{T-1}(c)   = E\left[c+ \frac{1}{T}(\bar{\xi}_T-c)\right] = \frac{T-1}{T} c + \frac{1}{T}\mvs.
\end{eqnarray*}}%
Assume {\small$ u_{t}(x;\pij) = \Gamma_t(c) = \frac{t}{T}c+ \frac{T-t}{T} \mvs  $}, for $t\ge n+1$. For  $t=n$, from \eqref{eqn_dp_2},

\vspace{-3mm}
{\small
\begin{eqnarray*}
  u_{n}(x;\pij) &=&  E_x[u_{n+1}(X^j_{n+1};\pij)]=   E_x[ \Gamma_{n+1}(C_{\sd,n+1})] =  E_x[C_{\sd,T}] =  \Gamma_n(c).
  \end{eqnarray*}}%
From above, \eqref{eqn_ck} and induction hypothesis,

\vspace{-3mm}
{\small  \begin{eqnarray*}
 \Gamma_n(c) &=& E\left[\Gamma_{n+1}\left(c+ \frac{1}{n+1}(\bar{\xi}_{n+1}-c)\right)\right], \\
 &=&E\left[\frac{n+1}{T}\left(c+ \frac{1}{n+1}(\bar{\xi}_{n+1}-c)\right) + \frac{T-n-1}{T} \mvs\right]=\frac{n}{T}c+ \frac{T-n}{T} \mvs.
\end{eqnarray*}}Hence, the lemma is true for all $t$. \eop

\vspace{2mm}

\noindent{\bf Proof of lemma \ref{lem_opt_fbar}:} The proof follows using backward induction. At stage $T$,   from \eqref{eqn_dp} and \eqref{eqn_ter_rew}, $u_T(x;\pij)=0$ if $x=(\fs,z,c)$ and  $u_T(x;\pij)=c$ if $x=(\fv,z,c)$. At stage $T-1$ and  $x=(\fs,z,c)$, from \eqref{eqn_dp_2} and lemma \ref{lem_vac_val}, 
 $u_{T-1}(x;\pij) = \inf_{p \in [0,1]} p( C_v + \Gamma_{T-1}(c)), $ 
since $g_z = 0$ if $z \ge \baz$. Thus the $(T-1)$-stage-BR against any $\pij$ equals $\{0\}$ and $u_{T-1}(x,\pij)=0$. Assume this is true for $t \ge n+1$ and consider $t=n$,  $x=(\fs,z,c)$. Again from \eqref{eqn_dp} and \eqref{eqn_ter_rew}, 
$u_{n}(x;\pij) = \inf_{p \in [0,1]} p (C_v+ \Gamma_n(c)). $%
Thus, the $n$-stage-BR equals $\{0\}$, and $u_{n}(x,\pij)=0$. The result follows.\eop

\vspace{2mm}

\noindent \textbf{Proof of theorem \ref{thm_eqlbm}:}
We first prove the if part   using backward induction. We begin with arguments  for $t=T-1$ which also lead to the required induction hypothesis. Towards this, let $\pi^*=(d_1^*,\dots,d_{T-1}^*)$ be a special strategy constructed using Algorithm \ref{NE_policy_algorithm}. Let all the opponents use strategy  $\pi^*$; denote such opponent strategy profile by  $\bpi^{-j*}$. Any best response strategy of agent $j$ against $\bpi^{-j*}$ is given by dynamic programming equations \eqref{eqn_dp}-\eqref{eqn_dp_2}. For $(\pi^*,\dots,\pi^*)$ to be equilibrium, we need to show that $\pi^*$ is a best response strategy against $\bpi^{-j*}$, i.e., for every $t,x$,  $d_t^{*}(x,\fv)$ is in $t$-stage-BR set against $\bpi^{-j*}$, i.e., it achieves the infimum in \eqref{eqn_dp_2}. In this proof, we briefly represent $E_{x,p}[u_{t+1}(X^j_{t+1};\bpi^{-j*})]$  as $E^{p,-j*}[u_{t+1}(X^j_{t+1})]$.

At stage $T-1$ and state $x=(\fs,z,c)$, from \eqref{eqn_ter_rew}, \eqref{eqn_rew}, \eqref{eqn_dp}-\eqref{eqn_dp_2} and by lemma \ref{lem_vac_val}, the DP equation \eqref{eqn_dp_2} simplifies to:

\vspace{-3mm}
{\small
\begin{eqnarray}\label{eqn_dp_t_1}
 u_{T-1}(\cs;\bpi^{-j*})&=& \inf_{p \in [0,1]} \{p(C_v - g_z)+E^{p,-j*}[r_T(X^j_{T})]\}, \nonumber\\
&& \hspace{-2.5cm} = \inf_{p \in [0,1]} \{p(C_v - g_z) + p \Gamma_{T-1}(c) + (1-p) E^{-j*}[\CI \indc{z + Y_{T} <\sbaz}]\}, \\    
&& \hspace{-2.5cm} = \inf_{p \in [0,1]} \{p(C_v - g_z + \Gamma_{T-1}(c) -E^{-j*}[\CI \indc{ Y_{T} <\sbaz - z}])+ E^{-j*}[\CI \indc{ Y_{T} <\sbaz - z}] \}. \nonumber     
 \end{eqnarray}}%
If $z\ge \baz$, the only choice for $d^*_{T-1}(x,\fv)$ in the special strategy (see      \eqref{eqn_cal_P_T_minus_1}, \eqref{eqn_cal_P_T_minus_1_for_M}) is 0. From lemma \ref{lem_opt_fbar} the optimizer in \eqref{eqn_dp_t_1}  equals $0$ and matches with $d^*_{T-1}(x,\fv)$.

\underline{Case 1, $t= T-1$, $x$ with $z<\baz<M$:}
 If further  $C_v + \Gamma_{T-1}(c)- g_z\ge \CI$, from \eqref{eqn_cal_P_T_minus_1}, $d^*_{T-1}(x,\fv)=0$ (in the special strategy $\pi^*$); thus none of the opponents vaccinate, hence $Y_T= 0$ a.s., and thus  $E^{-j*}[\CI \indc{Y_T<{\sbaz}-z}] = C_i$.  Then the 
$(T-1)$-stage-BR in  \eqref{eqn_dp_t_1} equals $\{0\}$  if  $C_v + \Gamma_{T-1}(c)- g_z> \CI$  and $[0,1]$ if $ C_v + \Gamma_{T-1}(c)- g_z= \CI$.  Both the sets  contain $0 = d^*_{T-1}(x,\fv) $.

On the other hand if $C_v + \Gamma_{T-1}(c)- g_z \le 0$, from \eqref{eqn_cal_P_T_minus_1} $d^*_{T-1}(x,\fv)=1$. In this case $z+Y_T=M-1 \ge \baz$ a.s., thus  $E^{-j*}[\CI \indc{Y_T<\sbaz-z}] =0$. Thus the $(T-1)$-stage-BR in \eqref{eqn_dp_t_1} equals $\{1\}$  if  $C_v + \Gamma_{T-1}(c)- g_z<0$  and $[0,1]$ if $ C_v + \Gamma_{T-1}(c)- g_z= 0$, both of which contain $1$. 

Lastly, if   $0<C_v + \Gamma_{T-1}(c)- g_z < \CI$,  then from \eqref{eqn_cal_P_T_minus_1} $d_{T-1}^*(x,\fv)=\tp$, which satisfies (see \eqref{eqn_tilde_p}),
$$C_v - g_z + \Gamma_{T-1}(c)  =E^{-j*}[\CI \indc{ Y_{T} <\sbaz - z}] =E^{\tp}_x[\CI \indc{Y_T < \sbaz-z}].$$%
Hence, $(T-1)$-stage-BR  in \eqref{eqn_dp_t_1} equals $[0,1]$ and contains $\tp$. 

\underline{Case 2, $t= T-1$, $x$ with $z<\baz=M$:} In this case $z + Y_{T} \le M-1 < \baz$ a.s., thus $E^{-j*}[\CI \indc{ Y_{T} <\sbaz - z}] = \CI$ . Now say $C_v + \Gamma_{T-1}(c)- g_z >\CI$ then from \eqref{eqn_cal_P_T_minus_1_for_M}  $d^*_{T-1}(x,\fv)=0$. From \eqref{eqn_dp_t_1}, the $(T-1)$-stage-BR also equals $\{0\}$. Secondly if  $C_v + \Gamma_{T-1}(c)- g_z <\CI$ then from \eqref{eqn_cal_P_T_minus_1_for_M}, 
 $d^*_{T-1}(x,\fv)=1$ and then the $(T-1)$-stage-BR equals $\{1\}$. Lastly, consider $C_v + \Gamma_{T-1}(c)- g_z = \CI$, then in special strategy $\pi^*$, $d^*_{T-1}(x,\fv)=p\in [0,1]$  from \eqref{eqn_cal_P_T_minus_1_for_M}.  
In this case,  $[0,1]$  is 
 the $(T-1)$-stage-BR, and contains any such $p$.

In all,  the best response against $\bpi^{-j*}$ contains a strategy in which decision rule $d^{*}_{T-1}$ 
 for each state $x$ is included in the corresponding $(T-1)$-stage-BR. Also,  by using
$d^*_{T-1}(x,\fv)$ in \eqref{eqn_dp_t_1}, it is easy to verify that $u_{T-1}(\cs;\bpi^{-j*}) = v_{T-1}(\cs)$  defined in \eqref{eqn_u_t_minus_1} for all the  sub-cases, i.e., for all  $x=(\fs,z,c)$. Further, one can trivially upper bound using \eqref{eqn_u_t_minus_1}: $
 u_{T-1}(\cs;\bpi^{-j*})   \le C_v+\Gamma_{T-1}(c)- g_z.
$

We now aim to prove exactly similar things for all $t < T-1$ using backward induction. 
 To be precise, assume for all 
 stages $t \ge k+1$ that $d^*_t$ is included in $t$-stage-BR against $\bpi^{-j*}$ and that (recall $v_t(x)$ is defined in \eqref{eqn_u_t_minus_1} and \eqref{eqn_vk}),
\begin{eqnarray}\label{eqn_uk_vk}
    u_t(x,\bpi^{-j*})=v_t(x) \le C_v + \Gamma_t(c) - g_z, \mbox{ for any } x=(\fs,z,c).
\end{eqnarray}Consider stage $t=k$. From \eqref{eqn_sk} of lemma \ref{lem_vac_val}, the DP equation \eqref{eqn_dp_2} corresponding to $k$-stage-BR  simplifies to:
{\small
\begin{eqnarray}
u_{k} (\cs;\bpi^{-j*})  &=&\inf_{p \in [0,1]} \{ p( C_v - g_z ) + E_x^{p,-j*}[u_{k+1}(X^j_{k+1})]\},\nonumber\\
\label{eqn_dp_k}
 &&\hspace{-1cm} = \inf_{p \in [0,1]} \{ p( C_v - g_z +\Gamma_k(c) ) + (1-p)E_x^{-j*}[u_{k+1}(\fs,z+Y_{k+1},C_{\sd,k+1})]\}, \hspace{2mm} \nonumber\\
  &&\hspace{-1cm} = \inf_{p \in [0,1]} p\Big( C_v - g_z +\Gamma_k(c) -E_x^{-j*}[u_{k+1}(\fs,z+Y_{k+1},C_{\sd,k+1}) ]\Big) \\
  &&  \hspace{4cm}+ E_x^{-j*}[u_{k+1}(\fs,z+Y_{k+1},C_{\sd,k+1})], \nonumber
\end{eqnarray}}

If $z\ge \baz$, from \eqref{eqn_thm_br_1} and \eqref{eqn_thm_br_1_zM}, $d^*_k(x,\fv)=0$, and from  lemma \ref{lem_opt_fbar}, the $k$-stage-BR in \eqref{eqn_dp_k} also equals $\{0\}$. Now we consider the remaining cases.

\underline{Case 3, $t= k$, $x$ with $z<\baz<M$:}
If $ C_v - g_z +\Gamma_k(c) \le 0$ , then from \eqref{eqn_thm_br_1}  $d_k^* (x,\fv) \in \{0, 1\}  \cup  \cp_{(k)}$ in the special policy.
Say $d^*_k(x,\fv)=0$ then $Y_{k+1}=0$ a.s., in \eqref{eqn_dp_k}; further, from the induction hypothesis,
\begin{eqnarray*}
 u_{k+1}(\fs,z,C_{\sd,k+1};\bpi^{-j*}) &\le& C_v+\Gamma_{k+1}(C_{\sd,k+1})- g_z \ \mbox{ a.s.},\\
 \mbox{ and hence, } E_x^{-j*}[u_{k+1}(\fs,z,C_{\sd,k+1})]  &\stackrel{a}{\le}& C_v+\Gamma_{k}(c)- g_z,
\end{eqnarray*}where inequality $a$ follows from lemma \ref{lem_vac_val}. Hence, the $k$-stage-BR  in \eqref{eqn_dp_k} equals $\{0\}$ when inequality $a$ is strict and $[0,1]$ otherwise, both of which contain $0$.

Now say $d_k^*(x,\fv)=1$, then $z+Y_{k+1}=M-1 \ge \baz$ a.s., and so $E^{-j*}[u_{k+1}(\fs,z+Y_{k+1},C_{\sd,k+1})] = 0$ from \Cref{lem_opt_fbar}. Thus, the $k$-stage-BR is the solution of,

\vspace{-3mm}
{\small
\begin{eqnarray*}
u^*_{k} (\cs;\bpi^{-j*}) &=&\inf_{p \in [0,1]}p( C_v - g_z+ \Gamma_k(c)),
\end{eqnarray*}}and equals $\{1\}$ if $C_v - g_z+ \Gamma_k(c)<0$ and $[0,1]$ if $C_v - g_z+ \Gamma_k(c)=0$; both of which contain $1$. 

Lastly, if $d^*_k(x,\fv) = \tp \in \cp_{(k)}$, then from \eqref{eqn_thm_br_1} and induction hypothesis \eqref{eqn_uk_vk}, such $\tp$ satisfies,

\vspace{-3mm}
{\small
$$C_v + \Gamma_{k}(c)- g_z= E^{\tp}_{x}[v_{k+1}(\fs,z+{Y}_{k+1},C_{\sd,k+1})]=E^{-j^*}_{x}[u_{k+1}(\fs,z+{Y}_{k+1},C_{\sd,k+1})] .$$}%
Hence, the $k$-stage-BR in \eqref{eqn_dp_k} equals $[0,1]$ and contains $\tp$. Similar arguments follow when $C_v - g_z +\Gamma_k(c) > 0$.

\underline{Case 4, $t= k$, $x$ with $z<\baz=M$:}
If $C_v + \Gamma_{k}(c)- g_z\le E_x[v_{k+1}(\Bar{X})]$ where $\Bar{X}=(\fs,M-1,C_{\sd,k+1})$, then from \eqref{eqn_thm_br_1_zM}  $d_k^* (x,\fv) \in \{0, 1\}  \cup \{p\}$ in the special policy, where $p$ satisfies \eqref{eqn_p_for_bark_M}.
Say $d^*_k(x,\fv)=0$ then using similar arguments as in case 3, the $k$-stage-BR contains $0$.

On the other hand, if  $d^*_k(x,\fv)=1$, then $z+Y_{k+1}=M-1$ a.s. Then from induction hypothesis,
$$E_x^{-j*}[u_{k+1}(\fs,z+Y_{k+1},C_{\sd,k+1})] =  E_x[v_{k+1}(\Bar{X})].$$
Thus, $k$-stage-BR  in \eqref{eqn_dp_k}  equals
 $\{1\}$ if $C_v - g_z+ \Gamma_k(c)<E_x[ v_{k+1} (\bar{X})]$ and $[0,1]$ if $C_v - g_z+ \Gamma_k(c)=E_x[ v_{k+1} (\bar{X})]$, and contains $1$.

Lastly,  if $d^*_k(x,\fv)=\tp\in(0,1)$, then from   \eqref{eqn_p_for_bark_M}, $k$-stage-BR in \eqref{eqn_dp_k}, equals $[0,1]$ and contains $\tp$. 

Also,  by using
$d^*_{k}(x,\fv)$ in \eqref{eqn_dp_k}, it is easy to verify that,
\begin{eqnarray*}
 u_{k}(\cs;\bpi^{-j*}) = \left\{\begin{array}{ll}
  E_x[u_{k+1}((\fs,z,C_{\sd,{k+1}});\bpi^{-j*})]   & \mbox{ if } d_k^*(x,\fv)=0,  \\
  C_v + \Gamma_{k}(c)- g_z   & \mbox{ else. }
\end{array}\right.
\end{eqnarray*}

From \eqref{eqn_vk},  $ u_{k}(\cs;\bpi^{-j*}) = v_{k}(\cs)  \le C_v+\Gamma_{T-1}(c)- g_z
$, which proves the induction hypothesis.
Hence, there exists a best response strategy $\pi^{j*}=(d^{j*}_1,\dots,d^{j*}_{T-1})$, such that $d^{j*}_{k}=d^{*}_{k}$ for any $k$. 
Hence, $(\pi^*, \dots, \pi^*)$ is a symmetric Nash equilibrium. 

 Now say $(\pi,\dots,\pi)$ is a symmetric Nash equilibrium, we prove that $\pi$ is a special policy. Denote opponent strategy profile by  $\pij$ when all the opponents use strategy  $\pi$. By definition of Nash equilibrium, the best response  of agent $j$ against $\pij$ must contain $\pi$, and should solve the dynamic programming equations \eqref{eqn_dp}-\eqref{eqn_dp_2}, i.e., for every $t,x$,  $d_t^{*}(x,\fv)$ should be contained in $t$-stage-BR set against $\pij$. Say $\pi =(d_1,\dots,d_{T-1})$.

The proof again follows from backward induction. Towards this, let  $t = T-1$ and $x=(\fs,z,c)$. From lemma \ref{lem_opt_fbar}, the $(T-1)$-stage-BR against  $\pij$ equals $\{0\}$ if $z\ge\baz$. Hence, we must have $d_{T-1}(x,\fv) = 0$ which equals $d_{T-1}^*(x,\fv)$ of algorithm \ref{NE_policy_algorithm}. Now we consider remaining the cases: 

\underline{Case 1, $t= T-1$, $x$ with $z<\baz<M$:}
Consider $x$ such that $C_v + \Gamma_{T-1}(c)- g_z\ge \CI$. In this case, if $d_{T-1}(x,\fv) \in (0,1]$ in the equilibrium policy, from \eqref{eqn_dp_t_1} the $(T-1)$-stage-BR equals $\{0\}$ which  contradicts $(\pi,\dots,\pi)$ being a symmetric NE. 
Hence $d_{T-1}(x,\fv)$ must  equal  0, and matches with \eqref{eqn_cal_P_T_minus_1} and \eqref{eqn_pi_T_minus_1}. 

Now consider $x$ such that $C_v + \Gamma_{T-1}(c)- g_z\le 0$. In this case, if $d_{T-1}(x,\fv) \in [0,1)$ in the equilibrium policy, from \eqref{eqn_dp_t_1} the $(T-1)$-stage-BR equals $\{1\}$ which again contradicts symmetric NE condition.
Hence $d_{T-1}(x,\fv)$ must  equal 1, and matches with \eqref{eqn_cal_P_T_minus_1} and \eqref{eqn_pi_T_minus_1}. 

Lastly, consider $0<C_v + \Gamma_{T-1}(c)- g_z<\CI$. In this case, if $d_{T-1}(x,\fv) = 1$, in the equilibrium strategy then $(T-1)$-stage-BR equals $\{0\}$ and if  $d_{T-1}(x,\fv) = 0$,  then $(T-1)$-stage-BR equals $\{1\}$, hence a contradiction. Now we are left with
$d_{T-1}(x,\fv) = p \in (0,1)$.  If $p$ is such that  $C_v - g_z + \Gamma_{T-1}(c) <  E^{p}_x[\CI \indc{Y_T < \baz-z}]$, then from \eqref{eqn_dp_t_1} the $(T-1)$-stage-BR equals $\{1\}$. On the other hand,   if  $C_v - g_z + \Gamma_{T-1}(c) >  E^{p}_x[\CI \indc{Y_T < \baz-z}]$, then $k$-stage-BR equals $\{0\}$, both of which lead to a contradiction. Hence the only choice for $p$ is the one  that satisfies equality, $C_v - g_z + \Gamma_{T-1}(c) =  E^{p}_x[\CI \indc{Y_T < \baz-z}]$, and such a choice matches with that given by algorithm \ref{NE_policy_algorithm} or \eqref{eqn_cal_P_T_minus_1} and \eqref{eqn_pi_T_minus_1}.

\underline{Case 2, $t= T-1$, $x$ with $z<\baz=M$:}
Consider $x$ such that $C_v + \Gamma_{T-1}(c)- g_z > \CI$. In this case, for any $d_{T-1}(x,\fv) \in [0,1]$, from \eqref{eqn_dp_t_1} the $(T-1)$-stage-BR equals $\{0\}$. 
Hence for $d_{T-1}(x,\fv)$ to be a part of symmetric NE, it must  equal  0, and matches with 
 \eqref{eqn_pi_T_minus_1} and \eqref{eqn_cal_P_T_minus_1_for_M}. 

In similar lines, if $C_v + \Gamma_{T-1}(c)- g_z < 0$, $d_{T-1}(x,\fv)$ must equal 1 as in  \eqref{eqn_pi_T_minus_1} and \eqref{eqn_cal_P_T_minus_1_for_M}. 

Lastly, if $C_v + \Gamma_{T-1}(c)- g_z=\CI$,   \eqref{eqn_pi_T_minus_1} and \eqref{eqn_cal_P_T_minus_1_for_M} prescribe any $p\in[0,1]$.

Hence, $d_{T-1}$ matches with that prescribed by algorithm \ref{NE_policy_algorithm} for all $x$. Also from part (i), $u_{T-1}(x;\pij)=v_{T-1}(x)$.

Assume that $d_{t}(x,\fv)$ of $\pi$ matches with $d^*_t(x,\fv)$ of algorithm \ref{NE_policy_algorithm} for $t = k+1,\dots,T-1$ and that,
$u_{k+1}((\fs, z,  c);\pij)=v_{k+1}(\fs, z, c) .$
Consider stage $t=k$.  From lemma \ref{lem_opt_fbar}, the $k$-stage-BR against  $\pij$ equals $\{0\}$ if $z\ge\baz$, and hence $d_k(x,
\fv)=0$ as that prescribed by Algorithm \ref{NE_policy_algorithm}.

\underline{Case 3, $t= k$, $x$ with $z<\baz<M$:} Consider $x$ such that $C_v + \Gamma_{T-1}(c)- g_z\le 0$ and  say $d_{k}(x,\fv) = q$.  If\footnote{if the decision $d_{k}(x,\sfv) \not \in \{0,1,p$\} where $p$ is unique and satisfies 
 the equality in $\cp_{(k)}$. In other words, $q \in (0,1) - \{p\}$.} $q\not \in   \{0, 1\}  \cup  \cp_{(k)}$ of \eqref{eqn_thm_br_1} and \eqref{eqn_vk} in $\pi$ then,

 \vspace{-3mm}
 {\small
\begin{equation}   \label{eqn_logic_case3}
C_v + \Gamma_{k}(c)- g_z \ne E^{-j}_{x}[u_{k+1}(\fs, z + Y_{k+1}, C_{\sd,k+1})]
= E^{q}_{x}[u_{k+1}(\fs, z + Y_{k+1}, C_{\sd,k+1})],\end{equation}}%
 as $v_{k+1}=u_{k+1}$ by induction hypothesis.
 From \eqref{eqn_dp_k} the $k$-stage-BR equals $\{0\}$ if the left hand side is bigger in the above, and  $\{1\}$ else. Hence a contradiction to $\pi$ being a 
 part of symmetric NE. 
Hence $d_{k}(x,\fv) \in  \{0, 1\}  \cup  \cp_{(k)} $ and matches with \eqref{eqn_thm_br_1} and \eqref{eqn_vk}. 

If $C_v + \Gamma_{T-1}(c)- g_z 
 > 0$, using similar logic
$d_{k}(x,\fv) \in  \{0\}  \cup  \cp_{(k)} $ as in \eqref{eqn_thm_br_1}.

\underline{Case 4, $t= k$, $x$ with $z<\baz=M$:} Now, consider $x$ such that $C_v + \Gamma_{T-1}(c)- g_z\le E_x[u_{k+1}(\fs,M-1,C_{\sd,k+1})]$. If $d_{k}(x,\fv) \not \in   \{0, 1\}  \cup  \{p\}$ of \eqref{eqn_vk} and 
\eqref{eqn_thm_br_1_zM} in $\pi$, then, equality in \eqref{eqn_p_for_bark_M} is not satisfied, and leads to contradiction as before.

Similar logic follows for remaining case  $C_v + \Gamma_{T-1}(c)- g_z 
 >  E_x[v_{k+1}(\fs,M-1,C_{\sd,k+1})]$. 
 

Further, upon substitution of $d_k(x,\fv)$ in \eqref{eqn_dp_k} for all $x$, it is easy to verify that,
\begin{eqnarray*}
 u_{k}(\cs;\bpi^{-j*}) = v_k(x)= \left\{\begin{array}{ll}
  E_x[u_{k+1}((\fs,z,C_{\sd,{k+1}});\bpi^{-j*})]   & \mbox{ if } d_k^*(x,\fv)=0,  \\
  C_v + \Gamma_{k}(c)- g_z   & \mbox{ else. }
\end{array}\right.
\end{eqnarray*}

Hence the induction hypothesis is true and the proof follows. 
\eop

%

\noindent\textbf{Proof of Corollary \ref{cor_wait_n_watch}:} Under wait-and-watch equilibrium, $d_t(x_t,\fv) \equiv 0$ for all $x_t$, $t<T-1$, therefore $X_t = (\fs,0,C_{\sd,t})$ for all $t\le T-1$.
Thus from theorem \ref{thm_eqlbm} and \eqref{eqn_vk}, the utility of a typical agent under this equilibrium, where initial state $x=(\fs,0,c_{\sd,1})$ equals (let $g:=g_0$):

\vspace{-3mm}
{\small
\begin{eqnarray}\label{eqn_wait_and_watch_util}
u_{1} (x;\bpi^{-j*}_w)&=&  v_{1}(\fs,0,c_{\sd,1}) = E_x[v_2(\fs,0,C_{\sd,2})], \nonumber \\
    &=& E_x[E_{X_2}[\dots E_{X_{T-2}}[v_{T-1}(\fs,0,C_{\sd,T-1})]\dots]],\\
 &=& E_x[v_{T-1}(\fs,0,C_{\sd,T-1})]= E_x[\min\{C_v- g + \Gamma_{T-1}(C_{\sd,T-1}), \CI\}], \nonumber
\end{eqnarray}}from \eqref{eqn_u_t_minus_1}.
Now consider any other symmetric NE, $\bpi = (\pi,\dots,\pi)$.  Define $B_t:= \{x : d_t(x,\fv) >0\}$,  from \eqref{eqn_vk} for any $x=(\fs,z,c)$ and $t \le T-2$:

\vspace{-3mm}
{\small
\begin{eqnarray}\label{eqn_v_term}
    v_t(x) = (C_v- g_z + \Gamma_t(c))\indc{x \in B_t} + E_{x}[v_{t+1}(\fs,z,C_{\sd,t+1})] \indc{x \in B^c_t}.
\end{eqnarray}}%
Thus from   theorem \ref{thm_eqlbm}, with $x=(\fs,0,C_{\sd,1})$, where $C_{\sd,1} = c_{\sd,1}$ a.s., we have,

\vspace{-3mm}
{\small
\begin{eqnarray*}
u_{1} (x;\bpi^{-j*}) &=& (C_v- g + \Gamma_1(C_{\sd,1}))\indc{x \in B_1} + E_{x}[v_{2}(\fs,0,C_{\sd,2})] \indc{x \in B^c_1}.
\end{eqnarray*}}%
Using \eqref{eqn_v_term} for the second term in the above, we get ($X_1$ represents $x$)

\vspace{-3mm}
{\small
\begin{eqnarray*}
u_{1} (x;\bpi^{-j*}) &&\\
 %
&& \hspace{-2cm} = (C_v- g + \Gamma_1(C_{\sd,1}))\indc{x \in B_1} + E_{x}[(C_v- g + \Gamma_2(C_{\sd,2}))\indc{X_2 \in B_2}\indc{x \in B^c_1}]\\
&& \hspace{2cm}  + E_{x}[ v_{3}(\fs,0,C_{\sd,3}) \indc{X_2 \in B^c_2} \indc{x \in B^c_1}],\\
  && \hspace{-2cm} \stackrel{a}{=} \sum_{\tau=1}^{T-1} E_x[(C_v - g+ \Gamma_\tau (C_{\sd,\tau}))  \indc{X_{\tau} \in B_\tau}\prod_{\ttau=1}^{\tau-1} \indc{X_{\ttau} \in B^c_\ttau}]\\
 && \hspace{2cm}  + E_{x}[v_{T-1}(\fs,0,C_{\sd,T-1})  \prod_{\ttau=1}^{T-1} \indc{X_{\ttau} \in B_\ttau^c}],
 \end{eqnarray*}}%
 where $a$ follows by repeated application of \eqref{eqn_v_term} and mathematical induction (also observe $X_\tau=(\fs,0,C_{\sd,\tau})$ on $\cap_{t < \tau}B_t^c$). Now using an obvious lower bound at $\tau = T-1$ in the first term, and using \eqref{eqn_u_t_minus_1} for the second term,

 \vspace{-3mm}
 {\small
 \begin{eqnarray} \label{eqn_appen_greater}
%
%
  && \hspace{-8mm}u_{1} (x;\bpi^{-j*}) 
\   \ge \ \sum_{\tau=1}^{T-2} E_x[(C_v - g+ \Gamma_\tau (C_{\sd,\tau}))  \indc{X_{\tau} \in B_\tau}\prod_{\ttau=1}^{\tau-1} \indc{X_{\ttau} \in B^c_\ttau}]\\
&& \hspace{3mm} + E_x[\min\{C_v - g + \Gamma_{T-1}(C_{\sd,T-1}),\CI \}  \indc{X_{T-1} \in B_{T-1}} \prod_{\ttau=1}^{T-2} \indc{X_{\ttau} \in B^c_\ttau}]\nonumber\\
&& \hspace{3mm} + E_{x}[\min\{C_v - g + \Gamma_{T-1}(C_{\sd,T-1}),\CI \}]    \indc{X_{T-1} \in B^c_{T-1}}\prod_{\ttau=1}^{T-2} \indc{X_{\ttau} \in B_\ttau^c}], \nonumber\\
 &=&  \sum_{\tau=1}^{T-2} E_x[(C_v - g+ \Gamma_\tau (C_{\sd,\tau}))  \indc{X_{\tau} \in B_\tau}\prod_{\ttau=1}^{\tau-1} \indc{X_{\ttau} \in B^c_\ttau}]\nonumber\\
&& \hspace{1cm} + E_x[\min\{C_v - g + \Gamma_{T-1}(C_{\sd,T-1}),\CI \} \prod_{\ttau=1}^{T-2} \indc{X_{\ttau} \in B^c_\ttau}].\nonumber
%
%
\end{eqnarray}}%
With sigma-algebra ${\cal F}_{\tau
}:=\sigma(X_t: t \le \tau)$, and using lemma \ref{lem_vac_val}, we get:

\vspace{-3mm}
{\small
\begin{eqnarray*}
&&\hspace{-5mm}u_{1} (x;\bpi^{-j*}) 
\   = \ \sum_{\tau=1}^{T-3} E_x[(C_v - g+ \Gamma_\tau (C_{\sd,\tau}))  \indc{X_{\tau} \in B_\tau}\prod_{\ttau=1}^{\tau-1} \indc{X_{\ttau} \in B^c_\ttau}]\\
&& \hspace{5mm} + E_x[(C_v - g+ E_{{\cal F}_{T-2}}[ \Gamma_{T-1} (C_{\sd,{T-1}})])  \indc{X_{{T-2}} \in B_{T-2}}\prod_{\ttau=1}^{T-3} \indc{X_{\ttau} \in B^c_\ttau}] \\
&& \hspace{1cm}  + E_x[\min\{C_v - g + \Gamma_{T-1}(C_{\sd,T-1}),\CI \} \indc{X_{T-2} \in B^c_{T-2}} \prod_{\ttau=1}^{T-3} \indc{X_{\ttau} \in B^c_\ttau}],\\
%
%
%
&&{\ge}\sum_{\tau=1}^{T-3} E_x[(C_v - g+ \Gamma_\tau (C_{\sd,\tau}))  \indc{X_{\tau} \in B_\tau}\prod_{\ttau=1}^{\tau-1} \indc{X_{\ttau} \in B^c_\ttau}]\\
&& \hspace{1cm} + E_x[\min\{C_v - g + \Gamma_{T-1}(C_{\sd,T-1}),\CI \} \prod_{\ttau=1}^{T-3} \indc{X_{\ttau} \in B^c_\ttau}],
\end{eqnarray*}}%
where above follows using similar logic as in \eqref{eqn_appen_greater}, along with the Tower property of conditional expectation (see \cite{protter}). Continuing in a similar manner, we get,

\vspace{-3mm}
{\small
\begin{eqnarray*}
 u_{1} (x;\bpi^{-j*}) &\ge& E_x[(C_v - g+ \Gamma_1(C_{\sd,1}))  \indc{X_1 \in B_1}\\
 && \hspace{1cm}+ E_x[\min\{C_v - g + \Gamma_{T-1}(C_{\sd,T-1}),\CI \}  \indc{X_{1} \in B^c_1}],\\
&\ge& E_x[\min\{C_v - g + \Gamma_{T-1}(C_{\sd,T-1}),\CI \}.
\end{eqnarray*}}%
Above lower bound is the utility under wait-and-watch \eqref{eqn_wait_and_watch_util}.
\eop

\noindent\textbf{Proof of lemma \ref{lem_monotonicity_of_p}:} We provide the proof for with respect to $g$, and the proof with respect to $c$ follows using exact similar arguments and the fact that $\Gamma_{T-1}(c)$ is linear and increasing in  $c$ (see \eqref{eqn_sk}).

When $\baz=M$, then clearly $p(g,c)$ is non-decreasing in $g$ from \eqref{eqn_outcome_vaccination}. 

Now, consider $\baz < M$ and see \eqref{eqn_outcome_vaccination}. Let $g_1>g_2$. If $g_1 \ge  C_v +  \Gamma_{T-1}(c)$ then $p(g_1,c)=1$, so $p(g_1,c) \ge p(g_2,c)$ trivially.   If $g_1 \le C_v + \Gamma_{T-1}(c) - \CI$, then $g_2 \le C_v + \Gamma_{T-1}(c) - \CI$, thus $p(g_1,c) =p(g_2,c)=0$. 
Define ${\cal G}:= (C_v+ \Gamma_{T-1}(c) - C_i, C_v+ \Gamma_{T-1}(c) )$. Then, we have proved $p(g,c)$ is non-decreasing on ${\cal G}^c$.

Next consider $p(g,c)$ on ${\cal G}$. Define the correspondence,   $D(g) = [0,1]$ for all $g$, and a function $H:[0,1] \times  {\cal G}   \to \mathbb{R}$  as below,

\vspace{-3mm}
{\small
$$H(p,g) := \left(F_{M-1}(\baz(\bnu)-1;p) -\frac{C_v + \Gamma_{T-1}(c)-g}{\CI} \right)^2.$$}%
Clearly, $H$ is a continuous function and $D$ is a constant compact-valued (hence also continuous) correspondence (see \cite{optimization}). Further $D^*(g):= \arg\min \{H(p,g): p \in D(g)\}=\{\tp\}$, by uniqueness of solution $\tp$ of \eqref{eqn_p_outcome_F}. 
Thus from the Maximum theorem (e.g., \cite{optimization}), $D^*$ is continuous single-valued correspondence on $ {\cal G} $, and hence $g\mapsto p(g,c)$ is continuous  on ${\cal G}$. Further, $p \mapsto F_M(\baz-1;p)$ is strictly decreasing and continuous, so continuous  inverse exists and is strictly decreasing. Thus, from \eqref{eqn_p_outcome_F},  $g\mapsto p(g,c)$  is strictly increasing on ${\cal G}$.

Now, we prove the continuity of $p(g,c)$ at boundary points of ${\cal G}$. Let $g_1 :=  C_v +\Gamma_{T-1}(c) - C_i$ and consider $g= g_1 + \epsilon$. From \eqref{eqn_p_outcome_F},
$F_{M-1}(\baz-1;p(g,c))= 1-\frac{\epsilon}{\CI}$.
Since $p \mapsto F_M(\baz-1;p)$ is strictly decreasing and continuous,  the inverse function exists and is continuous. Thus, as $\epsilon \to 0$, the inverse image  $p(g,c) \to 0$. Hence, we have continuity at $g_1$ further because $p(g,c)=p(g_1,c)=0$ when $g<g_1$. Similar arguments follow at $g=  C_v +\Gamma_{T-1}(c)$.

Hence, $p(g,c)$ is a non-decreasing function which is  continuous and strictly increasing on ${\cal G}$ when $\baz<M$.  \eop

\section{Proofs for  Leaders game}\label{appendix_leaders_game}

\newcommand{\Gmc}{\Gamma_C}
\textbf{Notations:} This appendix uses the following notations: $C_{\sd,T-1}$ is represented as $C$ and $\Gamma_{T-1}(C_{\sd,T-1})$ as $\Gmc$. Further, for any $\baz$, $p(g,C_{\sd,T-1})$ is represented as $P_g$, when we need to show dependence on $\baz$ explicitly, we show that using superscript, i.e., $P_g^{\sbaz}$. Further, $P_g(\omega)$, $\Gmc(\omega)$ denote the realisation of $P_g, \Gmc$ for any sample point $\omega$.

\begin{lemma} \label{lem_monotonicity_of_EF}
Define $\tg:=  \max\{C_v-C_i +\frac{c_{\sd,1}+\mvs}{T},0\}$. Assume {\bf A.}2, then the function $g\mapsto \Np(g,\baz)$ is  continuous and strictly decreasing for $g > \tg$; when $\baz<M$,  it is  continuous for all $g \ge 0$.  Further, when $\tg>0$, then $\Np(g,\baz)=1$ for  $g \in [0,\tg]$. 
\end{lemma}

\noindent\textbf{Proof:} To begin with, first observe by \eqref{eqn_ck}, \textbf{A.}2, and lemma \ref{lem_vac_val} that 
 $\Gmc$ is distributed as $p_{\Gamma} \dirac{\bar \Gamma} (\cdot) + (1-p_\Gamma)f_{\Gamma} (\cdot)dx$, with  $\dirac{\bar \Gamma}$ a dirac measure at  ${\bar \Gamma} := \frac{c_{\sd,1}+\mvs}{T}$, $p_\Gamma = (p_0)^{T-2}$ and  $f_\Gamma$  is the density  supported on $[\bar{\Gamma},\infty)$ that can be computed using $f_{\bar\xi}$ of assumption {\bf A}.2.

Now consider $\baz=M$. From \eqref{eqn_outcome_vaccination}  $\Np(g,\baz)  =   P(\Gmc \ge \CI + g - C_v  )$. Thus by the density component of $\Gamma_C$, $\Np(g,\baz)$ is continuous and strictly decreasing for all $g>\tg$. Further with $\tg>0$  clearly $\Np(g,\bnu)=1$ for all $g \in [0,\tg]$ (and there can be a discontinuity at $\tg$).

Consider $\baz<M$. From lemma \ref{lem_monotonicity_of_p}, $g \mapsto P_g$ is continuous point-wise, and hence $g \mapsto \Np(g,\baz)$ is also continuous by Lebesgue’s dominated convergence theorem (DCT) (e.g., see \cite{protter}).   

When $\tg > 0$, take any $g \in [0,\tg]$. From \eqref{eqn_outcome_vaccination} and \textbf{A.}2, $P_g = 0$ a.s., thus $\Np(g,\bnu)=1$. Now, for strict monotonicity, consider   $g_1>g_2\ge \tg$, let $ \varepsilon = g_1 - g_2$, and define,
 $${A}:=\{\omega:\CI+ g_2 - C_v<\Gmc(\omega)<\CI + g_1 - C_v - \Tilde{\varepsilon}\}.$$
In the above $\Tilde{\varepsilon}<\varepsilon$ is chosen such that $P_{g_1}(\omega) > \varepsilon_2$  for all $\omega\in A$ for some $\varepsilon_2>0$, which   exists by  lemma \ref{lem_monotonicity_of_p} and \eqref{eqn_outcome_vaccination}; in fact, $\varepsilon_2 =P_{g_1} (\Bar{\omega})$ where  $\Bar{\omega}$ is such that $\Gmc(\Bar{\omega})
 = \CI + g_1 - C_v-\Tilde{\varepsilon}$. Further, from the same lemma $P_{g_2}(\omega)=0$  for all $\omega\in A$. Since   $F_M(\baz-1;p)$ is strictly decreasing in $p$ for $p\in (0,1)$,  we have $F_M(\baz-1;P_{g_1}(\omega))<1-\varepsilon_3$  for some $\varepsilon_3>0$ (dependent on $\varepsilon_2$) for all $\omega \in A$.
Hence $F_M(\baz-1;P_{g_1}(\omega)) + \varepsilon_3 <1= F_M(\baz-1;P_{g_2}(\omega))$  for any $\omega \in A$. Further, again by lemma \ref{lem_monotonicity_of_p},   $F_M(\baz-1;P_{g_1}(\omega))\le F_M(\baz-1;P_{g_2}(\omega))$ for all $\omega$. By {\bf A.}2, $P(A)>0$ and hence $\Np(g_1,\baz)<\Np(g_2,\baz)$ when $g_1>g_2$ (see \eqref{eqn_lead_util_cons}). \eop

\textbf{Proof of theorem \ref{thm_lead_optimal}:} 
 From lemma \ref{lem_monotonicity_of_p}, the objective function $\Ui$ in  \eqref{eqn_lead_util_cons}-\eqref{eqn_leader_opt_wait} is strictly increasing in $g$.  Thus if $\Np(0,\baz)\le\delta$, then $g^*=0$ is the optimizer.

 Now say $\Np(0,\baz)>\delta$. 
When $\tg>0$, by lemma \ref{lem_monotonicity_of_EF}, any $g \in [0,\tg]$ is infeasible  as  $\Np(g,\baz) = 1$. Thus in all, any  $g \le \tg$ is infeasible. Further again by lemma \ref{lem_monotonicity_of_EF},  $\Np(g,\baz)$ is continuous and \underline{strictly decreasing} for all $g > \tg$. Also, 
$\Np(g,\baz) \downarrow 0$ as $g \to \infty$ by DCT, because $F_M(\baz-1,P_{g}(\omega)) \downarrow 0$ for almost all $\omega$ (see \eqref{eqn_outcome_vaccination}). Thus by intermediate value theorem, there exist a unique $g^*$ satisfying, $\Np(g^*,\baz)=\delta$, and by strict monotonicity of  $\Ui$, $g^*$ is the optimizer. 
\eop





\subsection{Optimizers for $\perf{c}$}
\label{sec_optimizers}
For $\perf{c}$, the  objective function and the 
 constraint in   \eqref{eqn_leader_opt_wait} simplify to {\small$\Ui(g,\baz)=Mgp(g,c)$} and {\small $\Np(g,\baz) = F_M\big(\baz-1;p(g,c) \big )$},
where $p(g,c)$ is defined in  \eqref{eqn_outcome_vaccination} with $\Gamma_{T-1}(c)=c:=\Gamma$. Again, $\Ui(g,\baz)$ is strictly increasing whereas $\Np(g,\baz)$ is decreasing  in  $g$. Thus, when $\Np(0,\baz) \le \delta$, $g^*=0$ is the optimizer. 
Now consider $\Np(0,\baz)>\delta$.  

Say  $\baz<M$. If $C_v + \Gamma - \CI \ge 0$, from \eqref{eqn_outcome_vaccination},  $p(g,c)=0$ for all $g \in [0,C_v + \Gamma - \CI]$,  and thus $\Np(g,\baz)=1$ making such $g$ infeasible. Further, $p(g,c)=1$ and thus $\Np(g,\baz)=0$ for all $g  \ge C_v + \Gamma$; furthermore $p(g,c) \not\in \{0,1\}$ for all $g \in ( C_v +\Gamma - \CI ,  C_v +\Gamma)$.  Hence using lemma \ref{lem_monotonicity_of_p},  $g\to \Np(g,\baz)$ is continuous and strictly decreasing on the interval  $(C_v +\Gamma - \CI,C_v+\Gamma)$. Thus by intermediate value theorem, there exists a unique $g^* \in (C_v +\Gamma - \CI,C_v+\Gamma)$ such that $\Np(g^*,\baz)=\delta$ for any $\delta \in (0,1)$, which is the optimizer further because    $\Ui$ is strict monotone in $g$.  If $C_v + \Gamma - \CI <0$, then using similar arguments,  there exists a unique $g^* \in (0,C_v+\Gamma)$ that satisfies $\Np(g^*,\baz)=\delta$ and  is the optimizer. As $\Np(g^*,\baz)=F_M(\baz-1,p(g^*,c))$, we get $p(g^*,c)=p^* $, where $p^* \in (0,1)$ is the 
 unique solution of  $F_M(\baz-1,p^*)=\delta$. Thus  from \eqref{eqn_outcome_vaccination}-\eqref{eqn_p_outcome_F}, $g^*= C_v + \Gamma - C_i F_{M-1}(\baz-1;p^*)$ which proves \eqref{eqn_prf_info}. 

One can easily derive the closed form expressions when $\baz =1$ provided in \eqref{eqn_perf_opt} by solving \eqref{eqn_prf_info}.  When $\baz=M$, $p(g,c) \in \{0,1\}$ for any $g$ from \eqref{eqn_outcome_vaccination}.Thus for $\perf{c}$ any feasible $g$ satisfies $p(g,c) =1$, and thus feasible $g\in( C_v + \Gamma -\CI,\infty)$. 
Hence  we only have $\epsilon$-optimality for $\baz=M$, and $\epsilon$-optimizer $\og{M,\epsilon} = C_v + \Gamma -\CI+ \epsilon$ from the strict monotonicity of $U$. \eop

\vspace{-2mm}
\begin{lemma}
\label{lem_del_zero}\label{lem_pk_per_inf}
Let $p^*_k$ be as in \eqref{eqn_prf_info} and consider any $1 \le k< M$. Then    for any $\delta \in (0,1)$,   $p^*_k \in (0,1)$ and $p^*_k<p^*_{k+1}$. Further
  as $\delta\to 0$,  $\frac{\delta}{(1-p^*_1)(1-p^*_k)} \to 0$. 
  \vspace{-2mm}
\end{lemma}
\noindent\textbf{Proof:} 
By definition of $p^*_k, p^*_{k+1}$, we have $\delta =  F_M(k-1;p^*_k) < F_M(k;p^*_{k})$, and hence $p^*_{k} < p^*_{k+1}$, as $p^*_{k+1}$ also satisfies  $\delta = F_M(k;p^*_{k+1})$. This proves first statement using results of Section \ref{sec_optimizers}. 

From \eqref{eqn_prf_info} $p^*_1$ satisfies,
{\small $
1-p^*_1 = \delta^{\frac{1}{M}}
$},  and $p_k^*$ satisfies,
{\small $$ \delta = (1-p^*_k)^{M-k+1}  \sum_{i=0}^{k-1} \binom{M}{i} (p^*_k)^i (1-p^*_k)^{k-1-i}  \le (1-p^*_k)^{M-k+1} \ 2^{M}.$$}%
Further {\small $\delta^{1-\frac{1}{M}}  = \frac{\delta}{(1-p^*_1) }   $} and {\small$(1-p^*_k) \le (1-p^*_1) = \delta^{\frac{1}{M}}$}, and so one can sandwich:

\vspace{-3mm}
{\small
\begin{eqnarray*}
  \frac{\delta^{1-\frac{1}{M}}}{\delta^{\frac{1}{M}}}  &\le& \frac{\delta}{(1-p^*_1)(1-p^*_k)}\le \frac{ \delta^{1-\frac{1}{M}}}{\left (\delta 2^{-M} \right )^{\frac{1}{M-k+1}}}, \mbox{\normalsize and hence the results follows as, }\\
   \delta^{1-\frac{2}{M}}  &\le& \frac{\delta}{(1-p^*_1)(1-p^*_k)}\le \delta^{1-\frac{1}{M}-\frac{1}{M-k+1}} \  2^{\frac{M}{M-k+1}}. \hspace{35mm} \mbox{ \eop }
\end{eqnarray*}} 

\noindent\noindent\textbf{Proof of theorem \ref{thm_perf_info_var_z}}: To begin with for any $l<M$,  from  \eqref{eqn_prf_info}:

\vspace{-2mm}
{\footnotesize
\begin{eqnarray}\label{eqn_FM_FM_1}
       F_{M-1}(l-1;p^*_l) &=& \hspace{-2mm}\sum_{i=0}^{l-1} \binom{M-1}{i} (p^*_l)^i (1-p^*_l)^{M-1-i} 
    =  \frac{1}{1-p^*_l}\sum_{i=0}^{l-1} \binom{M-1}{i} (p^*_l)^i (1-p^*_l)^{M-i},\nonumber\\
    &=& \frac{1}{1-p^*_l} \left(F_M(l-1;p^*_l) + \sum_{i=0}^{l-1} \left( \binom{M-1}{i} -\binom{M}{i} \right )(p^*_l)^i (1-p^*_l)^{M-i}\right),\nonumber\\ 
        &=& \frac{\delta}{1-p^*_l}  + \sum_{i=0}^{l-1} \left( \binom{M-1}{i} -\binom{M}{i} \right )(p^*_l)^i (1-p^*_l)^{M-i-1} \le \frac{\delta}{1-p^*_l}, 
\end{eqnarray}
}%
Using  \eqref{eqn_FM_FM_1} two times
(with equality only for $l=1$) and  from \eqref{eqn_prf_info}, 

\vspace{-1mm}
{\footnotesize
\begin{eqnarray*}
\frac{U_k^* - U_1^*}{M} = \og{k}p_k^* - \og{1} p_1^*  =  (C_v + \Gamma) (p^*_k -p^*_1) - \CI \left(  F_{M-1}(k-1;p^*_k) p^*_k- F_{M-1}(0;p^*_1)p^*_1 \right), \hspace{-115mm}\\
&>&  (C_v + \Gamma) (p^*_k -p^*_1) - C_i \left( \frac{\delta p^*_k}{1-p^*_k} -\frac{\delta p^*_1}{1-p^*_1}  \right) = (p^*_k -p^*_1) \left( C_v + \Gamma-\frac{C_i\delta }{(1-p^*_1)(1-p^*_k)} \right),
 \end{eqnarray*}}for any $1 < k  < M$. 
 By lemma \ref{lem_del_zero}, $\frac{\delta }{(1-p^*_1)(1-p^*_k)}  \to 0 $, as $\delta \to 0$, proving part(i), further because $p_k^*>p_1^*$ by the same lemma.  

From \eqref{eqn_prf_info}-\eqref{eqn_perf_opt},  {\small$\Ui^*_1 = M (C_v + \Gamma - \CI\delta^{\frac{M-1}{M}}) (1- \delta^{\frac{1}{M}})$}  and increases to {\small $M(C_v + \Gamma)$} as $\delta \to 0$, while 
{\small $\Ui^*_M = M(C_v +\Gamma - \CI)$}. Thus part(ii) follows using part(i). \eop

\noindent\textbf{Proof of theorem \ref{thm_gM_smallest}:} For part (i), we compare the case of  $\baz=M$ with $\baz=k$ for any fixed $k< M$.  The notations like   $P_g$ etc., correspond to   $\baz=k$. 

Let $\tg:=\max\{C_v - C_i + \frac{c_{\sd,1}+\mvs}{T},0\}$ as in lemma \ref{lem_monotonicity_of_EF}, and $\Bar{\delta}:=\min_{\sbaz}\Np(\tg,\baz)$; by density $f_{\bar{\xi}}$ of  \textbf{A.}2 with support $[0,\infty)$,  we have $\Bar{\delta}>0$. Thus by monotonicity  and theorem \ref{thm_lead_optimal}, for any $\delta<\bar{\delta}$ and any $\baz$, optimal $\og{\sbaz}$ satisfies $\Np(\og{\sbaz},\baz)=\delta$ (and in fact is more than $\tg$).

Again from theorem \ref{thm_lead_optimal}, now at $\baz=M$,  $\og{M}:= \og{M}(\sigma^2,\delta) $ satisfies,
\begin{eqnarray}\label{eqn_appen_gm}
\Np(\og{M},M) = P ( \Gamma_C \ge \og{M} - C_v + C_i) = \delta.
  \end{eqnarray}%
Hence    $\og{M} \uparrow \infty$ as $\delta \downarrow 0$ (similarly the same is true for any $k$). Let $f$ be the density component of $\Gmc$ --  given as $f_{\Gamma}$  in the proof of lemma \ref{lem_monotonicity_of_EF}. As $f$ has support on $[\Bar{\Gamma},\infty)$, it is possible to further reduce $\Bar{\delta}$ (if required) such that, 
\begin{equation}\label{eqn_f_decreasing}
    \gamma \mapsto f(\gamma)  \mbox{ is decreasing   for all } \gamma \ge {\bar g}-C_v \mbox{ where } {\bar g} := \og{M}(\sigma^2,{\bar \delta}).
\end{equation}
    
Now consider $\baz=k<M$. For any $\delta<\Bar{\delta}$, using \eqref{eqn_outcome_vaccination} and \eqref{eqn_appen_gm} we have (as $F_M (k-1, p) \le 1$ and equals 0 at $p=1$):

\vspace{-3mm}
{\small
\begin{equation*}
    \Np(\og{M}+C_i,k) \le 
    E[\indc{P_{\og{M}+C_i}<  1 }]=P( \Gamma_C - \og{M}- C_i > - C_v)
 = \delta.
    \end{equation*}}%
By theorem \ref{thm_lead_optimal}, and monotonicity of  lemma \ref{lem_monotonicity_of_EF},   $\og{k}:= \og{k} (\sigma^2,\delta) \le \og{M} + C_i
$.

Similar logic further provides     $\og{M} + \alpha C_i \le \og{k}$ for some  $\alpha>0$ and any $\delta<\Bar{\delta}$ as explained next. Let $g :=\og{M} + \alpha C_i  $, exact choice of $\alpha$ will be evident at the end of these arguments. 
We briefly represent  $F(k-1; P_g)$ as $F(g, \Gmc)$, as it depends on $C$ only via $\Gmc$. 
Then the  expected probability of non-eradication with $\baz =k$:
\begin{eqnarray}\label{eqn_appen_np}
   \Np(g,k)  &=&   E[F(g, \Gmc)] =
    E[F(g, \Gmc) \indc{P_g=0}] +  E[F(g, \Gmc)\indc{P_g\ne 0}],\nonumber\\
     &=& P (P_g =0)+ E[F(g, \Gmc)\indc{P_g\in (0,1)}],
    \end{eqnarray}
as $F(g, \Gmc)=0$ for all $\Gmc$ whenever $P_g=1$. From \eqref{eqn_outcome_vaccination}, the first term  equals,

\vspace{-3mm}
{\small
\begin{eqnarray}\label{eqn_appen_pg_0}
&&P (P_g =0) =   P ( \Gmc  \ge g - C_v + C_i) =P ( \Gmc \ge \og{M}  - C_v + (1+\alpha)C_i),
\\ 
&&= \delta - P\left(\og{M} - C_v + C_i \le \Gmc \le \og{M} - C_v + (1+\alpha)C_i\right), \mbox{ using \eqref{eqn_appen_gm}.}\hspace{10mm}\nonumber
\end{eqnarray}}
The second term in \eqref{eqn_appen_np} using \eqref{eqn_outcome_vaccination} equals (here $\gamma$ is a realisation of $\Gmc$):

\vspace{-3mm}
{\small\begin{eqnarray*} 
   \int_{g- C_v }^{g- C_v + C_i}  \hspace{-0.7cm} F(g, \gamma) f(\gamma) d\gamma =    \int_{\og{M}- C_v +  \alpha C_i}^{\og{M}- C_v +  C_i} \hspace{-0.7cm}  F(g, \gamma) f(\gamma) d\gamma +    \int_{\og{M}- C_v +   C_i}^{\og{M}- C_v + (1 + \alpha) C_i}\hspace{-0.7cm}   F(g, \gamma) f(\gamma) d\gamma.
\end{eqnarray*}}%
Hence, using the above  and \eqref{eqn_appen_pg_0} in \eqref{eqn_appen_np},

 \vspace{-3mm}
{\small 
    \begin{eqnarray}\label{eqn_appen_integ}
&& \hspace{-7mm}\Np(g,k) =    \delta    -   \int_{\og{M}- C_v +  C_i}^{\og{M}- C_v + (1 + \alpha) C_i} f(\gamma) d\gamma +    \int_{g -C_v}^{g- C_v +  C_i}  F(g, \gamma) f(\gamma) d\gamma, \\   
 &=& 
\delta +  \int_{\og{M}- C_v + \alpha C_i}^{\og{M}- C_v +  C_i} F(g, \gamma) f(\gamma) d\gamma -   \int_{\og{M}- C_v +  C_i}^{\og{M} - C_v + (1 + \alpha) C_i} (1- F(g, \gamma) )f(\gamma) d\gamma. \nonumber
\end{eqnarray}}
Observe from \eqref{eqn_outcome_vaccination}-\eqref{eqn_p_outcome_F} that $P_g = p(g,C)$ depends on $(g, C)$ only via  $\Gmc-g$, and $F(g, \Gmc)$ depends only on $P_g$ -- so $F(g, \Gmc )= F(\Gmc-g)$. 
Then by \eqref{eqn_appen_integ}, 

\vspace{-3mm}
{\small
\begin{eqnarray*}
\Np(g,k) 
     &=&  
\delta +  \int_{\og{M}+ \alpha C_i- C_v }^{\og{M}+  C_i- C_v } \hspace{-4mm} F(\gamma-g) f(\gamma) d\gamma -   \int_{\og{M}+  C_i- C_v }^{\og{M}+ (1 + \alpha) C_i - C_v } \hspace{-4mm}(1-  F(\gamma-g)  )f(\gamma) d\gamma. \nonumber
\end{eqnarray*}}
Define $D:=f(C_i-C_v+ \og{M})$.  By  change of variables $\Tilde{\gamma}=\gamma-\og{M}$,
 $\Np(g,k)$ equals (by choice of $\bar{\delta}$, for any $\delta<\Bar{\delta}$, we have $\og{M}>\Bar{g}$ given in  \eqref{eqn_f_decreasing}),

\vspace{-3mm}
{\small 
\begin{eqnarray*}
  &=& \delta +  \int_{\alpha C_i - C_v }^{\CI- C_v } \hspace{-3mm} F(\tgamma-\alpha \CI) f(\tgamma + \og{M}) d\tgamma -   \int_{\CI- C_v }^{  (1 + \alpha) C_i - C_v} \hspace{-7mm} (1- F(\tgamma-\alpha \CI) )f(\tgamma + \og{M}) d\tgamma,\\
  &\ge& \delta +  D \left(\int_{\alpha C_i - C_v }^{\CI- C_v } \hspace{-2mm}F(\tgamma-\alpha \CI) d\tgamma -   \int_{\CI- C_v }^{  (1 + \alpha) C_i - C_v} \hspace{-2mm} (1- F(\tgamma-\alpha \CI) )d\tgamma,\right).
\end{eqnarray*}}%
Now when $\alpha \to 0$, the second term in the above increases, while the last one decreases to 0, and both  the terms are independent of $\og{M}$ and $\delta$. Thus   there exists an $\alpha>0$, such that $ \Np( \og{M} + \Tilde{\alpha} C_i,\baz)\ge \delta$ for all $\Tilde{\alpha}\le \alpha$ and hence $\og{k}\ge \og{M}+ \alpha C_i$ for any $\delta
<\Bar{\delta}$ (by theorem \ref{thm_lead_optimal} and  monotonicity of lemma \ref{lem_monotonicity_of_EF}).

Now for any $\baz=k$, as $\delta \to 0$,  $\og{k}\to \infty$ and  thus from \eqref{eqn_outcome_vaccination} $P_{\og{k}} \to 1$  a.s. By DCT, $E[P_{\og{k}} ] \to 1$ and so  $\Ui^*_k - M \og{k} \to 0$ for any $k$. Consider an $\epsilon \in(0, M \alpha C_i)$, with $\alpha$ as in above.   Then, there exists a ${\bar \delta}$ (if required less than the current one) such that $\Ui^*_M < M \og{M} + \epsilon < M \og{M} + M \alpha C_i\le M \og{k}$ for all $\delta \le {\bar \delta}$. This proves part (i).

\renewcommand{\Gmcm}{\Gamma^*}

For  part (ii), fix $\baz<M$ and  reduce $\Bar{\delta}$ further if required to ensure $\Np(0,\baz)>\Bar{\delta}$ even for  $\perf{c_{\infty}}$ and consider $\delta<\Bar{\delta}$. We now use special notations to show explicit dependency on $\sigma^2$ -- let  $\Gmcs$, $\Gamma_{\Gmcs}$ respectively represent $C$ and $\Gamma_C$  when $\{\bar{\xi}_i\}$ in \eqref{eqn_ck} have variance $\sigma^2$. With slight abuse of notation, \textit{we let $\sigma^2$ represent the variance of ${C_\sigma}$ itself}. Similarly $\Np(g,\baz; \sigma^2)$  corresponds to  the case with variance $\sigma^2$. When $\sigma^2=0$, $C_\sigma = c^* := E[{C_\sigma}]$ a.s., and  thus $\Np(g,\baz; 0) = F_M(\baz-1,p(g,c^*))$.

Clearly $E[\mid C_\sigma - c^*\mid^2] = \sigma^2$,  and hence  when $\sigma^2 \to 0$, $C_\sigma$ converges to $c^*$  in  distribution. 
Further as before, $c \to F_M(\baz-1,p(g,c))$ is a  continuous and bounded function. Thus,
  $\Np(g,\baz; \sigma^2) \to F_M(\baz-1,p(g,c^*)) $  by weak convergence for any   $g$ as $\sigma^2 \to 0$. 
  
Let $\Gmcm:= C_\infty$ and observe $\Gmcm= \Gamma_{T-1}(c^*)$ by direct substitution using lemma \ref{lem_vac_val}. Let $g^*$ denote $\og{\sbaz}$ and $p^*$ denote $ p^*_{\sbaz}$ of \eqref{eqn_prf_info} for $\perf{c_\infty}$.   Since $\delta\in(0,1)$, by \eqref{eqn_prf_info} we have $p^* \in (0,1)$ which implies $0< C_v + \Gmcm - g^* < \CI$. Hence from \eqref{eqn_p_outcome_F} $
  p(g^*,c^*) = p^* $. Thus when $\sigma^2 \to 0$, again from  \eqref{eqn_prf_info}, we have:
 $$\Np(g^*,\baz; \sigma^2) \to  F_M(\baz-1,p(g^*,c^*))=F_M(\baz-1,p^*)=\delta . $$
Let  $g_{\sigma}$ be the optimal incentive scheme of the leader when  variance is $\sigma^2$. By theorem \ref{thm_lead_optimal},  $\Np(g_\sigma,\baz; \sigma^2) = \delta$ for all $\sigma^2>0$ and so, 
\begin{eqnarray*}
\mid \Np(g^*,\baz; \sigma^2) - \Np(g_\sigma,\baz; \sigma^2) \mid = \mid \Np(g^*,\baz; \sigma^2) - \delta \mid \to 0, \mbox{ as $\sigma^2\to 0$. }
\end{eqnarray*}%
Hence $g_\sigma \to g^*$, as  by lemma \ref{lem_monotonicity_of_EF}, $g\mapsto \Np(g ,\baz; \sigma^2)$ is continuous and strict monotone on $(\Bar{g},\infty)$ and thus has continuous inverse (recall, $\bar{g}\ge\tg$ by the choice of $\Bar{\delta}$). Now, it {\it suffices to show that  $ E[p(g_{\sigma},\Gmcs)] \to p^*$,} as then 
  $$\Ui_{\sbaz}^*(\sigma^2,\delta)=  M g_{\sigma} E[P_{g_{\sigma}}] =   M g_{\sigma} E[p(g_{\sigma},\Gmcs)] \to M g^* p^* = U^*_{\sbaz} \mbox{ of }\perf{c_\infty}.$$
  Towards proving the above, define $A_\epsilon:=\{\omega: \mid C_{\sigma}(\omega)-c^*
\mid \le \epsilon\} $ which by continuity of $p$ is a subset of $ \{\omega: \mid p(g_{\sigma},C_{\sigma}(\omega))-p(g_{\sigma},c^*)
\mid \le \hat{\epsilon}\}$ for some appropriate $ \hat{\epsilon}$ and  consider:

\vspace{-3mm}
{\small\begin{eqnarray}\label{eqn_app_convg}
   \mid E[p(g_{\sigma},\Gmcs)] - p^* \mid &  \le & \mid  E[p(g_{\sigma},\Gmcs)] - p(g_{\sigma},c^*) \mid + \mid  p(g_{\sigma},c^*) - p^* \mid  \\
   &  \le & \mid  E[p(g_{\sigma},\Gmcs) - p(g_{\sigma},c^*); A_\epsilon ]\mid + 2 P(A_\epsilon^c) +  \mid  p(g_{\sigma},c^*) - p^* \mid \nonumber \\
      &  \le &  \hat{\epsilon} + 2 P(A_\epsilon^c) +  \mid  p(g_{\sigma},c^*) - p^* \mid \to 0,  \nonumber 
\end{eqnarray}}because: 
by continuity $p(g_{\sigma},c^*) \to  p (g^*, c^*) =p^* $, $P(A_\epsilon^c) \to 0$ by weak convergence of $C_\sigma$ to  constant $c^*$ and finally $ \hat{\epsilon} \to 0 $ as $\epsilon \to 0.$ 

To complete part (ii), now consider  $\baz=M$ and let $\og{M,\epsilon}$ be as in \eqref{eqn_perf_opt} for $\perf{c_\infty}$. Thus by weak convergence,
{\small $$ P ( \Gamma_{C_\sigma} \ge \og{M,\epsilon} - C_v + C_i) {=}  P ( \Gamma_{C_\sigma}  \ge \Gamma^* + \epsilon) \downarrow 0 \mbox{ as }  \sigma^2 \to 0.$$}%
 Fix $\delta < \bar{\delta}$, and choose $\bar{\sigma}^2$ such that $ P ( \Gamma_{C_\sigma} \ge \og{M,\epsilon} - C_v + C_i) <\delta$ for all $\sigma^2 \le \bar{\sigma}^2$. 
 While 
 by \eqref{eqn_appen_gm},  at  optimizer 
 $P(\Gamma_{\Gmcs} \ge  \og{M}(\sigma^2,\delta) - C_v + C_i) = \delta$, and hence by strict monotonicity given by \textbf{A}.2, $ \og{M}(\sigma^2,\delta) < \og{M,\epsilon}$.  Also, 
$$U_M^*(\sigma^2,\delta) \le  M \og{M}(\sigma^2,\delta) < M \og{M,\epsilon} = U_M^{*} + M \epsilon, $$%
where $U^*_M$ corresponds to $\perf{c_{\infty}}$. \eop

\begin{lemma}\label{lem_min_psi}
 \begin{enumerate}[(i)] 
\item The vaccinated proportion $\ps{e}(\bnu)$ of population at eradicating ESSS, under any admissible VA rate strategy $\bnu$ is lower bounded, $\ps{e}(\bnu)>\nvdf$.
\item  For any $\epsilon > 0$, there exists $\teps>0$   such that  $\ps{e}(\bnu) < \nvdf  + {\epsilon}$ if  $\nu_b \in (0, b\rho \nvdf]$ and $\nu_e =  b \rho - \frac{\nu_b }{\nvdf} + \teps $. Further, as $\epsilon \to 0$ we have $ \teps \to 0$.
\end{enumerate}

\ignore{For any $\epsilon >0$ and $\nu_b \in (0, b\rho \nvdf]$ , there exists an admissible $\bnu^\epsilon = (\nu_b, \nu_e^\epsilon$ such that 
 $\nvdf  < \ps{e}(\bnu^\epsilon) < \nvdf  + {\epsilon}$ and $\nu_e = ( b \rho - \frac{\nu_b }{\nvdf} + \teps )$. Further, as $\epsilon \to 0$ we have $ \teps \to 0$.}


\end{lemma}
\noindent\textbf{Proof:} 
First observe that  the function  $\ps{e} (\bnu) $ defined in  table \ref{table_ess_candidate} can be extended to any $\bnu = (\nu_b, \nu_e) \in [0, \infty)^2$; this function provides   attractors when $\bnu$ is admissible (see \eqref{eqn_eradication_condition}). We first establish some properties of this function.

Firstly $\nu_b \mapsto \ps{e}(\bnu)$ is increasing, as the corresponding partial derivative:

\vspace{-3mm}
{\small
\begin{eqnarray}\label{eqn_appen_psi_nu_b}
  \pdv{ \ps{e}(\bnu)}{\nu_b} =  \frac{1}{2\nu_e} \left( \frac{(b +\nu_b + \nu_e)}{\sqrt{(b +\nu_b-\nu_e)^2 + 4 \nu_e \nu_b}} -1 \right)\indc{\nu_e>0} + \frac{b\indc{\nu_e=0}}{(\nu_b + b)^2} >0.
\end{eqnarray}}%
Now consider the partial derivative of $\ps{e}$ with respect to $\nu_e$ for  $\nu_e>0$,  

\vspace{-3mm}
{\small
\begin{eqnarray*}
     \pdv{ \ps{e}(\bnu)}{\nu_e} &=& \frac{ \nu_e 
     \left( 1 + \frac{\nu_b + \nu_e - b}{\sqrt{ (b+\nu_b-\nu_e)^2 + 4 \nu_b \nu_e}} \right) +(b+\nu_b - \nu_e) - \sqrt{ (b+\nu_b -\nu_e)^2 + 4 \nu_b \nu_e}}{2\nu_e^2}.
\end{eqnarray*}}%
The derivative of numerator in above with respect to $\nu_e$ equals,

\vspace{-3mm}
{\small
\begin{eqnarray*}
\left(\nu_e \frac{(b +\nu_b-\nu_e)^2 + 4 \nu_e \nu_b- (\nu_b +\nu_e-b)^2 }{((b +\nu_b-\nu_e)^2 + 4 \nu_e \nu_b)^{\frac{3}{2}}}\right) >0 \mbox{ for all } \nu_e>0,
\end{eqnarray*}}%
as  
$(b + \nu_b - \nu_e)^2 + 4 \nu_b \nu_e = (\nu_b + \nu_e - b)^2 + 4 b \nu_b $, and hence $\pdv{ \ps{e}(\bnu)}{\nu_e}>0$.

Further $\ps{e}$ is a continuous function of $\bnu$ on $[0, \infty)^2$, as clearly for any $\nu_b \ge 0 $ 
%
%
\begin{equation}\label{Eqn_nue_zero}
    \psi_e (\nu_b, 0)= \frac{\nu_b}{\nu_b + b} = \lim_{\nu_e\downarrow 0}  \ps{e}(\nu_b,\nu_e) \   \mbox{ (using L'Hopital rule). }
\end{equation}%
%
%
Now consider an admissible $\bnu$ such that $\nu_b> b\rho\nvdf$. Then from \eqref{eqn_eradication_condition},     we have $\nu_e \ge 0$ and hence using the monotonicity in $\nu_e$ and \eqref{Eqn_nue_zero}:

\vspace{-3mm}
{\small $$\ps{e}(\bnu) \ge \ps{e}(\nu_b,0) = \frac{\nu_b}{\nu_b + b} >  \frac{b \rho \theta^*}{b \rho \theta^* + b} = \nvdf.$$}%
If $\nu_b = b \rho \nvdf $,   $\bnu$ is admissible only if  $\nu_e > 0$, again  using the monotonicity:
$$ \ps{e}(\bnu)> \ps{e} \left (\nu_b, \frac{\nu_e}{2} \right )  \ge  \ps{e}(\nu_b,0) = \nvdf.
%
$$
Finally when $\nu_b < b \rho \nvdf$, $\bnu$ is admissible only if $\nu_e >  b \rho - \frac{\nu_b}{\nvdf}$.
Then\footnote{%
The square-root term in $\ps{e}(\nu_b,\nu_e')$ at $\nu_e'= b\rho -\frac{\nu_b}{\nvdf} $ equals:\\
$$(b +\nu_b+\nu_e')^2 - 4 b \nu_e' =(b\rho + \frac{\nu_e'}{\rho})^2  - 4 b \nu_e' =(b \rho - \frac{\nu_e'}{\rho})^2. 
   $$%
   So $ \ps{e}(\nu_b,\nu_e') =   \frac{-b -\nu_b + \nu_e'+b \rho - \frac{\nu_e'}{\rho}}{2 \nu_e' }=\nvdf.$}%
 $ \ps{e}(\bnu) > \ps{e}\left(\nu_b,   b \rho - \frac{\nu_b}{\nvdf}\right)=\nvdf$. Last
part  follows from the continuity and  strict monotonicity of $\ps{e}$ and that $\ps{e} \left(\nu_b,   b \rho - \frac{\nu_b}{\nvdf}\right) = \nvdf$ when $\nu_b \in [0, b\rho\nvdf]$.
\eop

\noindent\textbf{Proof of theorem \ref{thm_iff_conditionn_for_optimality}:}
Towards the first part, say $\frac{c_{v_2}}{\nvdf}<\Bar{c}_{v_2}.$ Then, 

\vspace{-3mm}
{\small$$L_k = \min\left\{ -  \frac{c_{v_2}}{\nvdf} \frac{M-k}{M},\frac{r+b}{r+2b } c_i  -  \bar{c}_{v_2}\frac{M-k}{M}\right\}, \mbox{ and let } O_k:= c_{v_1}-c_f(k) - L_k.$$}%
Note $L_k$ is strictly increasing and $O_k$ is strictly decreasing; also observe $c_{v_1}-c_f(0) > L_0 $  (as $c_f(0) = 0$) and  $c_{v_1}-c_f(M)\le L_M$ by \textbf{A.}1.   
Thus $O_k > 0$ at $k = 0$ and $O_k< 0$ at $k=M$, which establishes the existence and uniqueness of $k$ satisfying \eqref{eqn_k_eps_opt2}.

Now we construct  an $\epsilon$-optimal $\bnu^\epsilon$ such that $h_i(k;\bnu^\epsilon) \le 0$, $h_v(k;\bnu^\epsilon)<0$, and either $h_i(k-1;\bnu^\epsilon)>0$ or $h_v(k-1;\bnu^\epsilon) \ge 0$; by theorem \ref{thm_eradication_es}, $\baz(\bnu^\epsilon)=k$. Towards this, let $\bnu^\epsilon = (\nu_b^\epsilon,\nu_e^\epsilon) := (b\rho \nvdf -\epsilon_1, \ \ b \rho - \frac{\nu^\epsilon_b}{\nvdf} + \epsilon_2)$, where $(\epsilon_1, \epsilon_2)$ are chosen such that $\nvdf<\ps{e}(\bnu^\epsilon)<\nvdf+\epsilon$, possible by lemma \ref{lem_min_psi}. Hence from \eqref{eqn_h_func}, \eqref{eqn_k_eps_opt2}:

\vspace{-3mm}
{\small\begin{eqnarray*}
 h_i(k;\bnu^\epsilon)&=&c_{v_1} + \left(1-\frac{k}{M}\right)\bar{c}_{v_2} - \frac{r+b}{r+2b +b (\frac{\nu_b^{\epsilon}}{b \rho \nvdf}-1)} c_i - c_f(k) <0,\\  
  h_v(k;\bnu^\epsilon)&=&c_{v_1} + \left(1-\frac{k}{M}\right) \frac{{c}_{v_2}}{\ps{e}(\bnu^\epsilon)} - c_f(k) <0.  
\end{eqnarray*}}%
Further $h_v(k-1;\bnu^\epsilon)>0$ if {\tiny$-\frac{c_{v_2}}{\nvdf} \frac{M-k+1}{M}<\frac{r+b}{r+2b}c_i- \bar{c}_{v_2}\frac{M-k+1}{M}$}, else  $h_i(k-1;\bnu^\epsilon)> 0$,  by choosing $(\epsilon_1, \epsilon_2)$ further small (if required), possible again by lemma~\ref{lem_min_psi}.
Hence, $\baz(\bnu^\epsilon)=k$ by theorem \ref{thm_eradication_es}.

Now say $\frac{c_{v_2}}{\nvdf}\ge \Bar{c}_{v_2}.$ Then,
$L_k =  -  \bar{c}_{v_2}\frac{M-k}{M}$, 
 and let  $O_k:= c_{v_1}-c_f(k) - L_k$. 
 The arguments for existence and uniqueness of $k$ satisfying \eqref{eqn_k_eps_opt} follow in the exact similar manner as above. 
 From \eqref{eqn_k_eps_opt} and \eqref{eqn_h_func}, $h_i(k;\bnu) < h_v(k;\bnu)<0$ and $h_v(k-1;\bnu) \ge 0$ for any $\bnu$. Hence, the same is true for an  $\epsilon$-optimal $\bnu^\epsilon$ as well, which proves part (i).

Towards part (ii), first consider the case when $\frac{c_{v_2}}{\nvdf}< \Bar{c}_{v_2}$. Say $(\bg,\bnu)$ is incentive optimal policy, then $\baz(\bnu)=M$. From \eqref{eqn_h_func}, theorem \ref{thm_eradication_es}  and lemma \ref{lem_min_psi}, %

\vspace{-3mm}
{\small
\begin{eqnarray*}
 \mbox{\normalsize either } c_{v_1} - c_f(M-1) \ge  -\frac{{c}_{v_2}}{M \ps{e}(\bnu)} {>}-\frac{\bar{c}_{v_2}}{M} \ \  \left (\mbox{\normalsize for }h_v(M-1;\bnu)\ge 0\right ),\\
 \mbox{\normalsize or }   c_{v_1} - c_f(M-1) 
 > -\frac{\bar{c}_{v_2}}{M} + \frac{\lambda \nvdf}{\lambda \nvdf+\nu_b} c_{i}>-\frac{\bar{c}_{v_2}}{M} \  \  \left (\mbox{\normalsize for } h_i(M-1;\bnu)>0\right ),
\end{eqnarray*}}%
which gives \eqref{eqn_inc_opt_cond}. Now say \eqref{eqn_inc_opt_cond} holds. If $c_{v_1} + \frac{\Bar{c}_{v_2}}{M} - c_f(M-1) - c_i > 0$, then from \eqref{eqn_h_func} $h_i(M-1;\bnu)>0$  implying $\baz(\bnu)=M$ for any $\bnu$ by theorem \ref{thm_eradication_es}. If not, choose,

\vspace{-3mm}
{\small $$\nu_b > \frac{ \lambda \nvdf c_i}{c_{v_1} + \frac{\Bar{c}_{v_2}}{M} - c_f(M-1)} - \lambda \nvdf,$$}%
and a $\nu_e$ such 
that $\bnu =(\nu_b, \nu_e) $ is admissible. 
Then  $h_i(M-1;\bnu)>0$ implying $\baz(\bnu)=M$ by theorem \ref{thm_eradication_es}. Further let $\bg$ be the optimizer of \eqref{eqn_leader_opt} for $\baz=M$. Then clearly $(\bg,\bnu)$ is  incentive optimal.

Now consider the case when $\frac{c_{v_2}}{\nvdf}\ge \Bar{c}_{v_2}$. Say $(\bg,\bnu)$ is incentive optimal policy, then $\baz(\bnu)=M$. From \eqref{eqn_h_func}, theorem \ref{thm_eradication_es}  and lemma \ref{lem_min_psi}, %

\vspace{-3mm}
{\small
\begin{eqnarray*}
 \mbox{\normalsize either } c_{v_1} - c_f(M-1) \ge  -\frac{\bar{c}_{v_2}}{M} \ \  \left (\mbox{\normalsize for }h_v(M-1;\bnu)\ge 0\right ),\\
 \mbox{\normalsize or }   c_{v_1} - c_f(M-1) 
 > -\frac{\bar{c}_{v_2}}{M} + \frac{\lambda \nvdf}{\lambda \nvdf+\nu_b} c_{i}>-\frac{\bar{c}_{v_2}}{M} \  \  \left (\mbox{\normalsize for } h_i(M-1;\bnu)>0\right ),
\end{eqnarray*}}%
which gives \eqref{eqn_inc_opt_cond}. Now say \eqref{eqn_inc_opt_cond} holds, and let $\bg$ be the optimizer of \eqref{eqn_leader_opt} for $\baz=M$. Then clearly $(\bg,\bnu)$ is  incentive optimal for any $\bnu$, since $h_v(M-1;\bnu)\ge 0$ by \eqref{eqn_h_func}. Hence part (ii) is true.
\eop

\end{appendices}

\end{document}